\documentclass[a4paper,11pt]{article}

\usepackage{mathrsfs} 
\usepackage{authblk}
\usepackage{amsthm}
\usepackage{hyperref}
\usepackage{amssymb}
\usepackage{latexsym}
\usepackage{amsmath}
\usepackage{amsfonts}
\usepackage{mathabx}
\usepackage{hyperref}
\usepackage{mathtools}
\hypersetup{
    colorlinks=true,
    linkcolor=blue,
    filecolor=magenta,      
    urlcolor=cyan,
}
\usepackage[dvips]{graphicx}
\usepackage[utf8]{inputenc}

\usepackage[LGR,T1]{fontenc}%
\usepackage[francais,english]{babel}
\usepackage{url}
\usepackage{enumerate}
\usepackage{hyperref}
\usepackage{color}

\usepackage{bbm}

\usepackage{fancyhdr}
\usepackage{float}
\usepackage{graphicx}
\usepackage{cleveref}
\usepackage{diagbox}
\usepackage{multicol}
\usepackage{multirow}
\newcommand{\tabincell}[2]{\begin{tabular}{@{}#1@{}}#2\end{tabular}}

\usepackage{exscale}
\usepackage{relsize}
\usepackage{cancel}

\usepackage{algorithm}
\usepackage{algorithmicx}
\usepackage{algpseudocode}

\newtheorem{Theorem}{Theorem}[part]
\newtheorem{Definition}{Definition}[part]
\newtheorem{Proposition}{Proposition}[part]

\newtheorem{Assumption}{Assumption}[part]
\newtheorem{Lemma}{Lemma}[part]
\newtheorem{Corollary}{Corollary}[part]
\newtheorem{Remark}{Remark}[part]

\makeatletter \@addtoreset{equation}{section}

\@addtoreset{Definition}{section}

\@addtoreset{Theorem}{section}

\@addtoreset{Proposition}{section}

\@addtoreset{Property}{section}

\@addtoreset{Assumption}{section}

\@addtoreset{Corollary}{section}

\@addtoreset{Lemma}{section}

\@addtoreset{Remark}{section}

\@addtoreset{Example}{section}

\usepackage{amsmath}

\DeclareMathOperator*{\argmin}{arg\,min}

\DeclarePairedDelimiter\ceil{\lceil}{\rceil}
\DeclarePairedDelimiter\floor{\lfloor}{\rfloor}

\usepackage{xr}
\newcommand{\cA}{\mathcal{A}}
\newcommand{\cB}{\mathcal{B}}
\newcommand{\cC}{\mathcal{C}}

\newcommand{\cE}{\mathcal{E}}
\newcommand{\cF}{\mathcal{F}}

\newcommand{\cH}{\mathcal{H}}

\newcommand{\cL}{\mathcal{L}}
\newcommand{\cM}{\mathcal{M}}

\newcommand{\cO}{\mathcal{O}}

\newcommand{\cS}{\mathcal{S}}

\newcommand{\cU}{\mathcal{U}}

\newcommand{\cW}{\mathcal{W}}
\newcommand{\cX}{\mathcal{X}}
\newcommand{\cY}{\mathcal{Y}}
\newcommand{\cZ}{\mathcal{Z}}

\renewcommand{\P}{\mathbb{P}}

\newcommand{\R}{\mathbb{R}}

\newcommand{\Z}{\mathbb{Z}}

\def \proof{{\noindent \bf Proof. }}
\def \eproof{\hbox{ }\hfill$\Box$}

\newcommand{\ud}{\, \mathrm{d}}

\newcommand{\1}{{\bf 1}}
    
\newcommand{\set}[1]
    {\ensuremath{\{ #1 \}}}
    
\newcommand{\HP}[1] 
    {\ensuremath{\mathscr{H}^{#1}}}

\newcommand{\esp}[1]{\ensuremath{\mathbb{E} \! \left[#1\right] }}

\newcommand{\widehatbis}[1]{\ensuremath{#1}}

\title{A learning scheme by sparse grids and Picard approximations for semilinear parabolic PDEs}

\author{J.-F. Chassagneux \thanks{Universit\'e de Paris, Laboratoire de Probabilit\'es, Statistiques et Mod\'elisation, F-75013 Paris, France. Email: chassagneux@lpsm.paris.} \, ,  J. Chen  \thanks{Universit\'e de Paris, Laboratoire de Probabilit\'es, Statistiques et Mod\'elisation, F-75013 Paris, France. Part of the work was done during a visit at NUSRI, Suzhou, China, whose hospitality is greatly appreciated.  Email: juchen@lpsm.paris} \,, N. Frikha \thanks{Universit\'e de Paris, Laboratoire de Probabilit\'es, Statistiques et Mod\'elisation, F-75013 Paris, France. Email: frikha@lpsm.paris} \, , C. Zhou  \thanks{Department of Mathematics, City University of Hong Kong,  Kowloon Tong, Hong Kong. On leave from Department of Mathematics, National University of Singapore. Research supported by NSFC Grant No. 11871364 as well as Singapore MOE (Ministry of Educations) AcRF Grants R-146-000-271-112 and R-146-000-284-114. Email: chaozhou@cityu.edu.hk}.}

\begin{document}
\maketitle

\begin{abstract}

Relying on the classical connection between Backward Stochastic Differential Equations (BSDEs) and non-linear parabolic partial differential equations (PDEs), we propose a new probabilistic learning scheme for solving high-dimensional semilinear parabolic PDEs. This scheme is inspired by the approach coming from machine learning and developed using deep neural networks in Han and al. \cite{Han8505}. Our algorithm is based on a Picard iteration scheme in which a sequence of linear-quadratic optimisation problem is solved by means of stochastic gradient descent (SGD) algorithm. In the framework of a linear specification of the approximation space, we manage to prove a convergence result for our scheme, under some smallness condition.
In practice, in order to be able to treat high-dimensional examples, we employ sparse grid approximation spaces. In the case of periodic coefficients and using pre-wavelet basis functions, we obtain an upper bound on the global complexity of our method. It shows in particular that the curse of dimensionality is tamed in the sense that in order to achieve a root mean squared error of order $\varepsilon$, for a prescribed precision $\varepsilon$, the complexity of the Picard algorithm grows polynomially in $\varepsilon^{-1}$ up to some logarithmic factor $|\log(\varepsilon)|$ which grows linearly with respect to the PDE dimension. Various numerical results are presented to validate the performance of our method and to compare them with some recent machine learning schemes proposed in Han and al. \cite{weinan2017deep} and Hur\'e and al. \cite{hure2020deep}.

\end{abstract}


\section{Introduction}
In the present work, we are interested in the numerical approximation in high dimension of the solution to the semilinear parabolic PDE
\begin{equation}
\label{eq pde}
\left\{
\begin{array}{ll}
\partial_t u(t, x) + \mathcal{L} u(t, x) + f(u(t,x), \sigma^{\top}(x)\nabla_x u(t, x))  =0, &   (t, x) \in [0,T) \times \mathbb{R}^d, \\
 u(T, x) = g(x), \quad x \in \mathbb{R}^d  & 
\end{array}
\right.
\end{equation}

\noindent where $f: \mathbb{R} \times \mathbb{R}^d \rightarrow \mathbb{R}$, $g: \mathbb{R}^d \rightarrow \mathbb{R}$ are measurable functions and $\mathcal{L}$ is the infinitesimal generator of the forward diffusion process with dynamics 
\begin{equation}
\ud \cX_t  = b(\cX_t) \ud t +    \sigma(\cX_t) \ud \cW_t \, \label{eq forward sde}
\end{equation}

\noindent and defined, for a smooth function $\varphi$, by
\begin{equation}
\label{infinitesimal:operator:diffusion}
\mathcal{L} \varphi(t, x) :=  b(x)\cdot \nabla_x \varphi(t, x) +\frac12 \mathrm{Tr}[(\sigma\sigma^{\top})(x)\nabla^2_x \varphi(t, x)].
\end{equation}
\noindent Here, $\cW$ is a $d$-dimensional brownian motion defined on a complete probability space $(\mathfrak{O}, \cA,\P)$, $b : \R^d \rightarrow \R^d$ and $\sigma: \R^d \rightarrow \cM_d$ are measurable functions, $\cM_d$ being the set of $d\times d$ matrix. 
The initial condition $\cX_0$ is a square integrable random variable independent from the Brownian Motion $\cW$.
We denote by $(\mathcal{F}_{t})_{0\leq t \leq T}$ the  filtration generated by $\cW$ and $\cX_0$, augmented with $\P$ null sets. 
 
Developing efficient algorithms for the numerical approximation of high-dimensional non-linear PDEs is a challenging task that has attracted considerable attention from the research community in the last two decades.
We can quote various approaches (limiting to the "stochastic" ones) that have proven to be efficient in a high dimensional setting: branching methods, see e.g. \cite{henry-labordere2019}, machine learning methods (especially using deep neural networks), see e.g. \cite{Han8505}, and full history recursive multilevel Picard method (abbreviated MLP in the literature) see e.g. \cite{ewhutzenthaler2016multilevel}. This is a very active field of research, we refer to the recent survey papers \cite{weinan2020algorithms,beck2020overview} for more references and an overview of the numerical and theoretical results available.
We focus now more on one stream of research which uses the celebrated link between semilinear parabolic PDEs of the form \eqref{eq pde} and BSDEs. This connection, initiated in \cite{Pardoux:Peng:90}, reads as follows: denoting by $u$ a classical solution to \eqref{eq pde}, $(u(t, \cX_t), \sigma^{\top}(\cX_t) \nabla_x u(t, \cX_t))  = (\cY_t, \cZ_t)$ where the pair $(\cY, \cZ)$ is the $\mathbb{R}\times \mathbb{R}^d$-valued and $(\mathcal{F}_{t})$-adapted process solution to the BSDE with dynamics
\begin{align}
\cY_t & = g(\cX_T) + \int_t^Tf(\cY_s,\cZ_s) \ud s - \int_t^T \cZ_s \cdot \ud \cW_s\,,\;0\le t \le T\,,
\label{eq backward sde}
\end{align} 
so that, the original problem boils down to the numerical approximation of the above stochastic system. 
Various strategies have been used to numerically approximate the stochastic system $(\cX_t, \cY_t, \cZ_{t})_{t \in [0,T]}$. The most studied one is based on a time discretization of \eqref{eq backward sde} leading to a backward programming algorithm to approximate $(Y,Z)$, as exposed in e.g \cite{BOUCHARD2004175, zhang2004numerical} (see the references therein for early works). This involves computing a sequence of conditional expectations and various methods have been developed: Malliavin calculus based methods \cite{BOUCHARD2004175, CRISAN, hu2011}, optimal quantization methods \cite{BALLY20031,bally2003,PAGES2018847}, cubature methods \cite{crisan2012solving,crisan2014second,chassagneux2020cubature} and (linear) regression methods, see among others \cite{gobet2005regression, gobet2016approximation, gobet2016stratified}. It is acknowledged that such approaches will be feasible for problems up to dimension 10. This limitation is a manifestation of the so-called ``curse of dimensionality''. Recently, non-linear regression methods using deep neural networks were succesfuly combined with this approach and proved to be capable of tackling problems in high dimension \cite{hure2020deep}. However, other strategies have been introduced in the last five years or so to approximate \eqref{eq backward sde} trying to adopt a "forward point of view". Relying on Wiener chaos expansion and Picard iteration,  \cite{briand2014simulation, geiss2016simulation} introduced a method  that notably works in non-Markovian setting  but is still impacted by the curse of dimension. A key step forward has been realized by the so called \emph{deep BSDEs solver} introduced in \cite{Han8505}. Interpreting the resolution of a BSDEs as an optimisation problem, it relies on the expressivity of deep neural network and well established SGD algorithms to show great performance in practice. 
More precisely, in this approach, the $\cY$-process is now interpreted as a forward SDE controlled by the $\cZ$-process. Then, an Euler-Maruyama approximation scheme is derived in which the derivative of the solution $u$ appearing in the non-linear function $f$ (through the $\cZ$-process) is approximated by a multi-layer neural network. The optimal weights are then computed by minimizing the mean-squared error between the value of the approximation scheme at time $T$ and a good approximation of the target $g(\cX_T)$ using stochastic gradient descent algorithms. Again, this kind of deep learning technique seems to be very efficient to numerically approximate the solution to semi-linear parabolic PDEs in practice. However a complete theory concerning its theoretical performance is still not achieved \cite{beck2020overview}. One important observation is that, due to highly non-linear specification, the optimisation problem that has to be solved in practice, has no convexity property. The numerical procedure designed can only converge to local minima, whose properties (with respect to the approximation question) are still not completely understood. 

Inspired by this new forward approach, we introduce here an algorithm which is shown to converge to a global minimum. This, of course, comes with a price. First, we move from the deep neural networks approximation space to a more classical linear specification of the approximation space. However, due the non-linearity in the BSDE driver, the global optimisation problem to be solved is still non-convex. To circumvent this issue, we employ a \emph{Picard iteration procedure}. The overall procedure becomes then a sequence of linear-quadratic optimisation problems which are solved by a SGD algorithm. {Our first main result is a control of the global error between the \emph{implemented algorithm} and the solution to the BSDE which notably shows the convergence of the method under some smallness conditions, see Theorem \ref{th main conv result}. In particular, contrary to \cite{Han8505, han2020convergence} or \cite{hure2020deep}, our result takes into account the error induced by the SGD algorithm. In our numerical experiments, we rely on sparse grid approximation spaces which are known to be well-suited to deal with high-dimensional problems. Under the framework of periodic coefficients, we establish  as our second main result, an upper bound on the global complexity for our \emph{implemented algorithm}, see Theorem \ref{thm full complexity}. We notably prove that the curse of dimensionality is tamed in the sense that the complexity is of order $\varepsilon^{-p} |\log(\varepsilon)|^{q(d)}$, where $p$ is a constant which does not depend on the PDE dimension and $d\mapsto q(d)$ is an affine function. We also demonstrate numerically the efficiency of our methods in high dimensional setting.}      

The rest of the paper is organized as follows. In \Cref{Description of the numerical methods}, we first recall the \emph{deep BSDEs solver} of \cite{Han8505} but adapted to our framework. Namely, we use  a linear specification of the approximation space together with SGD algorithms. For sake of clarity, we denote this method: the \emph{direct algorithm}. Then, we introduce our new numerical method : the \emph{Picard algorithm}. We present our main assumptions on the coefficients and state our main convergence results. In \Cref{se main sparse grids}, we use sparse grid approximation with the \emph{direct and Picard algorithms}, using two types of basis functions: pre-wavelet \cite{bohn2018convergence} and modified hat function \cite{frommert2010efficient}. We discuss their numerical performances in practice through various test examples. We also compare them with some deep learning techniques \cite{Han8505,hure2020deep}. We also state our main theoretical complexity result. \Cref{study:optimization:problem} is devoted to the theoretical analysis required to establish our main theorems: all the proofs are contained in this section. Finally, we give a complete list of the algorithm parameters that have been used to obtain the numerical results in Appendix.

\paragraph{Notation:}
Elements of $\R^q$ are seen as column vectors. For $x\in\R^q$, $x_i$ is the $i$th component and $|x|$ corresponds to its Euclidian norm, $x \cdot y$ denotes the scalar product of $x$ and $y \in \R^q$. $\cM_q$ is the set of $q \times q$ real matrices. We denote by $\mathbf{e}^\ell$ the $\ell$th vector of the standard basis of $\R^q$. The vector $(1,\dots,1)^\top$ is denoted $\mathbf{1}$, $\mathrm{I}_d$ is the $d \times d$ identity matrix. We use the bold face notations $\bold{l} \in \mathbb{N}^{d}$ for multidimensional indices $\bold{l}:=(l_1, \cdots, l_d)$ with (index) norms denoted by $|\bold{l}|_{p}:=(\sum_{i=1}^{d} l_i^{p})^{1/p}$ and $|\bold{l}|_{\infty} := \max_{1\leq i \leq d} | l_i|$. For later use, for a positive integer $k$, we introduce the set $\mathbf{J}^\infty_{d, k}$ of multidimensional indices $\bold{l} \in \mathbb{N}^{d}$ satisfying $|\bold{l}|_{\infty} \leq k$.
For a finite set $A$, we denote by $|A|$ its cardinality. 

For a function $f:\R^d \rightarrow \R$, we denote by $\partial_{x_l} f$ the partial derivative function with respect to $x_l$, $\nabla f$ denotes the gradient function of $f$, valued in $\R^d$. We also use $\nabla^2 f = (\partial^2_{x_i, x_j} f)_{1\leq i, j \leq d}$ to denote the Hessian matrix of $f$, valued in $\cM_d$. For a sufficiently smooth real-valued function $f$ defined in $\mathbb{R}^d$, we let $D^{\bold{l}} f = \partial^{l_1}_{x_1}\cdots \partial^{l_d}_{x_d} f$ denote the differentiation operator with respect to the multi-index $\bold{l} \in \mathbb{N}^d$. For a fixed positive integer $k$ and a function $f$ defined on an open domain $ \mathcal{U} \subset \mathbb{R}^d$, we define its Sobolev norm of mixed smoothness
\begin{equation}
\label{H2:Sobolev:norm:mixed:smoothness}
\|f\|_{H^{k}_{mix}(\mathcal{U})} := \Big(\sum_{\bold{l} \in \mathbf{J}^\infty_{d, k} } \|D^{\bold{l}}f\|^2_{L_2(\mathcal{U})} \Big)^{\frac12}
\end{equation}

\noindent where the derivative $D^{\bold{l}}f$ in the above formula has to be understood in the weak sense and for a map $g:\mathcal{U} \mapsto \mathbb{R}$, $\|g\|^2_{L^2(\mathcal{U})} := \int_{\mathcal{U}} |g(x)|^2 dx$. The Sobolev space of mixed smoothness $H^{k}_{mix}(\mathcal{U})$ is then defined by 
\begin{equation}\label{Sobolev:space}
H^{k}_{mix}(\mathcal{U}) := \left\{ f \in L_2(\mathcal{U}) : \|f\|_{H^{k}_{mix}(\mathcal{U})} <\infty\right\}.
\end{equation}

For a positive integer $q$, the set $\cH^2_q$  is the set of progressively measurable processes $V$ defined on the probability space $(\mathfrak{O}, \cA,\P)$ with values in $ \mathbb{R}^q$ and satisfying $\esp{ \int_0^T |V_t|^2 \ud t}<+\infty$. The set $\cS^2_q$  is the set of adapted càdlàg processes $U$ defined on the probability space $(\mathfrak{O}, \cA,\P)$ with values in $ \mathbb{R}^q$ and satisfying $\esp{\sup_{t\in[0,T]}|U_t|^2} <+\infty$. We also define $\cB^2 := \cS^2_1\times \cH^2_d$.

\vspace{2mm}

\section{The \emph{direct} and  \emph{Picard  algorithms}}\label{Description of the numerical methods}

We describe here the  numerical methods studied in this work. The first one,  the \emph{direct algorithm} is an adaptation  of the \emph{Deep BSDEs solver} introduced in \cite{weinan2017deep} to the linear specification of the parametric space that we use here. The second one, the \emph{Picard algorithm}, is new and is the main contribution of our work. We  also give here the main general convergence results related to the \emph{Picard algorithm}. The complexity analysis is postponed to the next section. 

The methods we introduce below have for goal to compute an approximation of the value function $u$, satisfying the PDE \eqref{eq pde}, at the initial time on a given domain or at a specific point. This lead us to introduce the following setup for the initial value $\cX_0$:

\begin{Assumption}\label{ass X0} One of the two following cases holds:
\begin{enumerate}[(i)]
\item The law of $\cX_0$ has compact support and is absolutely continuous with respect to the Lebesgue measure.
\item The law of $\cX_0$ is a Dirac mass at some point $x_0 \in \R^d$.
\end{enumerate}
\end{Assumption}

\vspace{2mm}
Most of our numerical applications are done in the setting of Assumption \ref{ass X0}(ii), see next section. Then, obviously, the approximation of the value function is known only at the point $x_0$ at the initial time. However, one should note that it could also be interesting to work in the setting of Assumption \ref{ass X0}(i) if one seeks to obtain an approximation of the \emph{whole} value function (on the support of $\cX_0$) at the initial time.

\subsection{Assumptions on the coefficients and connection with the semilinear PDE}
In this subsection, we first give the assumptions on the BSDE coefficients that will be required for our approach and then recall the connection with semilinear PDEs. In particular, under these assumptions, the underlying PDE admits a unique classical solution. Under an additional regularity assumption on the coefficients,  the unique solution to the PDE admits smooth derivatives of enough order which are controlled on the whole domain by known parameters. 
This additional regularity, together with a periodicity assumption, will be used to obtain our theoretical complexity result, see Section \ref{subse picard periodic}. For sake of simplicity, it is also assumed that the coefficients $b, \, \sigma$ and $f$ do not depend on time and that $f$ does not depend on the space variable.

\begin{Assumption}\label{assumption:coefficient}
\begin{enumerate}
\item[(i)] The coefficients $b$, $\sigma$, $f$ and $g$ are bounded, Lipschitz-continuous with respect to all variables and $g \in \mathcal{C}^{2+\alpha}(\mathbb{R}^d)$, for some $\alpha \in (0,1]$. We will denote by $L$ the Lipschitz-modulus of the map $f$.

\item[(ii)] The coefficient $a=\sigma \sigma^{\top}$ is uniformly elliptic, that is, there exists $\lambda_0\geq 1$ such that for any $(x, \zeta) \in (\mathbb{R}^d)^2$ it holds
\begin{equation}
\label{uniform:ellipticity:condition}
\lambda^{-1}_0 |\zeta|^2 \leq  a(x)\zeta \cdot \zeta  \leq \lambda_0 |\zeta|^2.
\end{equation}
\item[(iii)] For any $(i, j) \in \left\{1, \cdots, d\right\}^2$, the coefficients $b_i$, $\sigma_{i, j}$, $g$ belong to $\mathcal{C}^{2d+1}(\mathbb{R}^d, \mathbb{R})$ and $f$ belongs $\mathcal{C}^{2d+1}(\mathbb{R}\times \mathbb{R}^d, \mathbb{R})$. Moreover, their derivatives of any order up and equal to $2d+1$ are bounded and Lipschitz continuous.
\item[(iv)] The coefficients $b$, $\sigma$, $f$ and $g$ are periodic functions.
\end{enumerate}
\end{Assumption}

\noindent From now on, we will say that Assumption \ref{assumption:coefficient} holds if and only if Assumption \ref{assumption:coefficient} (i), (ii),  (iii)  and \textcolor{black}{(iv)} are satisfied.

\vspace{2mm}
 Under Assumption \ref{assumption:coefficient} (i) and (ii), it is known (see e.g. \cite{Pardoux:Peng:90}) that for any square integrable initial condition $\cX_0$ there exists a unique couple $(\cY, \cZ)\in \mathcal{B}^2$ satisfying equation \eqref{eq backward sde} $\P$-a.s. Moreover, from \cite{ladyzenskaya} Chapter VI and \cite{friedman2008partial} Chapter 7, the PDE \eqref{eq pde} admits a unique solution $u \in \mathcal{C}^{1,2}([0,T] \times \mathbb{R}^d, \mathbb{R})$ satisfying: there exists a positive constant $C$, depending on $T$ and the parameters appearing in Assumption \ref{assumption:coefficient} (i) and (ii), such that for all $(t, x) \in [0,T]\times \mathbb{R}^d$ 
$$
|u(t, x)| + |\partial_t u(t, x)| + |\nabla_x u(t, x)| + | \nabla^2_x u(t, x)| \leq C.  
$$ 

From \cite{SHIGEPENG1991,Pardoux:Peng:92,pardoux:tang}, the semilinear PDE \eqref{eq pde} and the BSDE \eqref{eq forward sde}-\eqref{eq backward sde} are connected, namely, for all $(t,x) \in [0,T) \times \mathbb{R}^d$, it holds
$$
\cY_t=u(t,\cX_t), \quad  \cZ_t=\sigma^{\top}\!(\cX_t)\nabla_x u(t,\cX_t).
$$

Finally, under Assumption \ref{assumption:coefficient}, still from  \cite{ladyzenskaya} Chapter IV and \cite{friedman2008partial} Chapter III, the unique solution $u$ to the PDE \eqref{eq pde} is smooth, namely, setting $v_i = (\sigma^{\top}\nabla_x u)_i$, $1\leq i \leq d$, for any $\bold{l} \in \mathbf{J}^\infty_{d, 2}$, $D^{\bold{l}}v_i(t, .)$ exists and is bounded. In particular, there exists a positive constant, depending on $T$ and the known parameters appearing in Assumption \ref{assumption:coefficient} (i), (ii) and (iii) such that for all $\bold{l} \in \mathbf{J}^\infty_{d, 2}$ and all $(t, x) \in [0,T]\times \mathbb{R}^d$,  
\begin{equation}
\label{bound:multi:indices:derivatives:pde:solution}
\max_{1\leq i\leq d}|D^{\bold{l}}v_i(t, x)| \leq C.  
\end{equation}

\subsection{\emph{Direct algorithm}}
\label{subse direct algo}

We first consider the approximation of the forward component \eqref{eq forward sde}.
Given an equidistant grid $\pi:=\set{t_0=0<\dots<t_n<\dots<t_N = T}$ of the time interval $[0,T]$, $t_n=n h$, $n=0, \cdots N$, with time-step $h:=T/N$, we denote by $W := (\cW_{t_n})_{0 \le n \le N}$ the discrete-time version of the Brownian motion $\cW$ and define $\Delta W_n = W_{t_{n+1}}-W_{t_{n}}$, $0 \le n \le N-1$.

We then introduce a standard Euler-Maruyama approximation scheme of $\cX$ on $\pi$ defined by $X_0 = \cX_0$ and for $0 \le n \le N-1$,
\begin{align}\label{eq de euler X}
X_{t_{n+1}} = X_{t_n} + b(X_{t_n}) h + \sigma(X_{t_n}) \Delta W_n \,.
\end{align}
 
Before discussing the approximation of the backward component, we here state an important lemma concerning the existence of two-sided Gaussian estimates for the transition density of the above Euler-Maruyama approximation scheme. These estimates will prove very useful in the sequel, when studying the theoretical complexity of the \emph{Picard algorithm}. We denote by $p^{\pi}(t_i, t_j,  x, \cdot )$ the transition density function of the Euler-Maruyama scheme starting from the point $x$ at time $t_i$ and taken at time $t_j$, with $0\leq t_i < t_j \leq T$. We refer e.g. to \cite{menozzi2010} for a proof of the following result.

\begin{Lemma}\label{Aronson:two:sided:bounds} Assume that the coefficients $b$ and $\sigma$ satisfies Assumption \ref{assumption:coefficient} (i) and (ii). There exist constants $ \mathfrak{c}:= \mathfrak{c}(\lambda_0, b, \sigma, d)\in (0,1]$ and $ \mathfrak{C}:=\mathfrak{C}(T, \lambda, b, \sigma, d)\geq1$ such that for any $(x, x') \in (\mathbb{R}^d)^2$ and for any $0\leq i < j \leq N$
\begin{equation}
\label{Aronson:bounds:transition:density}
\mathfrak{C}^{-1} p( \mathfrak{c}(t_j-t_i), x - x') \leq p^{\pi}(t_i, t_j, x, x') \leq \mathfrak{C} p(\mathfrak{c}^{-1}(t_j-t_i), x'-x) 
\end{equation}
\noindent where for any $(t,x) \in (0, \infty) \times \mathbb{R}^d$, $p(t, x) := (1/(2\pi t))^{d/2}\exp(-|x|^2/(2t))$.
\end{Lemma}

We now turn to the approximation of the backward component \eqref{eq backward sde}. We first introduce a linear parametrization of the process $\cZ$. For each discrete date $t_n \in \pi\setminus \set{T}$, we consider  a parametric functional approximation space $\mathscr{V}^z_{n}$ generated by a set of basis functions $(\psi_n^{k})_{1 \le k \le K^z_n}$, for  $0  \le n \le N-1$ and some positive integer $K^z_n$. The measurable functions $\psi_n^{k}:\R^d \mapsto \R$ have at most polynomial growth. Note that, for $n  \ge 1$, the specification of the basis function could depend on the time $t_n$, but in order to simplify the discussion, we let the number of basis functions be the same and set to $K$. Namely, $K^z_n = K$, for all $n \ge 1$. For $n=0$, the specification will depend on the nature of $\cX_0$: if Assumption \ref{ass X0}(i) holds, then we will set $K^z_0 = K$; if Assumption \ref{ass X0}(ii) holds, then we simply set $K^z_0 = 1$ and $\psi_0^1$ is a function satisfying $\psi_0^1(x_0)=1$. For latter use, we set:
\begin{align}\label{eq de bar K}
\widebar{K}^z := \sum_{n=0}^{N-1} K^z_n = K^z_0 + (N -1)K\;.
\end{align}

\noindent Remark that there is no need to introduce an approximation space at $T$ since the function $g$ is explicitly known.
\\
For $0 \le n \le N-1,$ each component of $(\sigma^\top\nabla_x u)(t_n,\cdot)$ should be approximated in an optimal way by a function in $\mathscr{V}^z_n$. The process $\cZ$ appearing in the dynamics of the controlled process $\cY$, that has to be optimized, is parametrized using the spaces $(\mathscr{V}^z_n)_{0\leq n \leq N-1}$. Namely, the $\mathbb{R}^d$-valued random variable $\cZ_{t_n}$ will be approximated by 
\begin{align}\label{eq de Z}
\sum_{k = 1}^{K_n^z}  \psi_n^{k}(\widehatbis{X}_{t_n}) \mathfrak{z}^{n,k}\,,
\end{align}
where $\mathfrak{z}^{n,k} \in \R^d$ for any $1 \le k \le K_n^z$ and $0 \le n \le N-1$. Importantly, we denote, for later use, $\mathfrak{z}^\top := ((\mathfrak{z}^{n,k})^\top)_{0 \le n \le N-1, 1 \le k \le K^z_n}$ so that $\mathfrak{z} \in \R^{d\widebar{K}^z}$.\\

\begin{Definition}[Class of discrete control process]\label{definition:control:process}
We let $\cH^{\pi,\psi}$ be the set of discrete control process ${Z}$ defined by: for $\mathfrak{z} \in \R^{d\widebar{K}^z}$, 
\begin{align}\label{eq de control}
 {Z}_{t_n} :=   \  \sum_{k = 1}^{K^z_n}  \psi^{k}_n(\widehatbis{X}_{t_n}) \mathfrak{z}^{n,k},\; \text{ for } 0 \le n \le N-1,
\end{align}

\noindent and where we set ${Z}_{t}= {Z}_{t_n}  ,\; t_n \le t < t_{n+1}, 0 \le n \le N-1$ with the convention $Z_T = 0$.\\
\end{Definition}

\begin{Remark} 
We insist on the fact that for a given $Z \in \cH^{\pi,\psi}$, the $\mathbb{R}^d$-valued random variable $Z_{t_n}$ depends only on $\mathfrak{z}^n$, for any $0 \le n \le N-1$. The approximation space we consider is a finite dimensional vector space. This notably differs from the recent works \cite{beck2017machine, han2017overcoming, hure2020deep} where a non-linear approximation using neural network is used.
\end{Remark}

The dynamics of $\cY$ given by \eqref{eq backward sde}, in turn, has to be approximated. As previously mentioned in the introduction, we first rewrite it in forward form as follows
$$
\cY_t = \cY_0 - \int_0^t f(\cY_s, \cZ_s) \, ds + \int_0^t \cZ_s \cdot d\cW_s, \, t\in[0,T], \, \text{ with }\cY_0 = u(0,\cX_0).
$$
The main goal of the algorithm is to obtain a good estimate of $u(0,.)$ on the support of $\cX_0$. In order to do so, 
we define the starting point of $Y$, standing for the approximation of $\cY$, by using a linear functional approximation space denoted $\mathscr{V}^y$, namely
\begin{align}\label{eq de approx Y at 0}
{Y}_{0} :=    \sum_{k = 1}^{K^y}  \psi^{k}_y(\widehatbis{X}_{0}) \mathfrak{y}^{k} \quad\text{with}\quad \mathfrak{y} \in \R^{K^y}.
\end{align}
The specification of $\mathscr{V}^y$ will depend also on the nature of $\cX_0$. Namely, if Assumption \ref{ass X0}(i) holds, then we set $K^y = K$, while if Assumption \ref{ass X0}(ii) holds, then we simply set $K^y = 1$ and $\psi_y^1$ is a function satisfying $\psi_y^1(x_0)=1$.

\vspace{2mm}
 Then, employing a standard Euler scheme on $\pi$ together with the above approximation $Z \in \cH^{\pi,\psi}$ of the control process $\cZ$, we are naturally led to consider the following approximation scheme for $\cY$.
\begin{Definition}\label{de discrete Y}
\begin{enumerate}[i)]
\item Given $\mathfrak{u}=(\mathfrak{y},\mathfrak{z}) \in \R^{K^y} \times \R^{d \bar{K}^z}$, we denote by $Z^{\mathfrak{u}} \in \cH^{\pi,\psi}$   the  discrete control process as given in \eqref{eq de control}. Then, the discrete controlled process ${Y}^{\mathfrak{u}}$ is defined as follows.
\begin{enumerate}
\item Initialization: Set
\begin{align} \label{eq de Y starting point}
{Y}^{\mathfrak{u}}_0 = \sum_{k = 1}^{K}  \psi_y^{k}(\widehatbis{X}_{0}) \mathfrak{y}^{k}\,.
\end{align}
\item Discrete version: for any $0 \le n \le N-1$: 
\begin{align}\label{eq de bar Y}
{Y}^{\mathfrak{u}}_{t_{n+1}} = {Y}^{\mathfrak{u}}_{t_{n}} - h f({Y}^{\mathfrak{u}}_{t_{n}},Z^{\mathfrak{u}}_{t_{n}} ) 
+Z^{\mathfrak{u}}_{t_{n}}\cdot \Delta W_n 
\end{align}
\noindent where we recall that $ \Delta W_n = W_{t_{n+1}}-W_{t_{n}}$.
\item Continuous version: for any $0 \le n \le N-1$ and any $t_n \le t < t_{n+1}$,
\begin{align}\label{eq de bar Y:continuous:time}
{Y}^{\mathfrak{u}}_t  = {Y}^{\mathfrak{u}}_{t_n}
- (t-t_n) f(Y^\mathfrak{u}_{t_n},Z^{\mathfrak{u}}_{t_n})+  Z^{\mathfrak{u}}_{t_n} \cdot (\cW_{t}-\cW_{t_n})
\end{align}
\end{enumerate}
\item Based on the previous step, we define $\cB^{\pi,\psi} \subset \cB^2$ as the set of processes $(Y^\mathfrak{u},Z^\mathfrak{u})$, with $Z^{\mathfrak{u}} \in \cH^{\pi,\psi}$, $Y^{\mathfrak{u}}$ defined as above for some $\mathfrak{u} \in \R^{K^y} \times \R^{d \bar{K}^z}$.
\end{enumerate}
\end{Definition}

\begin{Remark} 
Let us note that the discrete process $(\widehatbis{X}_t,Y^\mathfrak{u}_t,Z^\mathfrak{u}_t)_{t \in \pi}$ depends on $\cX_0 $ and $(W_t)_{t \in \pi}$ but we omit these dependences in the notation. 
\end{Remark}

\noindent 
The main idea of  \emph{approximation by learning} methods is to force the discrete controlled process $Y^{\mathfrak{u}}_T$ at maturity $T$  to match the approximated terminal condition $g(\widehatbis{X}_{T})$, by minimizing a loss function. Here, we work with the quadratic loss function, so that one faces  the optimization problem
\begin{align}\label{eq de theoretical direct algo}
 \inf_{\mathfrak{u}=(\mathfrak{y}, \mathfrak{z})  \in \R^{K^y} \times \R^{d \bar{K}^z}} \mathfrak{g}(\mathfrak{u}) :=\esp{\mathrm{G}(\cX_0,{W},\mathfrak{u})} 
\;\text{ with }\; \mathrm{G}(\cX_0,{W},\mathfrak{u}) = |g(\widehatbis{X}_T) - {Y}^{\mathfrak{u}}_T |^2\;.
\end{align}
 
\noindent However, one has to come up with a numerical procedure to compute the solution in practice. 

In order to numerically compute a solution to the optimization problem \eqref{eq de theoretical direct algo} (if any exists), one generally employs a stochastic approximation scheme such as a SGD algorithm. For an  overview of the theory of stochastic approximation, the reader may refer to \cite{Duflo1996}, \cite{Kushner2003} and \cite{benveniste2012adaptive} and to \cite{beck2017machine, weinan2017deep, han2017overcoming, hure2020deep, guselectnet} for applications to deep learning approximation of PDEs.

\vspace{2mm}

\noindent
We now describe the SGD algorithm that we implement in order to compute a solution $(\mathfrak{y},\mathfrak{z}) \in \R^{K^y} \times \mathbb{R}^{d \bar{K}^z}$ to the optimization problem \eqref{eq de theoretical direct algo}.

For a prescribed positive integer $M$ representing the number of steps in the stochastic algorithm and two deterministic non increasing sequences of positive real number $(\gamma^y_m)_{m\geq1} $ and $(\gamma^{z}_m)_{m\geq1}$  representing the learning rates, we design the following \emph{direct algorithm}.
\begin{Definition}(Implemented direct algorithm)\label{de implemented direct algo}
\begin{enumerate}
\item Simulate $M$ independent discrete paths of the Brownian motion $\mathfrak{W} = ({W}^{m})_{1 \le m \le M}$ and $M$ independent  samples of the initial condition $(\cX_0^m)_{1\le m \le M}$.
\item Initialization: select a random vector $\mathfrak{u}_0=(\mathfrak{y}_0,  \mathfrak{z}_0)$ with values in $\mathbb{R}^{K^y}\times \mathbb{R}^{d \bar{K}^z}$, independent of $\mathfrak{W}$ and $(\cX_0^m)_{1\le m \le M}$, and such that $\mathbb{E}[|\mathfrak{u}_0|^2]<\infty$.
\item Iteration: For $0 \le m \le M-1$, compute
\begin{align}
\mathfrak{y}_{m+1} &= \mathfrak{y}_m - \gamma^{y}_{m+1} \nabla_{\mathfrak{y}} \mathrm{G}\!\left(\cX_0^{m+1},{W}^{m+1},\mathfrak{u}_m \right),
\\
\mathfrak{z}^{n}_{m+1} & = \mathfrak{z}^{n}_m - \gamma^{z}_{m+1} \nabla_{\mathfrak{z}^{n}} \mathrm{G}\!\left(\cX_0^{m+1},{W}^{m+1},\mathfrak{u}_m\right),
\end{align}
for $0\le n \le N-1$.
\end{enumerate}
The output of the algorithm is then $\mathfrak{u}_M=(\mathfrak{y}_M, \mathfrak{z}_M)$.
%
\end{Definition}

\begin{Remark}
In order to analyse the asymptotic properties of stochastic approximation schemes, one usually chooses the learning sequences $(\gamma_m)_{m\geq1} = (\gamma^y_m)_{m\geq1} $ or  $(\gamma_m)_{m\geq1} = (\gamma^{z}_m)_{m\geq1} $ such that
 \begin{equation}
 \label{eq step decr}
 \sum_{m\geq1} \gamma_m = \infty \quad \mbox{and} \quad \sum_{m\geq1} \gamma^2_m < \infty,
 \end{equation}

\noindent see e.g. \cite{Duflo1996, Kushner2003, benveniste2012adaptive}.
\end{Remark}

\noindent The following lemma, whose proof is postponed to Section \ref{subse application direct algo}, provides the analytic expression of the local gradient functions $\nabla_\lambda \mathrm{G}\!\left(\cX_0,W,\mathfrak{u} \right)$, $\lambda \in \left\{ \mathfrak{y}, \mathfrak{z}^{n},  0\leq n\leq N-1\right\} $ appearing in the above SGD algorithm. It shows that $(\mathfrak{y}_{m+1}, \mathfrak{z}_{m+1})$ can be easily computed once $(Y^{\mathfrak{u}_m}, Z^{\mathfrak{u}_m})$ have been simulated for any $0 \le m \le M-1$.
\begin{Lemma}\label{pr expression deriv direct algo}
For $\lambda \in\left\{ \mathfrak{y}^k, 1\le k \le K^y \right\} \bigcup \left\{ \mathfrak{z}^{n, k},   1\le k \le K_n^z, 0\leq n\leq N-1\right\}$ and $\mathfrak{u} = (\mathfrak{y}, \mathfrak{z}) \in \R^{K^y} \times \R^{d\bar{K}^z}$, it holds
\begin{align}
\nabla_{\lambda} \mathrm{G}\!\left(\cX_0, {W},\mathfrak{u} \right) = - 2 (g(X_T) - {Y}^{\mathfrak{u}}_T) \nabla_{\lambda}{Y}^{\mathfrak{u}}_T 
\end{align}
with 
\begin{align*}
\nabla_{\mathfrak{y}^k}  Y^{\mathfrak{u}}_{T} & = \psi_y^k (X_0) \prod_{l = 0}^{N-1} \left(1- h \nabla_y f( Y^{\mathfrak{u}}_{t_l}, Z^\mathfrak{u}_{t_l}) \right),  \; 1\leq k \leq K^y, 
\end{align*}
\noindent and for any $0  \le n \le N-1$ and any $1 \le k \le K_n^z$,
\begin{align*}
 \nabla_{\mathfrak{z}^{n,k}} Y^{\mathfrak{u}}_{T} =   \psi_n^{k}(X_{t_n})   \left( \Delta W_{n}^{\top}  - h\nabla_z f( Y^{\mathfrak{u}}_{t_{n}}, Z^\mathfrak{u}_{t_{n}}) \right ) \prod_{l = n+1}^{N-1} \left(1- h \nabla_y f( Y^{\mathfrak{u}}_{t_l}, Z^\mathfrak{u}_{t_l}) \right)
\end{align*}
with the convention $\prod_{\emptyset}=1$.
\end{Lemma}
\vspace{2mm}
\noindent Under  Assumption \ref{ass X0}, Assumption \ref{assumption:coefficient} (i) and Assumption \ref{as psi function}, the well-posedness of Algorithm \ref{de implemented direct algo}, that is, the fact that it holds
\begin{align*}
\argmin_{\mathfrak{u} \in \R^{K^y}\times \R^{d \bar{K}^z}} \mathfrak{g}(\mathfrak{u}) \neq \emptyset \;,
\end{align*}
is proved in \Cref{le well-posedness}. \\
Additionally, for any $\mathfrak{u}^\star \in \argmin_{\mathfrak{u} \in \R^{K^y}\times \R^{d \bar{K}^z}} \mathfrak{g}(\mathfrak{u})$, we show in Proposition \ref{pr conv direct algo} that
\begin{align*}
\quad \esp{|u(0,\cX_0)-Y_0^{\mathfrak{u}^\star}|^2 +  h \sum_{n=0}^{N-1}  |\cZ_{t_{n}}-Z^{\mathfrak{u}^\star}_{t_{n}}|^2}\le C\left(  \cE_{\pi} + \cE_{\psi}
\right) \,,
\end{align*}
for some positive constant $C$. The quantities $ \cE_{\pi} $ and $\cE_{\psi}$ represents the discrete-time error and the error due to the approximation in the functional spaces $(\mathscr{V}^y,\mathscr{V}^z_n)$, respectively. They are defined by
\begin{align}
\cE_{\pi} &:=  \esp{ \sum_{n=0}^{N-1}  \int_{t_{n}}^{t_{n+1}}\left( |\cY_s-\cY_{t_{n}}|^2 + |\cZ_s-\cZ_{t_{n}}|^2 +
\textcolor{black}{ |\cX_{t_{n}} - \widehatbis{X}_{t_n}|^2} \right)ds } \, \label{eq error disc time}
\end{align}
and
\begin{align}
\cE_{\psi} &:= \inf_{\mathfrak{u} \in \R^{K^y}\times \R^{d \bar{K}^z}} \esp{|u(0,\cX_0) - Y^\mathfrak{u}_0|^2 + \sum_{n=0}^{N-1} h |(\sigma^\top\nabla_x u)(t_n,\widehatbis{X}_{t_n})- Z^\mathfrak{u}_{t_n}|^2} \,.\label{eq error approx space brut}
\end{align}
Let us mention, for later use, that, in the setting of  Assumptions \ref{ass X0} and \ref{assumption:coefficient}(i), 
\begin{align}\label{eq control disc error}
\cE_{\pi} &\le C h\;, 
\end{align}
for some positive constant $C$, see e.g. Ma and Zhang \cite{ma2002path} and Pag\`es \cite{pages2018numerical}. 

\vspace{2mm}
\noindent We shall not seek to obtain theoretical convergence results for the \emph{direct algorithm} itself. However, we illustrate its performance numerically in Section \ref{se main sparse grids} when using sparse grids approximations \cite{bungartz2004sparse}.

\begin{Remark}\label{re num complexity direct algo}
The numerical complexity $\cC$ will be measured by the number of coefficients update realized to obtain the approximation. From the previous description, we obtain straightforwardly that the complexity at worst satisfies
\begin{align}\label{eq num complexity direct algo}
\cC = O_d(NKM)\;.
\end{align}
\end{Remark}

\subsection{A \emph{Picard algorithm}}
\label{subse picard algo}
An issue with the above algorithm comes from the fact that the optimization problem \eqref{eq de theoretical direct algo} is generally not convex. Even though $\mathfrak{u} \mapsto (Y_0^\mathfrak{u},Z^\mathfrak{u})$ is linear for our choice of parametrisation, in general the mapping $\mathfrak{u} \rightarrow Y^{\mathfrak{u}}$ is non-linear since $f$ itself is non-linear. As a consequence, in practical implementation, we have no guarantee that the algorithm converges to local or global minima. On top of practical problems, this renders the theoretical analysis of the implemented \emph{direct algorithm} difficult, in particular if one wants to obtain rates of convergence to assess precisely the numerical complexity of the method. 

In this section, we introduce a \emph{Picard algorithm} which transforms this non-convex optimisation problem into a sequence of linear-quadratic optimization problems. This is done by using the special structure of the original problem. Indeed, it is well known that the solution of the BSDE \eqref{eq backward sde} itself is obtained as the limit of a sequence of Picard iterations, see e.g. \cite{el1997backward} and \cite{bender2007forward} from a numerical perspective. \\


\subsubsection{Theoretical Picard algorithm}

\noindent Our \emph{Picard algorithm} is  based on the iteration of the following operator:
\begin{align}\label{eq de Phi}
\R^{K^y}\times \R^{d \bar{K}^z} \ni \tilde{\mathfrak{u}}  \mapsto \Phi(\tilde{\mathfrak{u}} ) :=  \check{\mathfrak{u}} \in  \R^{K^y}\times \R^{d \bar{K}^z} \end{align}
where, 
\begin{align}\label{eq optim step}
\check{\mathfrak{u}} = \argmin_{\mathfrak{u} \in \R^{K^y}\times \R^{d \bar{K}^z}  }\esp{|g(\widehatbis{X}_T) - {U}^{\tilde{\mathfrak{u}},\mathfrak{u}}_T|^2}\,.
\end{align}

\noindent In the above expectation, the process $X$ is the Euler-Maruyama approximation scheme on the time grid $\pi$ with dynamics \eqref{eq de euler X} and ${U}^{\tilde{\mathfrak{u}},\mathfrak{u}}$ (simply denoted as ${U}$ below) is given by the following {decoupling approximation scheme}:
\begin{enumerate}
\item For $\tilde{\mathfrak{u}} \in \R^{K^y}\times \R^{d \bar{K}^z}$, we first consider $({Y}^{\tilde{\mathfrak{u}}}, Z^{\tilde{\mathfrak{u}}}) \in \cB^{\pi,\psi} $ as introduced in Definition \ref{de discrete Y}.  
\item Then, for any $\mathfrak{u} \in \R^{K^y}\times \R^{d \bar{K}^z}$, consider the discrete control process $Z^\mathfrak{u} \in \cH^{\pi,\psi}$ as introduced in \eqref{eq de control} of Definition \ref{definition:control:process} and define the control process ${U}^{\tilde{\mathfrak{u}},\mathfrak{u}} $ by
\begin{align} \label{eq de U starting point}
{U}^{\tilde{\mathfrak{u}},\mathfrak{u}} _0 = Y_0^\mathfrak{u} 
\end{align}
recall \eqref{eq de Y starting point}
and for any $0 \le n \le N-1$, 
\begin{align}\label{eq scheme picard}
{U}_{t_{n+1}}^{\tilde{\mathfrak{u}},\mathfrak{u}} = {U}_{t_{n}}^{\tilde{\mathfrak{u}},\mathfrak{u}} - h f({Y}^{\tilde{\mathfrak{u}}}_{t_n}, Z^{\tilde{\mathfrak{u}}}_{t_n} ) 
+ Z^{\mathfrak{u}}_{t_n} \cdot (W_{t_{n+1}}-W_{t_{n}})\,,
\end{align}
and its continuous version, for any $t_n\le t \le t_{n+1}$,
\begin{align*}
{U}_{t}^{\tilde{\mathfrak{u}},\mathfrak{u}} = {U}_{t_{n}}^{\tilde{\mathfrak{u}},\mathfrak{u}} - (t-t_n) f({Y}^{\tilde{\mathfrak{u}}}_{t_n}, Z^{\tilde{\mathfrak{u}}}_{t_n} ) 
+ Z^{\mathfrak{u}}_{t_n} \cdot (\cW_t-W_{t_{n}})\,.
\end{align*}
\end{enumerate}
Note that under Assumption \ref{assumption:coefficient}(i) it holds $\esp{\sup_{t \in [0,T]} |{U}_{t}^{\tilde{\mathfrak{u}},\mathfrak{u}}|^2} < +\infty$.

\begin{Definition}[Theoretical Picard algorithm] \label{de picard scheme theo} For a prescribed positive integer $P$:
\begin{enumerate}
\item Initialization: set $\mathfrak{u}^0 \in \R^{K^y}\times \R^{d \bar{K}^z}$.
\item Iteration: for $1 \le p \le P$, compute: ${\mathfrak{u}}^p = \Phi({\mathfrak{u}}^{p-1})$. 
\end{enumerate}
The output of the algorithm is then ${\mathfrak{u}}^P$.
\end{Definition}


\subsubsection{Well-posedness of the theoretical algorithm} 
The main novelty compared to the optimization problem \eqref{eq de theoretical direct algo} comes from the fact that the map $\mathfrak{u} \mapsto {U}^{\tilde{\mathfrak{u}},\mathfrak{u}}$ is now linear. This linearity is achieved by freezing the driver $f$ in the dynamics of the control process along the process $(Y^{\tilde{\mathfrak{u}}}, Z^{\tilde{\mathfrak{u}}}) \in \cB^{\pi,\psi}$. The parameter $\tilde{\mathfrak{u}}\in \R^{K^y}\times \R^{d \bar{K}^z}$ is then updated through the Picard iteration procedure. This is of course the main purpose of this \emph{Picard algorithm}, compared to the \emph{direct algorithm} of Section \ref{subse direct algo}. At this stage, the above algorithm is theoretical and its solution (if any exists) still needs to be numerically approximated. This will be discussed in full details in the next section.

We here discuss the well-posedness of the optimization problem \eqref{eq optim step}. We first introduce some notations  that will be useful in the sequel to study the \emph{Picard algorithm} as iterated least-square optimization problems.

\noindent First, to clarify the linear structure, we introduce the following notations 
\begin{enumerate}[i)]
\item  For $1 \le k \le K^y$, $\theta^{k} := \psi_y^{k}(\widehatbis{X}_0)$.

\item For $0 \le n \le N-1$, $1 \le k \le K$, the $\mathbb{R}^d$-valued random vectors $\omega^{n , k}$ is defined by
\begin{align} \label{eq de omega}
\omega^{n, k}=  \Psi^{n, k} \Delta W_n\, \text{ with }  \Psi^{n, k} = \psi^k_n(X_{t_n}) \;,
\end{align}
and we set $\omega^\top := ((\omega^{n,k})^\top)_{0 \le n \le N-1, 1 \le k \le K^z_n}$ (so that $\omega$ is an $\R^{d\bar{K}^z}$-valued random vector).
\item the random vector $\Omega = (\theta^\top, \omega^\top)^\top$ which takes values in $\R^{K^y}\times \R^{d \bar{K}^z}$.
\end{enumerate}
Note that both $\omega$ and $\Omega$ depends on $W$ and $\cX_0$, but we will omit this in the notation. Then, we rewrite
\begin{align}\label{eq rewrite objective}
g(\widehatbis{X}_T) - U^{\tilde{\mathfrak{u}},\mathfrak{u}}_T & = \mathfrak{G}^{\tilde{\mathfrak{u}}} - \mathfrak{u} \cdot \Omega
\end{align}
where 
\begin{align} \label{eq de Gamma}
\mathfrak{G}^{\tilde{\mathfrak{u}}} = \mathfrak{G}^{\tilde{\mathfrak{u}}}(\cX_0,W) := g(\widehatbis{X}_T) {+} \sum_{n=0}^{N-1} hf(Y^{\tilde{\mathfrak{u}}}_{t_n},Z^{\tilde{\mathfrak{u}}}_{t_n})\;.
\end{align}
Thus, the optimization problem \eqref{eq optim step} is given by
\begin{align}\label{eq optim step bis}
\check{\mathfrak{u}}  = \mathrm{argmin}_{\mathfrak{u} \in \R^{K^y}\times \R^{d\bar{K}^z}} \mathfrak{H}(\tilde{\mathfrak{u}},\mathfrak{u})
\quad
\text{with}
\quad
 \mathfrak{H}(\tilde{\mathfrak{u}},\mathfrak{u}) : = \esp{|\mathfrak{G}^{\tilde{\mathfrak{u}}} - \mathfrak{u} \cdot \Omega|^2}
\end{align}
\noindent and simply reads as a Linear-Quadratic optimization problem. Classically, we  introduce semi-norms on the parameter spaces.
\begin{Definition}\label{de double - triple norm}
For $\mathfrak{u}=(\mathfrak{y},\mathfrak{z}) \in \R^{K^y}\times \R^{d\bar{K}^z} $, we define
\begin{align*}
\| \mathfrak{y} \|_{y}^2 &:=  \esp{| \mathfrak{y} \cdot \theta |^2} , \;\| \mathfrak{z} \|_{z}^2 :=  \esp{|\mathfrak{z} \cdot \omega |^2} \text{ and } \vvvert \mathfrak{u} \vvvert^2 :=  \esp{|\mathfrak{u} \cdot \Omega |^2}.
\end{align*}
Let us insist on the fact that these quantities depend on the choices of $\pi$,$\psi$ though this is not reflected in the notation.
\end{Definition}

\begin{Remark}
\begin{enumerate}[i)]
\item
Observe that from the very definition of the random vector $\Omega$, for any $\mathfrak{u} =(\mathfrak{y},\mathfrak{z})\in \R^{K^y}\times \R^{d \bar{K}^z}$, it holds
\begin{align}\label{eq link triple double}
\vvvert \mathfrak{u} \vvvert^2 = \| \mathfrak{y} \|_{y}^2 + \| \mathfrak{z} \|_z^2.
\end{align}

\item With the notations of Section \ref{subse direct algo}, the following relations hold
\begin{align} \label{eq rewritte norms}
\vvvert \mathfrak{u} \vvvert^2 &= \esp{\left|Y^\mathfrak{u}_0 + \int_{0}^T Z^\mathfrak{u}_t \ud \cW_t \right|^2},\;
\| \mathfrak{y} \|_y^2 = \esp{|Y^\mathfrak{u}_0|^2}\;
\text{and}
\;
\| \mathfrak{z} \|_z^2 = \sum_{n=0}^{N-1} h \esp{|Z^{\mathfrak{u}}_{t_n}|^2}\;,
\end{align}

\noindent for any $\mathfrak{u} =(\mathfrak{y},\mathfrak{z})\in \R^{K^y}\times \R^{d \bar{K}^z}$.
\item For later use, see Section \ref{subsubse full algo}, we also note that $\|\cdot\|_{z}$ develops as 
\begin{align}
\| \mathfrak{z} \|_z^2 
&= \sum_{n=0}^{N-1} h\sum_{l=1}^d (\mathfrak{z}^{n,\cdot}_l)^\top \esp{\Psi^{n,\cdot}(\Psi^{n,\cdot})^\top}\mathfrak{z}^{n,\cdot}_l\;
\end{align}

\noindent by using the independence of the increments $(\Delta W_{n})_l$, for $0 \le n \le N-1$ and $l \in \set{1,\dots,d}$ and where we used the notations $\mathfrak{z}^{n, \cdot}_{l}=(\mathfrak{z}^{n, 1}_l, \ldots, \mathfrak{z}^{n,K^z_n}_l)$ and $\Psi^{n,\cdot}=(\Psi^{n,1}, \ldots, \Psi^{n,K^z_n})$.

\end{enumerate}
\end{Remark}

\noindent A key assumption to ensure the well-posedness of our approach is the following.
\begin{Assumption}\label{as psi function}
There exist  two positive constants $\kappa_{K} \ge 1 \ge \alpha_K$  such that for any $(\mathfrak{y},\mathfrak{z}) \in   \R^{K^y}\times \R^{d \bar{K}^z}$
\begin{align*}
\alpha_{K} |\mathfrak{y}|^2 \le \|  \mathfrak{y} \|_{y}^2 \le \kappa_{K} |\mathfrak{y}|^2 
\; \text{ and } \; h\alpha_{K} |\mathfrak{z}|^2 \le \|  \mathfrak{z} \|_{z}^2 \le h \kappa_{K} |\mathfrak{z}|^2\;.
\end{align*}
\end{Assumption}


\begin{Lemma} \label{le wellposedness picard} For all $(\tilde{\mathfrak{u}},\mathfrak{u}) \in \left( \R^{K^y}\times\R^{d\bar{K}^z} \right)^2 $, it holds
 \begin{align}\label{eq carac second deriv}
\mathfrak{u}^\top \nabla^2_{\mathfrak{u}}  \mathfrak{H}(\tilde{\mathfrak{u}},\mathfrak{u}) \mathfrak{u}
&= 2 \vvvert \mathfrak{u} \vvvert^2\,,
 \end{align}
 where $\nabla^2_{\mathfrak{u}} \mathfrak{H}$ denotes the Hessian of the function $\mathfrak{u} \mapsto \mathfrak{H}(\tilde{\mathfrak{u}},\mathfrak{u}) $.
 \\
 Moreover, under Assumption \ref{as psi function}, the optimization problem \eqref{eq optim step} admits a unique solution and for any $\tilde{\mathfrak{u}} \in \R^{K^y}\times \R^{d \bar{K}^z}$ and $(\mathfrak{y},\mathfrak{z}) \in \R^{K^y}\times \R^{d \bar{K}^z}$, it holds
 \begin{align}
 & (  \mathfrak{y} - \check{\mathfrak{y}}) \cdot \nabla_{\mathfrak{y}} \mathfrak{H}(\tilde{\mathfrak{u}},\mathfrak{u})  = 2 \| \mathfrak{y} - \check{\mathfrak{y}}\|_{y}^2 \geq 2 \alpha_{K} |\mathfrak{y} - \check{\mathfrak{y}}|^2
 \label{eq convex prop:1}\\
& (  \mathfrak{z} - \check{\mathfrak{z}}) \cdot \nabla_{\mathfrak{z}} \mathfrak{H}(\tilde{\mathfrak{u}},\mathfrak{u})  = 2 \| \mathfrak{z} - \check{\mathfrak{z}}\|_{z}^2 \geq 2 h \alpha_{K} |\mathfrak{z} - \check{\mathfrak{z}}|^2 \label{eq convex prop:2} 
 \end{align}
where $\check{\mathfrak{u}}=(\check{\mathfrak{y}},\check{\mathfrak{z}}) \in \R^{K^y}\times \R^{d \bar{K}^z}$ is the unique solution to \eqref{eq optim step}.\\

\end{Lemma}

\proof From \eqref{eq optim step bis}, we straightforwardly compute
\begin{align}
\nabla_\mathfrak{u} \mathfrak{H}(\tilde{\mathfrak{u}}, \mathfrak{u} ) &= -2\esp{(\mathfrak{G}^{\tilde{\mathfrak{u}}} - \mathfrak{u} \cdot \Omega)\Omega} \text{ and } \nabla^2_{\mathfrak{u}} \mathfrak{H}(\tilde{\mathfrak{u}}, \mathfrak{u} ) = 2\esp{\Omega\Omega^\top}\;
\end{align}

\noindent which, recalling Definition \ref{de double - triple norm}, directly yields \eqref{eq carac second deriv}.
In particular, we have
\begin{align}
\nabla_{\!\mathfrak{y}} \mathfrak{H}(\tilde{\mathfrak{u}}, \mathfrak{u}) & = 
-2\esp{(\mathfrak{G}^{\tilde{\mathfrak{u}}} - \mathfrak{u} \cdot \Omega)\theta} 
= -2\esp{(\mathfrak{G}^{\tilde{\mathfrak{u}}} - \mathfrak{y} \cdot \theta)\theta}
\label{eq nabla y}
\end{align}
and, for any $0 \le n \le N-1$ and any $1 \leq l \leq d$,
\begin{align}
\nabla_{\!\mathfrak{z}^{n,\cdot}_l} \mathfrak{H}(\tilde{\mathfrak{u}}, \mathfrak{u}) 
=-2 \esp{(\mathfrak{G}^{\tilde{\mathfrak{u}}} - \Omega \cdot \mathfrak{u})\omega^{n,\cdot}_l}
=-2 \esp{(\mathfrak{G}^{\tilde{\mathfrak{u}}} - \omega^{n,\cdot}_l \cdot \mathfrak{z}^{n,\cdot}_l)\omega^{n,\cdot}_l}\;.
\label{eq nabla theta}
\end{align}
Under Assumption \ref{as psi function}, we deduce from \eqref{eq carac second deriv} that the problem is strictly convex and has a unique (global) minimum $\check{\mathfrak{u}} = (\check{\mathfrak{y}},\check{\mathfrak{z}})$. 
{The inequalities \eqref{eq convex prop:1} and \eqref{eq convex prop:2} then follow from \eqref{eq nabla y} and \eqref{eq nabla theta} combined with the fact that $\nabla_\mathfrak{z} \mathfrak{H}(\tilde{\mathfrak{u}}, \check{\mathfrak{u}}) = 0$.}
\eproof

\subsubsection{Algorithm implementation}
\label{subsubse full algo}

From a practical point of view, the sequence of theoretical Linear-Quadratic optimization problem described in the previous section has to be approximated. Due to the possibly high dimension of the matrix $\mathbb{E}[\Omega \Omega^{\top}]$, we will rely on a SGD algorithm\footnote{Since it is also the procedure used for the \emph{direct algorithm}, the numerical comparison between the two will be more relevant.} to compute the unique solution to \eqref{eq optim step}. Indeed, for a fixed vector $\tilde{\mathfrak{u}}$  in $\R^{K^y}\times \R^{d \bar{K}^z}$, the key point is to observe that the unique minimizer $\check{\mathfrak{u}}$ is the unique solution to the equation
\begin{align}\label{eq optimality condition}
\nabla_\mathfrak{u} \mathfrak{H}(\tilde{\mathfrak{u}}, \mathfrak{u}) =0.
\end{align}

We importantly remark, using \eqref{eq nabla y} and \eqref{eq nabla theta} that the above relation \eqref{eq optimality condition} holds true
if and only if
\begin{align} 
\esp{H^y(\cX_0,W,\tilde{\mathfrak{u}},\mathfrak{y})} = 0, \quad \text{ and } \quad
\esp{H^{n,l}(\cX_0,W,\tilde{\mathfrak{u}},\mathfrak{z}^{n,\cdot}_l)} = 0, \label{eq optimality condition:detailed}
\end{align}

\noindent where $H^y$ is a map from $\R^d\times(\mathbb{R}^d)^{N} \times (\R^{K^y}\times \R^{d \bar{K}^z}) \times \R^{K^y}$ to $\mathbb{R}^{K^y}$  and defined by 
\begin{align}\label{def:funct:H:1st:part}
H^y(\cX_0,W,\tilde{\mathfrak{u}},\mathfrak{y}) :=-\frac{2}{{\beta_K}}  (\mathfrak{G}^{\tilde{\mathfrak{u}}} - \mathfrak{y} \cdot \theta)\theta,
\end{align}
and $H^{n, l}$ are maps defined on $\R^d\times(\mathbb{R}^d)^{N} \times (\R^{K^y}\times \R^{d \bar{K}^z}) \times \R^{K_n^z}$ taking values in $\mathbb{R}^{K_n^z}$ and given by 
\begin{align}\label{def:funct:H:2nd:part}
H^{n, l}(\cX_0,W,\tilde{\mathfrak{u}},\mathfrak{z}^{n,\cdot}_l) := -\frac{2}{{\beta_K \sqrt{h}}} (\mathfrak{G}^{\tilde{\mathfrak{u}}} - \omega^{n,\cdot}_l \cdot \mathfrak{z}^{n,\cdot}_l)\omega^{n,\cdot}_l, \quad  0 \le n \le N-1, \quad 1 \leq l \leq d.
\end{align}

\noindent We importantly point out that we abuse the notation in \eqref{def:funct:H:1st:part} and \eqref{def:funct:H:2nd:part} since the variable $(\cX_0,W)$ stands for a vector of $\R^d\times(\mathbb{R}^d)^{N}$ and $(\cX_0,W)\mapsto \mathfrak{G}^{\tilde{\mathfrak{u}}} = \mathfrak{G}^{\tilde{\mathfrak{u}}}(\cX_0,W)$ is also defined by \eqref{eq de Gamma} while in \eqref{eq optimality condition:detailed} the random vector $W=(\cW_{t_n})_{1\leq n \leq N}$ stands for the discrete path of the Brownian motion $\cW$ and $\cX_0$ for the starting value of $\widehatbis{X}$.

\noindent In \eqref{def:funct:H:1st:part} and \eqref{def:funct:H:2nd:part}, the deterministic constant {$\beta_K$ corresponding to a normalizing factor is introduced in order to control the $L^{2}(\mathbb{P})$-moment of the random vectors $(\mathfrak{G}^{\tilde{\mathfrak{u}}} - \mathfrak{y} \cdot \theta)\theta$ and $(\mathfrak{G}^{\tilde{\mathfrak{u}}} - \omega^{n,\cdot}_l \cdot \mathfrak{z}^{n,\cdot}_l)\omega^{n,\cdot}_l$.} \textcolor{black}{Namely, we select $\beta_K$ large enough so that 
\begin{align}\label{eq de beta psi}
(\beta_{K})^2 \ge (1+\esp{|\theta|^4}) \vee \max_{0 \le n \le N-1, 1 \le l \le d} (1+\esp{|\tilde{\omega}^{n,\cdot}_l|^4})
\end{align}
with $\tilde{\omega}^{n,\cdot}_\ell = \frac{{\omega}^{n,\cdot}_\ell }{\sqrt{h}}$.} 
Let us insist on the fact that the chosen $\beta_K$ above should be uniform for all time grid $\pi$. It depends only on the level of approximation coming from the definition of the approximation spaces. This qualitative level of approximation is controlled by the number of basis function per time step, namely $K$.

\noindent For latter use, comparing \eqref{def:funct:H:1st:part} to \eqref{eq nabla y} and \eqref{def:funct:H:2nd:part} to 
\eqref{eq nabla theta}, we remark that
\begin{align}\label{eq link h and H}
\esp{H^y(\cX_0,W,\tilde{\mathfrak{u}},\mathfrak{y})}  = \frac{1}{\textcolor{black}{\beta_K}} \nabla_{\!\mathfrak{y}} \mathfrak{H}(\tilde{\mathfrak{u}}, \mathfrak{u}) \,
\textnormal{ and } \,
\esp{H^{n, l}(\cX_0,W,\tilde{\mathfrak{u}},\mathfrak{z}^{n,\cdot}_l)} = 
\frac{1}{\textcolor{black}{\beta_K \sqrt{h}}} \nabla_{\!\mathfrak{z}^{n,\cdot}_l} \mathfrak{H}(\tilde{\mathfrak{u}}, \mathfrak{u}) \,,
\end{align}
for $0 \le n \le N-1$, $l \in \set{1, \dots, d}$ and  $(\tilde{\mathfrak{u}},\mathfrak{u})\in (\R^{K^y}\times\R^{d\bar{K}^z})^2$.

\vspace{4mm}
The implemented \emph{Picard algorithm} is obtained by iterating a stochastic gradient operator which is the counterpart of $\Phi$ defined by \eqref{eq de Phi} obtained by the numerical approximation that we now introduce.

\begin{Definition}\label{de sgd for phi}
Let $M$ be a positive integer. Let $\mathfrak{W} := ({W}^{m})_{1 \le m \le M}$, be $M$ discrete paths along the time grid $\pi$ of the Brownian motion $\cW$,  $\mathfrak{X}_0 := (\cX_0^m)_{1\le m \le M}$, be $M$ independent  samples of the initial condition (and independent from $\mathfrak{W}$) and $(\gamma_m)_{m\ge1}$ a deterministic sequence of positive real numbers satisfying: 
 \begin{equation}
 \label{eq step decr}
 \sum_{m\geq1} \gamma_m = \infty \quad \mbox{and} \quad \sum_{m\geq1} \gamma^2_m < \infty. 
 \end{equation}
 We set, for all $m\ge 1$, 
 \begin{align}\label{eq de gamama}
 \gamma^y_m = \gamma_m \quad \text{ and }  \quad \gamma^z_m = \frac{\gamma_m}{\sqrt{h}} \;. 
 \end{align}
 {\color{black}
  Let $ \mathfrak{u}_0 = (\mathfrak{y}_0, \mathfrak{z}_0)$ be a random vector taking values in $\R^{K^y}\times\R^{d\bar{K}^z}$, independent of $(\mathfrak{X}_0,\mathfrak{W})$ and such that $\mathbb{E}[|\mathfrak{u}_0|^2]<\infty$.
 
The operator $\Phi_M$, parametrized by $(\mathfrak{u}_0 ,\mathfrak{X}_0,\mathfrak{W})$, is given by
\begin{align}\label{eq de Phi:SG:approximation}
\R^{K^y}\times \R^{d \bar{K}^z} \ni \tilde{\mathfrak{u}}  \mapsto \Phi_M(\mathfrak{u}_0, \mathfrak{X}_0,\mathfrak{W},\tilde{\mathfrak{u}} ) = \mathfrak{u}_M
\end{align}
where $\mathfrak{u}_M$ is the output of the SGD algorithm after $M$ steps and is obtained as follows:
\begin{enumerate}
\item The initial value is set to $\mathfrak{u}_0$.
\item Iteration: For $0 \le m \le M-1$, compute 
\begin{align}
\mathfrak{y}_{m+1} & = \mathfrak{y}_{m} - \gamma^y_{m+1} H^{y}(\cX_0^{m+1},W^{(m+1)}, \tilde{\mathfrak{u}},\mathfrak{y}_{m}) \;, \label{RM iteration 1}
\\
&\textnormal{and} \nonumber
\\
\mathfrak{z}^{n,\cdot}_{l,m+1} & = \mathfrak{z}^{n,\cdot}_{l,m} - \gamma^z_{m+1} H^{n,l}(\cX_0^{m+1},W^{(m+1)}, \tilde{\mathfrak{u}},\mathfrak{z}^{n,\cdot}_{l,m}) \;, \label{RM iteration 2}
\end{align}
for any $0\le n \le N-1$ and any $1\leq  l \leq d$.
\end{enumerate}
}
\end{Definition}

\begin{Definition}[Implemented \emph{Picard algorithm}] \label{de picard scheme full} For a prescribed positive integer $P$:
\begin{enumerate}
{\color{black}
\item Initialization: Select a random vector $\mathfrak{u}^0_0$ taking values in $\R^{K^y}\times\R^{d\bar{K}^z}$ such that $\mathbb{E}[|\mathfrak{u}^0_0|^2]<\infty$. Set $\mathfrak{u}^0_M:=\mathfrak{u}^0_0$.
\item Iteration: for $1 \le p \le P$, simulate independently a set of $M$ independent discrete paths $\mathfrak{W}^p$ of the Brownian motion $\cW$, independent initial condition $\mathfrak{X}_0^p$ and an initial starting point $\mathfrak{u}^p_0$ (independently also of $\mathfrak{u}^0_0$, $\mathfrak{u}^j_0$ and of the previous $(\mathfrak{X}_0^j,\mathfrak{W}^j)$, $1 \leq j\leq p-1$), and compute $\mathfrak{u}^p_M := \Phi_M(\mathfrak{u}^{p}_0,\mathfrak{X}_0^p,\mathfrak{W}^p,\mathfrak{u}^{p-1}_M)$ as in Definition \ref{de sgd for phi}.
}
\end{enumerate}
The output of the algorithm is then $\mathfrak{u}^P_M $.
\end{Definition}

\begin{Remark}\label{re num picard algo}
\begin{enumerate}[i)]

\item 
The choice of the learning sequence $\gamma$ for the SGD algorithm \eqref{RM iteration 1}, \eqref{RM iteration 2} might be delicate in practice, see e.g. Section \ref{subsubse numerics one}. 
\item \textcolor{black}{The initialisation is random in the above algorithm. We do not always follow this procedure in our numerical experiments, see Section \ref{se main sparse grids}.}
\item 
The numerical complexity $\cC$ of the full algorithm is the sum of the local complexity of each SGD algorithm so that 
\begin{align}\label{eq num complexity direct algo}
\cC = O_d( PNKM) \,.
\end{align}
\end{enumerate}
\end{Remark}

\noindent Using the output $\mathfrak{u}^P_M $ of the \emph{Picard algorithm}, we set the approximating function at time $0$ to be:
\begin{align}\label{eq de approx u}
\cU^{P}_M(x) := \sum_{k=1}^{K^y} (\mathfrak{y}^P_M)^k \psi_y^k(x)\;,
\end{align}
recalling \eqref{eq de approx Y at 0}.

We then aim to control the following mean squared error:
\begin{align}\label{eq de error algorithm}
\cE_{\mathrm{MSE}} := \esp{\left|\cU^{P}_M(\cX_0) - u(0,\cX_0)\right|^2} \;.
\end{align}

We obtain an explicit upper bound on the mean-squared error when specifying the parameters of the algorithm as follows.
For $\gamma>0$, $\rho \in (\frac12,1)$, we set $\gamma_m := \gamma m^{-\rho}$, $m \ge 1$ and we assume that the number of steps $M$ in the SGD algorithm satisfies, for some $\eta \ge 0$,
\begin{align}\label{eq constraint alpha-beta-M 1 intro}
\gamma\frac{\alpha_{K}}{\beta_{K}}M^{1-\rho} \ge \frac{\sqrt{2}}2\quad 
\text{ and }\quad
 \frac{\kappa_{K}}{h \wedge \alpha_{K}}\left(e^{- 2\sqrt{2}\ln(2)\gamma\frac{\alpha_{K}}{\beta_{K}}M^{1-\rho}} 
+  \frac{ \beta_{K}}{\alpha_{K} M^\rho }
\right) \le \eta\;.
\end{align}

\begin{Theorem} \label{th main conv result}Let Assumption \ref{ass X0}, Assumption \ref{assumption:coefficient} (i), (ii), Assumption \ref{as psi function} and \eqref{eq constraint alpha-beta-M 1 intro} hold. If $LT^2$ and $\eta$ are small enough, then there exists $\delta < 1$ such that
\begin{align}\label{eq bound full error}
\cE_{\mathrm{MSE}} \le C_{\rho,\gamma}\left( \delta^P  + 
 \frac{\kappa_{K}}{h \wedge \alpha_{K}} \left(e^{- 2\sqrt{2}\ln(2) \gamma \frac{\alpha_{K} }{\beta_{K}}M^{1-\rho}} 
+  \frac{ \beta_{K}}{\alpha_{K} M^\rho }
\right)  +  h +\cE_{\psi} \right).
\end{align}
for some positive constant $C_{\rho,\gamma}$, where we recall that $\mathcal{E}_{\psi}$ is given by \eqref{eq error approx space brut}.
\end{Theorem}

\begin{Remark}
\begin{enumerate}
\item As expected, the above upper bound is the sum of the error due to the Picard iteration, the error induced by the SGD algorithm, the discrete-time approximation error, recall \eqref{eq control disc error}, and the error $\mathcal{E}_{\psi}$ generated by the approximation in the functional spaces $(\mathscr{V}^y,\mathscr{V}^z_n)$.
\item The smallness condition on $LT^2$ is precisely given in the statement of Proposition \ref{pr main abstract error}. This condition should not come as a surprise since we use Picard iteration. The smallness condition on $\eta$ is not restrictive in practice as the quantity it controls should go to zero to obtain the convergence of the numerical procedure.
\end{enumerate}
\end{Remark}

\noindent To deduce a rate of convergence from \eqref{eq bound full error}, one has to chose the approximation basis functions  $(\psi_y^{k})_{1\leq k \leq K^y}$ and $(\psi_n^{k})_{0\leq n \leq N-1, 1\leq k \leq K^z_n}$ and to set optimally the algorithm's parameters. The choice of the basis function has a dramatic impact on the complexity of the algorithm. In the next section, we work with sparse grid  approximation and we are able to show that the complexity is controlled both theoretically and in practice under Assumption \ref{assumption:coefficient}.

\section{Convergence results for sparse grid approximation}
\label{se main sparse grids}
Both the implemented \emph{direct algorithm}, see Definition \ref{de implemented direct algo}, and the implemented \emph{Picard algorithm}, see Definition \ref{de picard scheme full}, rely on the choice of the approximation spaces $\mathscr{V}^y$ and $\mathscr{V}^z_{n}$, $0 \le n \le N-1$ and the choice of the related basis functions  $(\psi_y^{k})_{1\leq k \leq K^y}$ and $(\psi_n^{k})_{0\leq n \leq N-1, 1\leq k \leq K^z_n}$. The impact is both theoretical, in terms of convergence rate and numerical complexity, and practical in terms of computational time. 
%
%
\noindent We choose here to use sparse grid approximations. This will allow us to obtain interesting numerical complexity results in the setting of Assumption \ref{assumption:coefficient},  see \Cref{thm full complexity}.
We carefully investigate the convergence of the implemented \emph{Picard Algorithm}. It is not the first time that sparse grid approximations are investigated in the context of linear regression. We will use the framework introduced in \cite{bohn2018convergence}. Note however that some restriction in the choice of sparse grid approximations are introduced by Assumption \ref{as psi function}.

\vspace{2mm}
The basis functions are built using elementary bricks that
have a compact support included in the bounded domain 

\begin{align}\label{definition:domain:On}
\cO_n = \prod_{l=1}^d [\mathfrak{a}^n_l, \mathfrak{b}^n_l]  
\quad \mbox{ where } \quad  \mathfrak{a}^n_l < \mathfrak{b}^n_l\; \text{ for } l \in \set{1,\dots,d}.
%
\end{align}

The domain specification strongly depends on the applications under study. We will consider two main cases in this work.
\begin{enumerate}
\item For all $1 \le n \le N-1$, 
\begin{align}\label{de domain fixed}
\cO_n = \prod_{l=1}^d [\mathfrak{a}_l, \mathfrak{b}_l]=:\cO .
\end{align} 
Namely, the coefficients $\mathfrak{a}$ and $\mathfrak{b}$ do not depend on $n$. This will be the case in Section \ref{subse picard periodic} where we consider coefficient functions that are $\cO$-periodic.
\item Alternatively, the coefficient $\mathfrak{a}$ and $\mathfrak{b}$ are functions of the time-step but also of the diffusion coefficients $(b,\sigma)$, recall \eqref{eq forward sde}, and the PDE dimension $d$. Namely
\begin{align}\label{eq general case a-b}
\mathfrak{a}^n:= \mathfrak{a}(t_n,b,\sigma,d) \quad\text{ and }\quad\mathfrak{b}^n:= \mathfrak{b}(t_n,b,\sigma,d). 
\end{align}
\end{enumerate}
However, In both cases the basis functions are obtained by a transformation of the domain $[0,1]^d$ on which we define the primary basis using \emph{sparse grids}. The transformation is defined as follows:
\begin{align}\label{eq de tau}
\tau_n:\mathbb{R}^d \ni x \mapsto \tau_n(x) = \Big(\frac{x_1 - \mathfrak{a}^n_1 }{\mathfrak{b}^n_1- \mathfrak{a}^n_1}, \dots, 
\frac{x_d - \mathfrak{a}^n_d }{\mathfrak{b}^n_d- \mathfrak{a}^n_d}
\Big)^\top \in \mathbb{R}^d\,.
\end{align}
We will introduce two types of basis functions: the first one, based on pre-wavelet basis, follows  from \cite{bohn2018convergence} and the second one, based on hat functions modified at the boundary of the domain, follows from \cite{frommert2010efficient}.

\subsection{Convergence results for the pre-wavelet basis}

\subsubsection{Definition of the pre-wavelet basis}

We describe here the elementary bricks that are used to build the basis functions of the approximation spaces.

For a level $l \in \mathbb{N}$ and an index $i \in \set{0, \dots, 2^l}$, we first consider the family of hat functions given by
\begin{align}
\phi^{l,i}(x) = \phi(2^lx-i) \text{ with }  \phi(x)   =\left\{
\begin{array}{rcl}
1- |x|     &      & {if \; -1<x<1}\\
0       &      & {otherwise\,.}
\end{array} \right.
\end{align}
The univariate pre-wavelet basis functions $\chi^{l,i}:\R\rightarrow \R$ are defined by
\begin{align}\label{de level 0 1 unidim}
\chi^{0,0}  = 1_{[0,1]}\,, \chi^{0,1} = \phi^{0,1}\,, \chi^{1,1} = 2 \phi^{1,1} - 1
\end{align}
and for $l \ge 2$, $i \in I_{l}\setminus \set{1, 2^l-1}$ with $I_l := \set{1 \le i \le 2^l-1\,|\,i \text{ odd}}$
\begin{align}\label{de any level unidim inside}
\chi^{l,i} = 2^{\frac{l}2}\left(\frac1{10} \phi^{l,i-2} - \frac6{10}\phi^{l,i-1} + \phi^{l,i}- \frac6{10}\phi^{l,i+1} + \frac1{10}\phi^{l,i+2} \right).
\end{align}
For the boundary points $i \in \set{1, 2^l-1}$, we set
\begin{align}\label{de any level unidim boundary}
\chi^{l,1} = 2^{\frac{l}2} \left( -\frac65\phi^{l,0} +\frac{11}{10}\phi^{l,1} -\frac35 \phi^{l,2} + \frac1{10}\phi^{l,3}\right) \text{ and } \chi^{l,2^l-1}(x) = \chi^{l,1}(1-x)\,,\, x \in \R \;.  
\end{align}
The multivariate pre-wavelet function on $\R^d$ are obtained by a classical tensor-product approach.
For a multi-index level $\mathbf{l}=(l_1,\dots,l_d)$ and a multi-index position $\mathbf{i}=(i_1,\dots,i_d)$,
\begin{align}\label{de any level multidim tensor prod:def}
\chi^{\mathbf{l},\mathbf{i}}(x) = \prod_{j = 1}^d \chi^{l_j,i_j}(x^j).
\end{align}
In this multivariate case, the index sets are given by
\begin{align}
\mathbf{I}_{\mathbf{l}} = \left\{ \mathbf{i} \in \mathbb{N}^d \, \Big|
\begin{array}{rcl}
0 \le i_j \le 1     &      & {if \; l_j = 0},\\
i_j \in I_{l_j}       &      & {if \; l_j > 0},
\end{array} 
\text{ for all } 1 \le j \le d
\right\}\;.
\end{align}
The hierarchical increment spaces are then defined for $\mathbf{l} \in \mathbb{N}^d$ by
\begin{align*}
\mathscr{W}_{\mathbf{l}} := \mathrm{span} \set{\chi^{\mathbf{l},\mathbf{i}} \,| \, \mathbf{i}\in \mathbf{I}_{\mathbf{l}} }\;.
\end{align*}
The sparse grid space approximation at level $\ell$ is given by
\begin{align}
\mathscr{S}_\ell := \bigoplus_{\mathbf{l} \in \cL_\ell} \mathscr{W}_{\mathbf{l}}, \quad \cL_\ell := \set{\mathbf{l} \in \mathbb{N}^d\,, \zeta_d(\mathbf{l}) \le \ell} \label{sub-grid space}
\end{align}

\noindent with $\zeta_d(\mathbf{0}):=0$ and for \textcolor{black}{$\mathbf{l} \neq \mathbf{0}$}
$$ \zeta_d(\mathbf{l}) = |\mathbf{l}|_1 - d + |\set{j|l_j = 0}| + 1\,,$$
where for a multi-index $\mathbf{l}\in \mathbb{N}^d$ we recall that $|\mathbf{l}|_1 = \sum_{\ell=1}^d l_\ell$ and that $|A|$ is the cardinality of $A$.\\
The key point here is that the dimension of $\mathscr{S}_\ell$ satisfies
\begin{align}\label{eq dim approx space}
\mathrm{dim}(\mathscr{S}_\ell) = \textcolor{black}{ O(2^\ell\ell^{d-1})} \;,
\end{align}

\noindent so that the curse of dimensionality only appears with respect to the level $\ell$, see \cite{feuersanger2010sparse} (and also in the constant related to the notation $O(.)$). The key point now is that the approximation error when using the sparse space is also controlled if the function to be approximated is smooth enough. To this end, for the fixed open domain $(0,1)^d$, we consider the space of function with mixed derivatives $H^{2}_{mix}((0,1)^d)$ (see the section Notation for a precise definition). Then, for any $v \in H^{2}_{mix}((0,1)^d)$, it holds 
\begin{align}\label{eq control sparse basic}
\inf_{\xi \in \mathscr{S}_\ell} \|\xi-v\|^2_{L^2((0,1)^d)}  \le C 2^{-4\ell} \ell^{d-1}\|v\|^2_{H^{2}_{mix}((0,1)^d)}
\end{align}

\noindent for some positive constant $C:=C(d)$. We refer e.g. Theorem 3.25 in \cite{bohn2017error} for a proof of this result. Again, we importantly emphasize that in the above control of the error the curse of dimensionality only appears with respect to the level $\ell$.

\vspace{2mm}

\begin{Remark}\label{re indexation}
The number of basis functions is thus $K = \mathrm{dim}(\mathscr{S}_{\ell})$. We denote by $k:\cC \mapsto \set{1,\dots,K}$ any bijection enumerating $\cC$. We will often slightly abuse the notation and write directly $(\psi_{n}^k)_{1 \le k \le K}$ instead of $(\psi_n^{(\mathbf{l},\mathbf{i})})_{(\mathbf{l},\mathbf{i}) \in \cC}$ to be consistent with the notation introduced in the previous section.
\end{Remark}

\subsubsection{The \emph{Picard Algorithm} in the case of periodic coefficients}

\label{subse picard periodic}

In this section, we work under the setting of Assumption \ref{assumption:coefficient} (iv). To alleviate the notation -- but without loss of generality -- we  assume that the coefficients are $1$-\emph{periodic} in the following sense: for $\lambda = b, \sigma \text{ or } g$
\begin{align}
\lambda(x + q ) = \lambda(x)\,, \text{ for all } (x, q) \in \R^d \times \Z^d\,,
\end{align}
which implies the same property for the value function $u$ and its derivatives.
\\
We thus consider here that $ \cO = [0,1]^d$, recall \eqref{de domain fixed} and  $\tau = I_d$, recall \eqref{eq de tau}.  Here, we are looking for an approximation $\cU^P_M(\cdot)$ of $u(0,\cdot)$ on the whole domain $\cO$, recall \eqref{eq de approx u}. We thus set $\cX_0$ to be uniformly distributed on $(0,1)^d$, which means that Assumption \ref{ass X0}(i) holds true.
\\
For sake of clarity, we summarize the current setting in the following assumption:
\begin{Assumption}\label{ass periodic}
Let Assumption \ref{ass X0}(i) and Assumption \ref{assumption:coefficient}  hold true. Moreover, set
$ \cO = [0,1]^d$ and $\cX_0 \sim \cU((0,1)^d)$.
\end{Assumption}

\vspace{2mm}
To take into account the periodic setting in our approximation, let us first define the $1$-\emph{periodisation} of a compactly supported function $\varphi$ by
\begin{align}
\widecheck{\varphi}(x) := \sum_{  q \in \Z^d }\varphi(x+ q)\;, \text{ for all } x \in \R^d \label{eq de one periodisation}\;.
\end{align}

\noindent The basis functions $\psi$ are then given by $\psi = \widecheck{\chi}$. Namely, for any $0 \le n \le N-1$, for an approximation date $t_n$, we introduce the set of functions
\begin{align}\label{defintion:space:basis:function:translated}
\mathscr{V}^z_n := \set{\xi: \mathbb{R}^d \mapsto \mathbb{R} \, | \, \xi(x) = \widecheck{v}(x), \, \textnormal{ for some } \,v \in \mathscr{S}_{\ell}} \;.
\end{align}
Moreover, at the initial time, the approximation of $u(0,\cdot)$ will also be computed in
\begin{align}\label{defintion:space:basis:function:for:y}
\mathscr{V}^y := \set{\xi: \mathbb{R}^d \mapsto \mathbb{R} \, | \, \xi(x) = \widecheck{v}(x), \, \textnormal{ for some } \,v \in \mathscr{S}_{\ell}} \;.
\end{align}

\begin{Remark}
We could have set an approximation level different for each time step, however we shall not use this possibility in our theoretical or numerical convergence results. We thus simply consider a fixed positive level $\ell$ of approximation, that, obviously, will be chosen later in an optimal way. 
\end{Remark}


\noindent Let also introduce the function 
\begin{align}
 \R^d \ni x \mapsto \widehat{x} \in [0,1)^d
 \end{align}
such that $\widehat{x}$ and $x$ belong to the same equivalence class in $\R^d/\Z^d$. Denoting by $\P_{X_{t_n}}$ the probability measure on $\mathbb{R}^d$ associated to the random vector $X_{t_n}$ given by the Euler-Maruyama scheme \eqref{eq de euler X} taken at time $t_n$ and starting from $\cX_0$ at time 0 and using Lemma \ref{Aronson:two:sided:bounds}, we remark that the boundary of the domain $\cO$ has null $\P_{X_{t_n}}$-measure. We thus deduce
\begin{align}
\widecheck{\psi}(X_{t_n}) = \psi(\widehat{X}_{t_n})\quad\P-a.s.\;,\label{identity:periodic:function:euler:scheme}
\end{align}
and in practice we should work with the latter quantity. 
Namely, we  construct our approximation scheme using:
\begin{align}\label{de eq approx periodic}
{Y}^{\mathfrak{u}}_{0} :=    \sum_{k = 1}^{K^y}  {\psi}^{k}_y(\widehat{X}_{0}) \mathfrak{y}^{k} \quad\text{and}\quad 
 {Z}^{\mathfrak{u}}_{t_n} :=    \sum_{k = 1}^{K^z_n}  {\psi}^{k}_n(\widehat{X}_{t_n}) \mathfrak{z}^{n,k},\; \text{ for } 0 \le n \le N-1\,,
\end{align}
with $\mathfrak{u}=(\mathfrak{y},\mathfrak{z}) \in \R^{K^y}\times\R^{d \bar{K}^z}$.

\textcolor{black}{
Under the current setting of periodic coefficients and sparse grid approximation, we take benefit of the convergence results given in Theorem \ref{th main conv result} to obtain our main theoretical result on the complexity of the \emph{Picard algorithm}. Indeed, the next theorem shows that the curse of dimensionality is tamed by using the sparse grid approximation.
}
\begin{Theorem} \label{thm full complexity}
Let Assumption \ref{ass periodic} hold and assume that $LT^2$ is small enough. For a prescribed $\varepsilon>0$, the complexity $\cC_\varepsilon $, defined in Remark \ref{re num picard algo}, of the full \emph{Picard algorithm} in order to achieve a global error $\cE_{\mathrm{MSE}}$ of order $\varepsilon^2$, recall \eqref{eq de error algorithm},  satisfies
{
$$
\cC_\varepsilon = O_d(\varepsilon^{-\frac52(1+2\iota)}|\log_{2}(\varepsilon)|^{1+ \frac{45+50\iota}{36}(d - 1)})\;
$$
\noindent for any $1<\iota <\frac95$.
}
\end{Theorem}

\noindent The proof of this theorem is given in Section \ref{se con and com analysis} where the algorithm's parameters are optimally set with respect to $\varepsilon$. 

\vspace{2mm}
\paragraph{Periodic example}
We consider here 1-periodic coefficients on $\mathbb{R}^d$. The coefficients of the forward SDE \eqref{eq forward sde} are given by, for $x \in \R^d$,
\begin{align*}
b_i(x) = 0.2 \sin(2\pi x_i)\,,\; \sigma_{i, j}(x) = \frac{1}{ \sqrt{d \pi} } ( 0.25 + 0.1 \cos(2\pi x_i)) \1_{\{i = j\}}\,,\; 
1\le i,j \le d\,.
\end{align*}
The coefficients of the BSDE reads
\begin{align*}
g(x) &= \frac{1}{\pi} \left( \sin\left(2\pi  \sum_{i=1}^d  x_i\right) + \cos\left(2 \pi  \sum_{i=1}^d  x_i\right)  \right) ,  x \in \R^d,  
\\
f(t, x, y, z) &=  2\pi^2 y \sum_{i=1}^d  (\sigma_{i, i}(x))^2 - \sum_{i=1}^d    b_i(x)  \frac{z_i}{ \sigma_{i,i}(x) } + h(t, x),  t\in [0, T], x\in \R^d, y\in \R , z\in \R^d,   
\end{align*}
where $ h(t, x) = 2 \left( \cos(2\pi  \sum_{i=1}^d x_i  + 2\pi (T-t) ) - \sin(2\pi \sum_{i=1}^d  x_i  + 2\pi (T-t) ) \right)  $. The explicit solution is given by 
\[ u(t, x)  = \frac{1}{\pi} \left( \sin\left(2\pi  \sum_{i=1}^d x_i  + 2\pi (T-t) \right) +  \cos\left(2\pi \sum_{i=1}^d  x_i  + 2\pi (T-t) \right) \right) ,  x\in \R^d. \]

\begin{figure}[H]
	\centering
	\includegraphics[width=14cm]{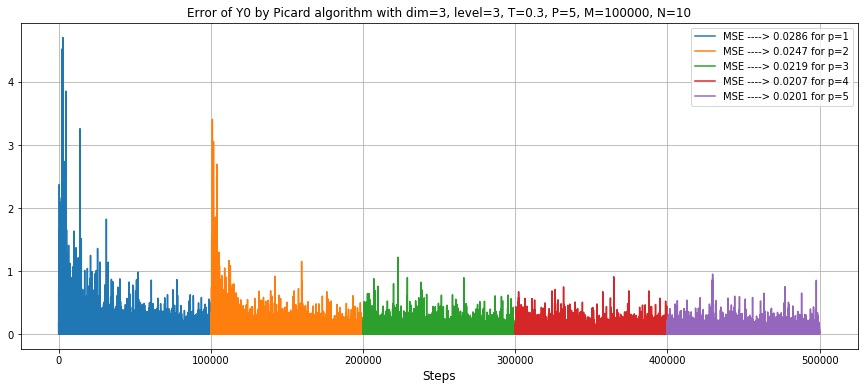}\caption{$m\mapsto |\hat{Y}_0^m - u(0, \cX_0^m)|^2$ for the \emph{Picard algorithm}, $d=3$. The MSE is computed by the mean of the last 10000 steps of each Picard iteration.}\label{fig Y_0_periodic 3d Picard}
\end{figure}

We perform the test for $d=3$ and $M=100000, N=10, T=0.3, level=3$ by \emph{Picard Algorithm} with $P=5$, then there are $K^y = K_n^z = K =225$ basis functions. We obtain a mean square error $\cE_{\mathrm{MSE}}  = 0.0201$ at the 5-th Picard iteration: See \Cref{fig Y_0_periodic 3d Picard} displaying the learning performance. The parameters of the test are shown in \Cref{Parameters periodic} in the appendix.

\subsubsection{Numerical convergence of the \emph{Picard} and \emph{direct Algorithm}}
 \label{subsubse numerics one}
 
We will now investigate numerically the behavior of the \emph{Picard algorithm} and \emph{direct algorithm}
on ``test'' examples that have been already considered in the literature. In particular, this will allow us also to compare our methods to existing methods as the ones investigated in \cite{weinan2017deep,hure2020deep}.

\vspace{2mm}
For this section, we work in the setting of  Assumption \ref{ass X0}(ii). This means that at the initial time, the output $(\hat{y}_0,\hat{z}_0)$ of the  algorithms: $(\hat{y}_0,\hat{z}_0) = (\mathfrak{y}_M,\mathfrak{z}_M^0)$, for the \emph{direct algorithm}, recall Definition \ref{de implemented direct algo}, or $(\hat{y}_0,\hat{z}_0) = (\mathfrak{y}^P_M,(\mathfrak{z}^P_M)^0)$, for the \emph{Picard algorithm}\footnote{\textcolor{black}{Deviating slightly from Definition \ref{de picard scheme full}, we will use for the initialization of the current SGD step, the last value computed at the previous step instead of a random value.
} }, recall Definition \ref{de picard scheme full}, are approximating $(u(0,x_0),\sigma^\top\nabla_xu(0,x_0)) \in \R\times\R^d$. Since, these are one point values, there is no need to introduce basis functions at the initial time and the approximating spaces are just $\mathscr{V}^y = \R$ and $\mathscr{V}^z_0 = \R^d$. Then, for any $1 \le n \le N-1$, for the discrete time $t_n$, we set the approximating space as follows:\begin{align}\label{defintion:space:basis:function:translated:numerics-prewavelet}
\mathscr{V}^z_n := \set{\xi: \mathbb{R}^d \mapsto \mathbb{R} \, | \, \xi(x) = {v}(\tau_n(x)), \, \textnormal{ for some } \,v \in \mathscr{S}_{\ell}} \;,
\end{align}
recall \eqref{eq de tau}. In particular, the basis functions are given by $\psi^k_n(x) = \chi^k(\tau_n(x))$, recall \eqref{de any level multidim tensor prod:def} and Remark \ref{re indexation}.

\vspace{2mm}
We now report more specifically the various algorithms parameters that have been used in practice. The first thing to note is that we are able to obtain good results with a low level of approximation. Indeed, in all our numerical tests, we set the level $\ell = 3$. The \Cref{table L2 sparse grids boundary} below indicates the number of basis functions that have theoretically to be considered when including boundary function.

\begin{table}[H]
	\centering
	\begin{tabular}{|c|c|c|c| }
		\hline
		\diagbox{dimensions}{levels} & $\ell \le 3$ &  $\ell\le 4$  & $\ell\le 5$ \\
		\hline
		d=2 & 49 & 113 & 257  \\
		\hline
		d=3 &225 &593 &1505 \\
		\hline
		d=4 &945 & 2769  & 7681  \\
		\hline
		d=5 &3753 &12033 &36033\\
		\hline
	\end{tabular}
	\caption{The number of functions in the sparse grid approximation with boundary for different dimensions and levels.}
\label{table L2 sparse grids boundary}
\end{table}

Next, we need to define the domain \textcolor{black}{ $\cO_n, 1\le n\le N -1, $} where the approximation will be computed, which depends on the underlying process, recall \eqref{definition:domain:On}-\eqref{eq general case a-b}. We will consider two cases in our simulations, each component of the forward SDE is given by a Brownian motion with drift $\mu$ and volatility $\sigma$: $t \mapsto x_0+ \mu t + \sigma W_t$ or a geometric Brownian motion: $t \mapsto x_0 \exp((\mu-\sigma^2/2)t+\sigma W_t)$.
\begin{enumerate}
\item For the Brownian motion with drift, we set
\begin{align}\label{eq de BM domain}
 \cO_n = x_0 + [\mu t_n - r \sigma \sqrt{t_{n}}, \mu t_n + r \sigma \sqrt{t_{n}}  ]^d, \quad \text{for some } r\in \R^+. 
\end{align}
\item For the geometric Brownian motion, we set
\begin{align}\label{eq de GBM domain}
{ \mathcal{O}_n = [x_0e^{R -  r\sigma \sqrt{t_{n}}},  x_0 e^{R +  r\sigma \sqrt{t_{n}} }], \quad R =  (\mu - \frac{1}{2} \sigma^2 )t_{n} ,  \quad \text{for some } r\in \R^+. }
\end{align}
\end{enumerate}

Finally, a delicate parameter to chose is the the learning rate. Empirically, it was set to:

\begin{align}  
\label{eq de learning rate}
\gamma_{m}(\lambda, t_n, \alpha, \beta_0, \beta_1, m_0) = \frac{\beta_1(\lambda) n+ \beta_0(\lambda)}{ 1+ (m + m_0(\lambda))^{\alpha(\lambda) } }, \quad 1\le m \le M,  \; \lambda \in \{ \mathfrak{y}, \mathfrak{z}^{0, \cdot}  \}_{n=0} \cup  \{\mathfrak{z}^{n, \cdot} \}_{1\le n\le N-1},  
\end{align}
where $\beta_0, \beta_1 \in \mathbb{R}^+, m_0 \in \mathbb{N}^+, \alpha \in \left(\frac12, 1\right] $.

\begin{Remark}
	i) $m_0$ is a suitable positive number to decrease the learning rates  for  avoiding a big jump of the estimated $\lambda$ in the beginning steps of the algorithm. \\
	ii) The parameter {$ r \in \R^+$ is a suitable number to balance the running time and the errors of the algorithms.}\\
	iii) Both $\beta_0$ and $\alpha$ can be used to adjust the converge speed and the variance of the estimated $\lambda$. Suitable parameters make the algorithm more stable, converge faster and reduce the  variance of the estimated $\lambda$.\\
	iv) Usually, we increase the value of $\alpha$  or decrease the value of $(\beta_0, \beta_1)$ gradually to decrease the convergence rate  with the increase of step $p, 1\le p \le P$ for \emph{Picard algorithm}.
\end{Remark}
Concerning the number of steps in the SGD algorithm, we make the following remark.
\begin{Remark}
We used two techniques to  control $M$ in order to  reduce the computational cost: \\
	i) We use $\alpha \in \left(0, \frac12\right]$  which still works well as the SGD algorithm can converge faster. \\
	ii) If $\beta_1(\lambda) \equiv 0$, for $M$ large enough, the algorithm eventually converge, but $ \{ \mathfrak{z}^{n, k} \}_{1\le n\le N-1}^{1\le k \le K}$ convergence becomes slower and slower with the increase of $n$ (the time step).  We  thus choose $\beta_1(\lambda) >0$  in practice to make all $ \{ \mathfrak{z}^{n, k} \}_{1\le n\le N-1}^{1\le k \le K}$ converge altogether with a smaller $M$ (thanks to the learning rates increase with $n$).
\end{Remark}

The remaining parameters are precised in the examples below. We refer also Section to the Appendix for the collection of all algorithm parameters values used in the numerical simulation.

\paragraph{Quadratic model}
First, we consider the quadratic example, whose driver is set to
\begin{align} \label{eq de quad model driver}
f(y, z) = a|z|^2 = a(z_1^2 + z_2^2 + \cdots + z_d^2), \quad  y\in \mathbb{R}, \; z\in \mathbb{R}^d, 
\end{align}
where $a \in \R$ is a constant, and the terminal condition to 
\begin{align} \label{eq de quad model term cond}
g(x) = \log\left(\frac{1+ |x|^2 }{2}\right), \quad x\in \mathbb{R}^d. 
\end{align}
 The explicit solution can be obtained through the Cole-Hopf transformation(see e.g. \cite{chassagneux2016numerical,weinan2017deep}):
$$  y_t = u(t, x) = \left \{
\begin{array}{rl}
	\dfrac{1}{a} \log \esp{ e^{a\cdot g(x + \cW_{T-t})} },  & {a \neq 0} \\
	\esp{g(x+\cW_{T-t})},  & {a = 0}
\end{array} \right.$$
and
$$  z_t^i = \partial_{x_i} u(t, x) = \dfrac{ \esp{\partial_{x_i} g(x+\cW_{T-t}) e^{a\cdot g(x + \cW_{T-t})} }}{\esp{ e^{a\cdot g(x + \cW_{T-t})}}  },  \quad i=1, 2, \cdots, d.$$
Thus, to obtain a numerical reference solution and $95\%$ confidence interval for $y_0$ and $z_0^i, i=1, 2,\cdots , d$, we use classical Monte Carlo estimation of the expectations. 

\begin{figure}[H]
\centering
\begin{minipage}[t]{0.42\textwidth}
\centering
\includegraphics[width=0.8\textwidth]{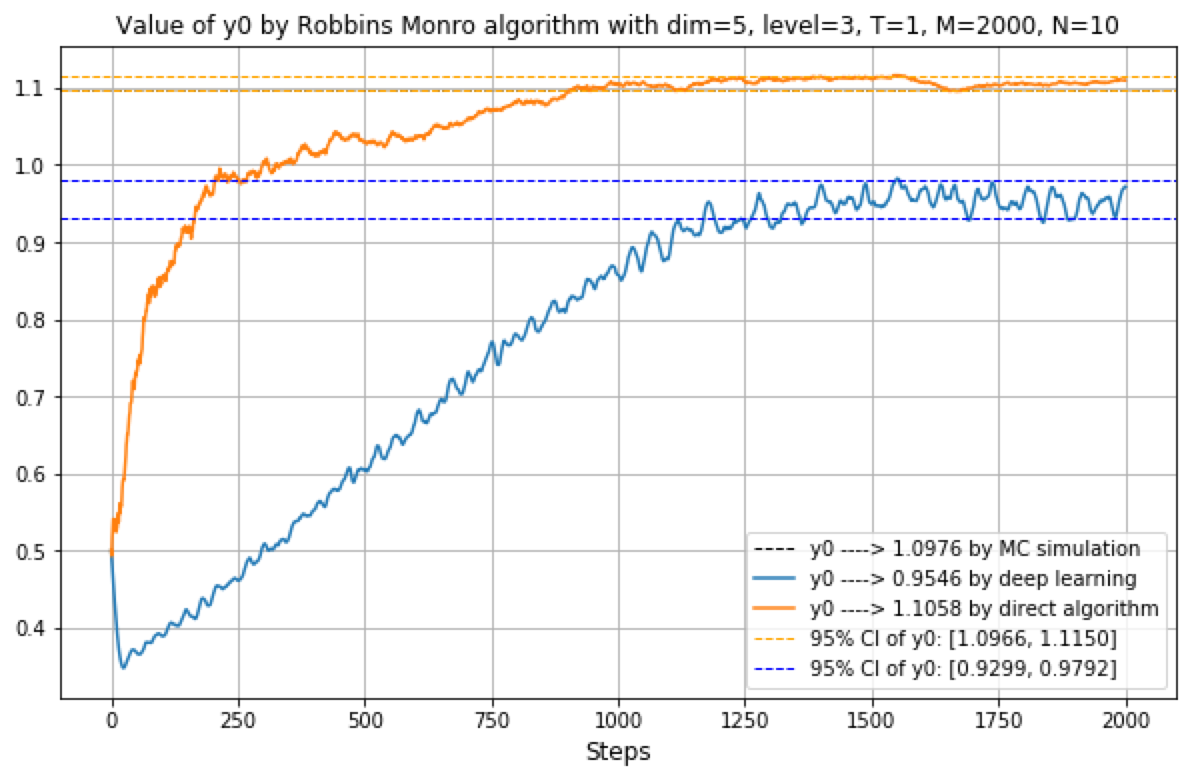}\caption{$\hat{y}_0$ for the quadratic model with d=5 and T=1 by \emph{direct algorithm} and deep learning algorithm.}\label{fig y0 quadratic}
\end{minipage}
\begin{minipage}[t]{0.01\textwidth}
\centering
\; \; \;\;\;\quad
\end{minipage}
\begin{minipage}[t]{0.47\textwidth}
\centering
\includegraphics[width=\textwidth]{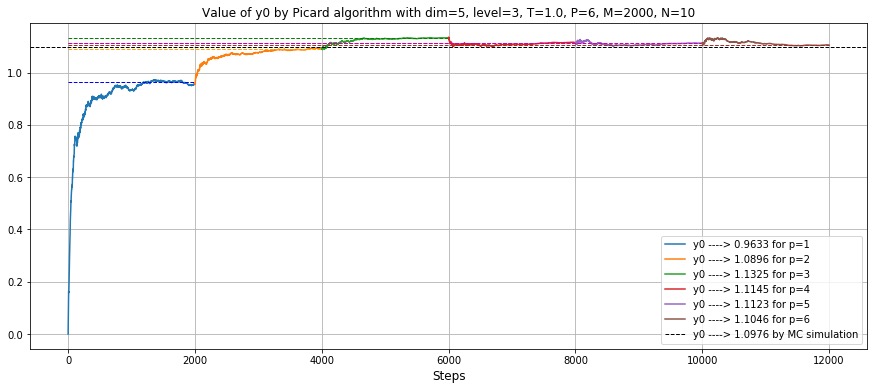}\caption{The value of $\hat{y}_0$ by \emph{Picard algorithm}  with d=5, level=3, T=1, P=6, M=2000.}
	\label{fig y0 quadratic L2 prewavelet picard} 
\end{minipage}
\end{figure}

The underlying diffusion $\cX$ is given by the Brownian motion $\cW$, and the parameters are selected as follows: $a = 1, M = 2000, N = 10,  T = 1$. We compute a reference solution $\bar{y}_0 = 1.0976$ with $95\%$ confidence interval $(1.0943, 1.1009)$ when $d=5$ by  Monte Carlo method using $10^5$ simulation paths. \Cref{fig y0 quadratic} shows the {numerical approximation of $y_0$ and its 95\% confidence interval by the same color line} of the 5-dimensional quadratic model by \emph{direct algorithm} and the \emph{deep learning algorithm} introduced in \cite{weinan2017deep}. The difference of $\hat{y}_0$  between our SGD algorithm and Monte Carlo simulation is less than $10^{-2}$.   It turns out that for this ``low'' dimensional example and with this set of parameter, it is more precise than the \emph{deep BSDEs solver}.

For $d = 5,  M = 2000, N = 10,  P=6, T = 1, a = 1$, \Cref{fig y0 quadratic L2 prewavelet picard} shows that $\hat{y}_0 $ converges for each Picard iteration, and overall  $\hat{y}_0 \rightarrow 1.1046$. We can observe that $\hat{y}_0$ is very close to the reference solution $\bar{y}_0$ when the number of iteration $p$ is greater or equal to $4$.


\paragraph{A financial model} \label{financial model} We now report our numerical results for a model with a financial flavour. The underlying process $\cX$ follows a d-dimensional geometric Brownian motion, for $\mu \in \R $, $\sigma>0$, namely
\[ d\cX_t^i = \cX_t^i (\mu dt + \sigma dW_t^i ), \quad i=1, 2, \cdots d,\quad \cX_0 = x_0 \in (\R_+)^d \;.\]
The driver of the BSDE is given by, for $(y, z)\in \mathbb{R}\times \mathbb{R}^d$,
\begin{align*}
	f(y, z) = - R^l y - \dfrac{\mu - R^l}{\sigma}\sum_{i=1}^d z_i +(R^b - R^l) \max{\left\{ 0, \dfrac{\sum_{i=1}^d z_i}{\sigma} -y  \right\} },  
\end{align*}
and the terminal condition
\begin{align*}
g(x) =  \max\left\{ \left[ \max_{1\le i\le d} x_i  \right]- K_1, 0 \right\} - 2 \max\left\{ \left[ \max_{1\le i\le d} x_i  \right]- K_2, 0 \right\}.
\end{align*}
	Hence, for all $t\in [0, T), x \in \mathbb{R}^d $, it holds that $u(T, x) = g(x)$ and
	\begin{equation} \label{eq 7.2}
		\dfrac{\partial u}{\partial t} + \dfrac{\sigma^2}{2} \sum_{i=1}^d x_i^2 \dfrac{\partial^2 u}{\partial x_i^2} - \min\left\{ R^b (u-\sum_{i=1}^d x_i \dfrac{\partial u}{\partial x_i}  ), R^l (u-\sum_{i=1}^d x_i \dfrac{\partial u}{\partial x_i}  ) \right\}=0\,.
	\end{equation}
This is a typical example of ``non-linear market'' specification, where there are two different interest rates for borrowing and lending money, see Bergman \cite{bergman1995option} and e.g. \cite{weinan2017deep, gobet2005regression, briand2014simulation, crisan2012solving, bender2007forward}, where this example has been used as a test example for numerical methods for BSDEs.

\vspace{2mm}
	
\noindent In our numerical test below, we set the parameters as follows: $N=10, M=6000$, $\mu = 0.06, \sigma =0.2$, $R^l = 0.04, R^b =0.06$, $K_1=110,$ $ K_2=130$,  $T=0.5$ and $x_0=(100, \cdots, 100)$.

\begin{figure}[H]
\centering
\begin{minipage}[t]{0.44\textwidth}
\centering
\includegraphics[width=\textwidth]{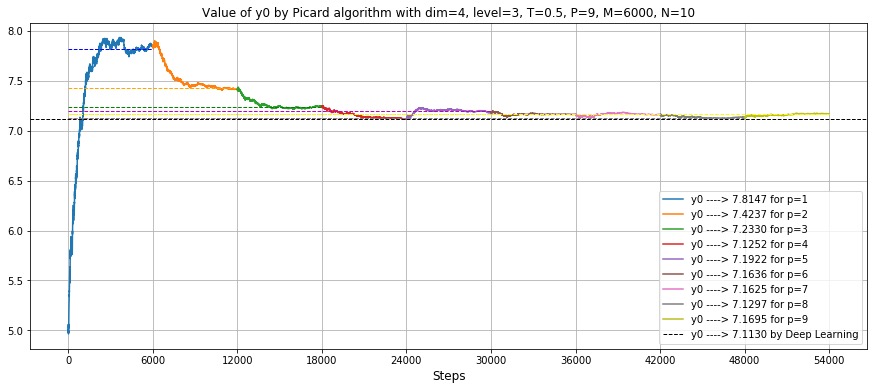}\caption{Approximation $\hat{y}_0$ by \emph{Picard algorithm} when d=4 and T=0.5.} 
\label{L2 sparse prewavelets financial model}
\end{minipage}
\begin{minipage}[t]{0.05\textwidth}
\centering
\; \; \;\;\;\quad
\end{minipage}
\begin{minipage}[t]{0.44\textwidth}
\centering
\includegraphics[width=\textwidth]{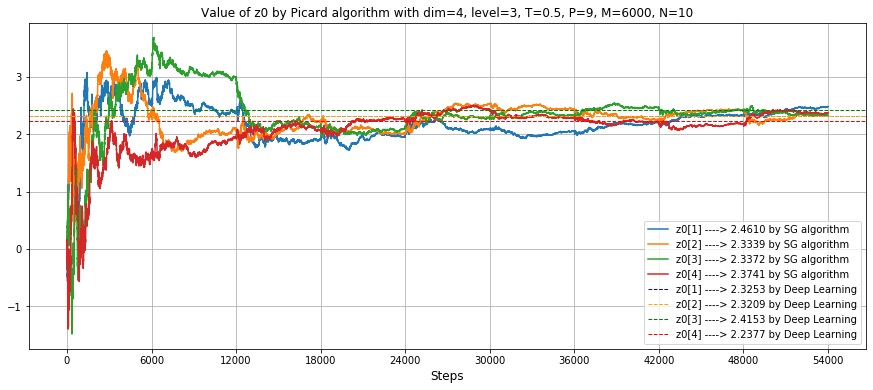}\caption{Approximation $\hat{z}_0$ by \emph{Picard algorithm} when d=4 and T=0.5.} 
\label{L2 sparse prewavelets financial model z0}
\end{minipage}
\end{figure}

\noindent \Cref{compare 2 algos L2 sparse}  compares the results of the \textcolor{black}{ \emph{direct algorithm}} and the \emph{deep learning algorithm} \cite{weinan2017deep}  when $d = 2, 4$. The approximated value of $y_0$ obtained by the two methods are very close. \Cref{L2 sparse prewavelets financial model} 
and \Cref{L2 sparse prewavelets financial model z0}
 show the performance of the \emph{Picard algorithm} with parameters $d=4, P=9$,  \textcolor{black}{and $\hat{y}_0$ converges to $7.1695$ at the last step.}

\begin{table}[htp]
	\centering
	\begin{tabular}{|c|c|c|c|c| }
		\hline
		\multicolumn{1}{|c|}{\multirow{2}*{dimensions}}& \multicolumn{2}{c|}{SGD algo with Sparse grids} &\multicolumn{2}{c|}{Deep learning scheme \cite{weinan2017deep}  }\\
		\cline{2-5}
		\multicolumn{1}{|c|}{}&$y_0$ &  95\% CI of $y_0$  &  $y_0$  &  95\% CI of $y_0$ \\
		\hline
		d=2 & 4.3332  &  [4.2921, 4.3743]    & 4.3516 & [4.3420, 4.3612]   \\
		\hline
		d=4 &  7.0960 & [7.0432,  7.1487]   & 7.1130 & [7.0649,  7.1611] \\
		\hline
	\end{tabular}
	\caption{Comparison of the \emph{direct algorithm} and the \emph{deep learning algorithm} for the financial model.}
	\label{compare 2 algos L2 sparse}
\end{table}

\subsubsection{Limits of the \emph{Picard Algorithm}}

We now illustrate on a numerical example that the smallness assumption may be necessary to obtain the convergence of the \emph{Picard Algorithm}.
To this end, we consider the following model. For a given $a\in \R$, the BSDE driver is given by
\[	f(y, z) := \arctan(ay) + \sum_{j=1}^{d} z_j, \quad (y,z) \in \mathbb{R}\times \mathbb{R}^d, \]
and the terminal condition
\[ g(x) := \dfrac{e^{1 + \1 \cdot x}}{1+e^{1 + \1 \cdot x}}, \quad x\in \mathbb{R}^d .\]
The underlying process $\cX$ is simply equal to the Brownian motion $\cW$, namely $b=0$, $\sigma = I_d$. We set the terminal time $T=1$ and the dimension $d=2$. 

\noindent We study numerically the above model for different value of $a$, which controls the Lipschitz constant of $f$, in the case of the \emph{Picard Algorithm}. The value obtained are compared to the ones obtained by two other methods: a multistep scheme in \cite{chassagneux2014linear}  and the \emph{deep BSDEs solver} of \cite{weinan2017deep}. The values obtained by these two methods are considered to be close to the true solution.

\begin{figure}[H]
\centering
\begin{minipage}[t]{0.44\textwidth}
\centering
\includegraphics[width=\textwidth]{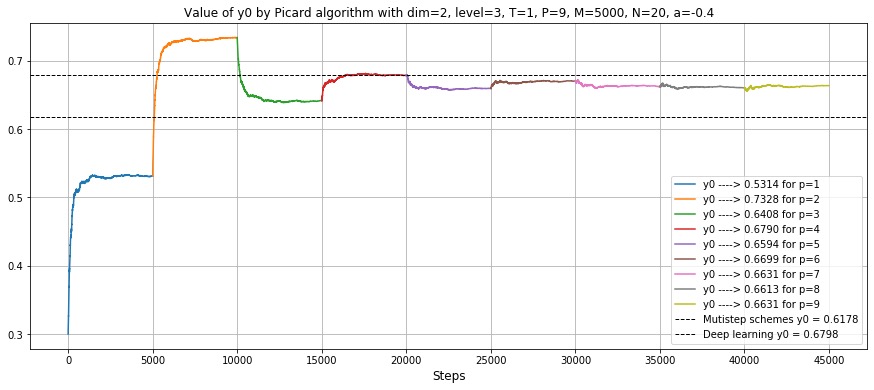}\caption{The value of $\hat{y}_0$ by  \emph{Picard algorithm} with d=2, level=3, T=1, P=9, M=5000, a=-0.4.}
	\label{fig a equal to minus 0.4 }
\end{minipage}
\begin{minipage}[t]{0.05\textwidth}
\centering
\; \; \;\;\;\quad
\end{minipage}
\begin{minipage}[t]{0.44\textwidth}
\centering
\includegraphics[width=\textwidth]{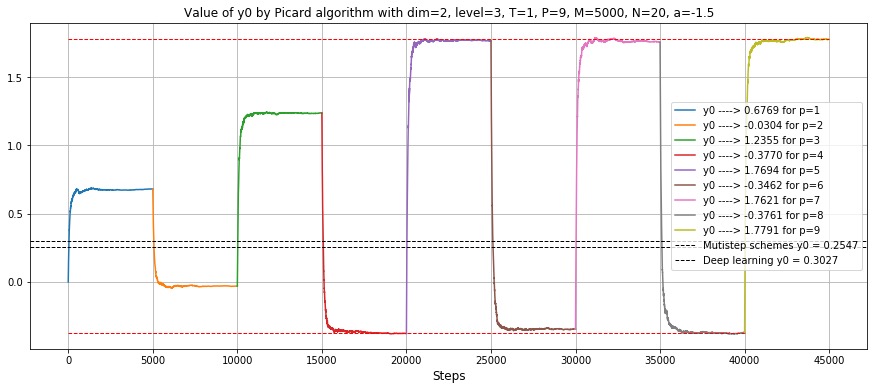}\caption{The value of $\hat{y}_0$ by  \emph{Picard algorithm} with d=2, level=3, T=1, P=9, M=5000, a=-1.5}
	\label{fig a equal to minus 1.5 }
\end{minipage}
\end{figure}
\begin{figure}[H]
	\centering
	\includegraphics[width=11.cm]{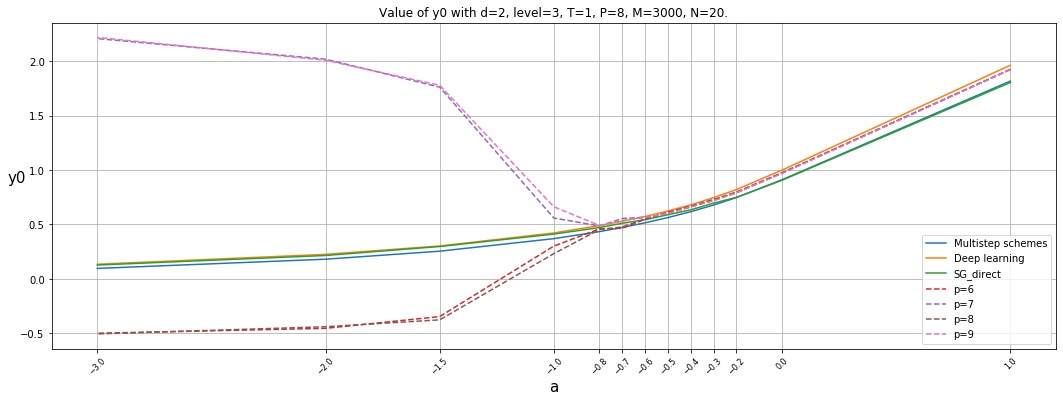}\caption{The value of $\hat{y}_0$ by \emph{Picard algorithm}, \emph{direct algorithm}  and deep learning method  with d=2, level=3, T=1, P=9, M=5000. The last four steps are shown for the \emph{Picard Algorithm} illustrating the bifurcation phenomenon. Note that the \emph{direct algorithm} does not exhibit such behaviour.}
	\label{fig compare 3 methods}
\end{figure}

\noindent When $a= -0.4$,  \Cref{fig a equal to minus 0.4 } shows that $ \hat{y}_0 $ converges. However, this is not the case anymore when $a = -1.5$, see \Cref{fig a equal to minus 1.5 }, as $\hat{y}_0$ oscillates between two values.
Actually, we see on \Cref{fig compare 3 methods} that a bifurcation occurs for the \emph{Picard Algorithm} around $a=-0.8$.

\subsection{Numerical results with the modified hat functions basis}

In the previous section, using the pre-wavelet basis, we were able to establish a theoretical upper-bound on the global complexity for the \emph{Picard algorithm} and to show that both the \emph{Picard algorithm} and \emph{direct algorithm} converge in practice too. However, the number of basis functions, even though we use a sparse approximation, is still quite important \textcolor{black}{which prevents us from dealing effectively with high-dimensional PDE}. In particular, the number of basis functions used to capture what happens on the boundary of the domain is large. In this section, we use, the so-called ``modified hat functions'' that allows to get rid of the boundary basis. 
\subsubsection{Definition of the basis functions}
\label{subsec Numerics L2 hat}
{\color{black}
The modified hat functions are defined by the following method (which corresponds to equation (2.16) in \cite{frommert2010efficient}),
\begin{equation} \label{eq 1.7}
	 \tilde{\phi}^{l,i}(x) :=\left\{
\begin{array}{ll}
 1 
     & {\textnormal{if} \quad l=1\wedge i =1}\\
\left\{
\begin{array}{ll} 
 1 - 2^{l -1} \cdot x  & \qquad\quad {\textnormal{if} \quad x\in [0 , 2h_{l} ] }\\
 0 & \qquad\quad \textnormal{otherwise}
\end{array} \right\}
        & {\textnormal{if} \quad  l>1\wedge i =1}\\
\left\{
\begin{array}{ll} 
  2^{l-1} \cdot x  + (1-i)/2 &  {\textnormal{if} \quad  x\in [1 - 2 h_{l}, 1] }\\
 0 &  {\textnormal{otherwise}}
\end{array} \right\}
          &  {\textnormal{if} \quad l>1\wedge i = 2^{l} - 1}\\
\phi^{l,i}(  x)         & {\textnormal{otherwise}\,,}
\end{array} \right.
\end{equation}

The multivariate hat function on $\R^d$ are obtained by a classical tensor-product approach. For these basis functions, we can remove the points on the boundary of the space so that all the components $l_j, j=1, \cdots, d$, are positive for a multi-index level $\mathbf{l}=(l_1,\dots,l_d)$ and a multi-index position $\mathbf{i}=(i_1,\dots,i_d)$,
\begin{align}\label{de any level multidim tensor prod}
\tilde{\phi}^{\mathbf{l},\mathbf{i}}(x) = \prod_{j = 1}^d \tilde{\phi}^{l_j,i_j}(x^j).
\end{align}
In this multivariate case,  the index set are given by
\begin{align}
\mathbf{I}_{\mathbf{l}} = \left\{ i \in \mathbb{N}^d \,|\,
i_j \in I_{l_j}    
\text{ for all } 1 \le j \le d
\right\}\;.
\end{align}
}

\Cref{table L2 sparse grids without boundary} shows the number of points in the  sparse grids without boundary. In particular, we observe that it is much less than sparse grids with boundary for the same dimensions and levels, recall \Cref{table L2 sparse grids boundary}.

\begin{table}[h]
	\centering
	\begin{tabular}{|c|c|c|c| }
		\hline
		\diagbox{dimensions}{levels} & $l\le 3$ &  $l\le 4$  & $l\le 5$ \\
		\hline
		d=2 & 17 & 49  & 129  \\
		\hline
		d=4 & 49  & 209  & 769  \\
		\hline
		d=5 & 71 & 351 & 1471\\
		\hline
		d=10 &  241 & 2001  & 13441 \\
		\hline
		d=20 & 881  & 13201 & 154881 \\
		\hline
		d=25 & 1351 & 24751 & 352351 \\
		\hline
		d=50 & 5201 & 182001 & 4867201 \\
		\hline
		d=100 & 20401 & 1394001 & \diagbox{}{}\\
		\hline
	\end{tabular}
	\caption{The number of points in the sparse grid approximation without boundary functions for different dimensions and levels.}
\label{table L2 sparse grids without boundary}
\end{table}

\begin{figure}[H]
\centering
\begin{minipage}[t]{0.42\textwidth}
\centering
\includegraphics[width=6cm]{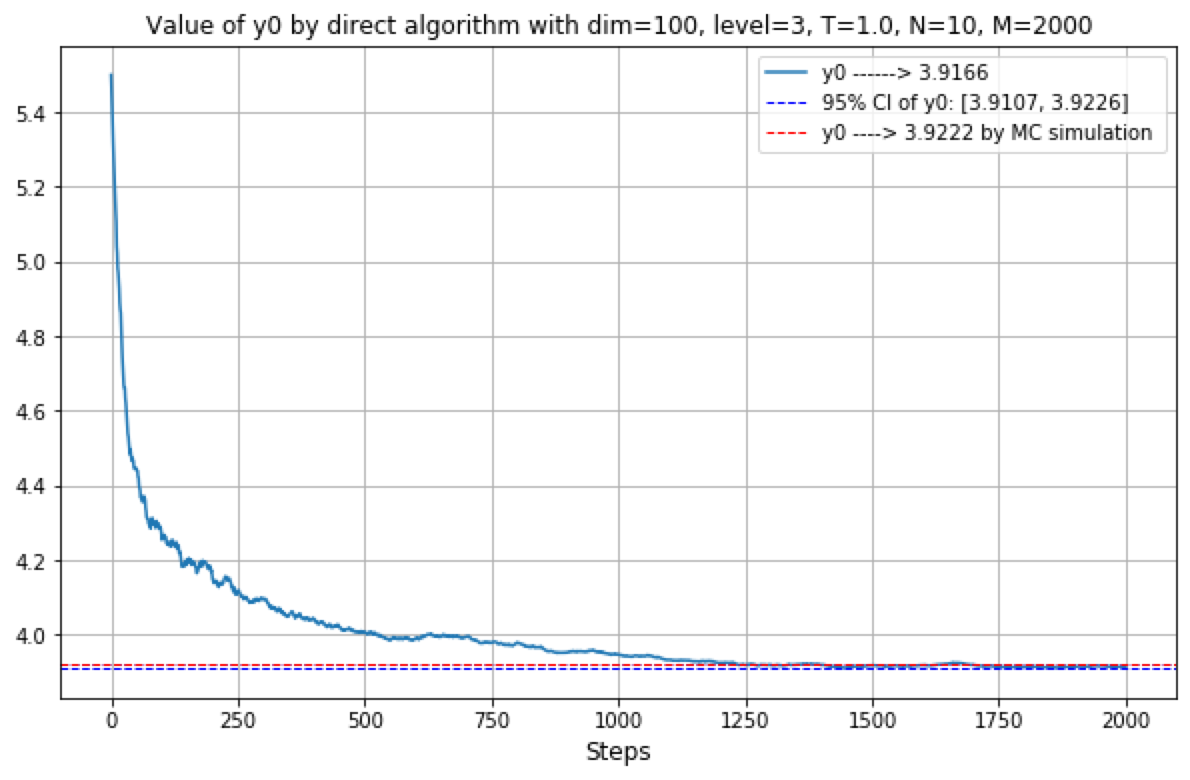}\caption{$\hat{y}_0$ for the quadratic model with d=100 and T=1}\label{fig y0 quadratic L2 hat 100}
\end{minipage}
\begin{minipage}[t]{0.1\textwidth}
\centering
\; \; \;\;\;\quad
\end{minipage}
\begin{minipage}[t]{0.42\textwidth}
\centering
\includegraphics[width=6cm]{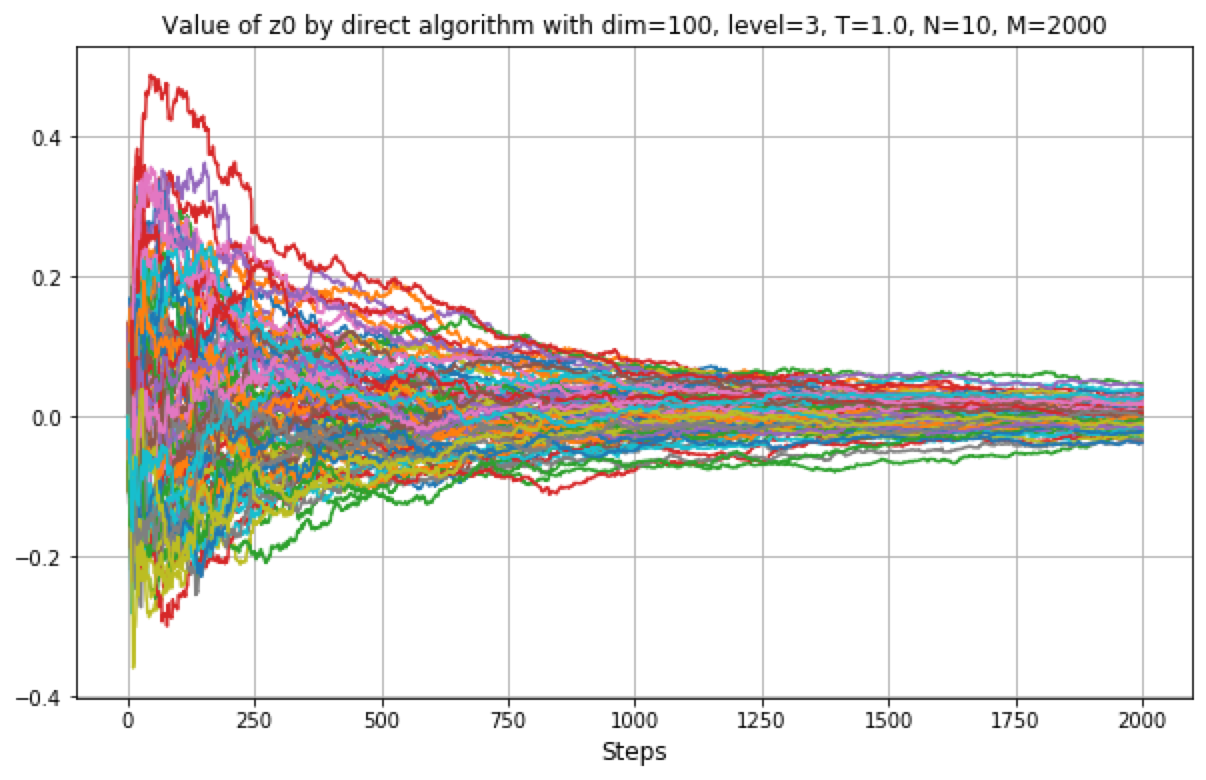}\caption{$\hat{z}_0^i$ for the quadratic model with d=100 and T=1}\label{fig z0 quadratic L2 hat 100}
\end{minipage}
\end{figure}

\begin{figure}[H]
	\centering
	\includegraphics[width =0.8 \textwidth]{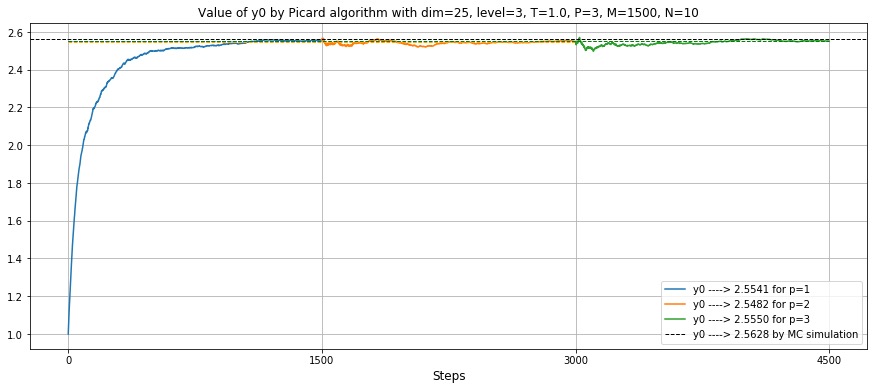}\caption{$\hat{y}_0$ for the quadratic model by \emph{Picard algorithm} with d=25, T=1, P=3, N=10, M=1500}\label{fig y0 quadratic L2 hat picard}
\end{figure}

\paragraph{The quadratic model} We come back to the quadratic model introduced in \eqref{eq de quad model driver}-\eqref{eq de quad model term cond}.
In this setting, we can test the 100-dimensional version of this model.
Let $M=2000, N=10, T=1, a=1$, the convergence of $\hat{y}_0$ and $\hat{z}_0$, when using the \emph{direct algorithm}, is shown in \Cref{fig y0 quadratic L2 hat 100} and \Cref{fig z0 quadratic L2 hat 100}: 3819 seconds were spent on this test. The error for $\hat{y}_0$ appears to be less than 0.01. For $Z$, the true solution is \textcolor{black}{$ \bar{z}_0^i = \frac{2\cW_0^i}{\esp{1+|\cW_T|^2}} = 0, \forall i =1, ..., d$}. The gain in computational time is important in comparison with the pre-wavelet specification of the last section. Not only less basis functions are used, but
one should also note that the computational cost of a hat function is less than a pre-wavelet function up to a factor 5.
We do also test the \emph{Picard algorithm} in  25-dimensional setting. We set $M=1500, N=10, T=1, a=1, P= 3$ and get $\hat{y}_0 \approx 2.5481 $ quite quickly, see \Cref{fig y0 quadratic L2 hat picard} (all the initial value of $ \mathfrak{z}^{n,k}, 1\le n\le N-1, 1\le k\le K^z_n $ are set to $0$ in this test).


%
%
%
%
%

\paragraph{The financial model of  \eqref{eq 7.2}} 
	
\begin{figure}[H]
\centering
\begin{minipage}[t]{0.42\textwidth}
\centering
\includegraphics[width=0.8\textwidth]{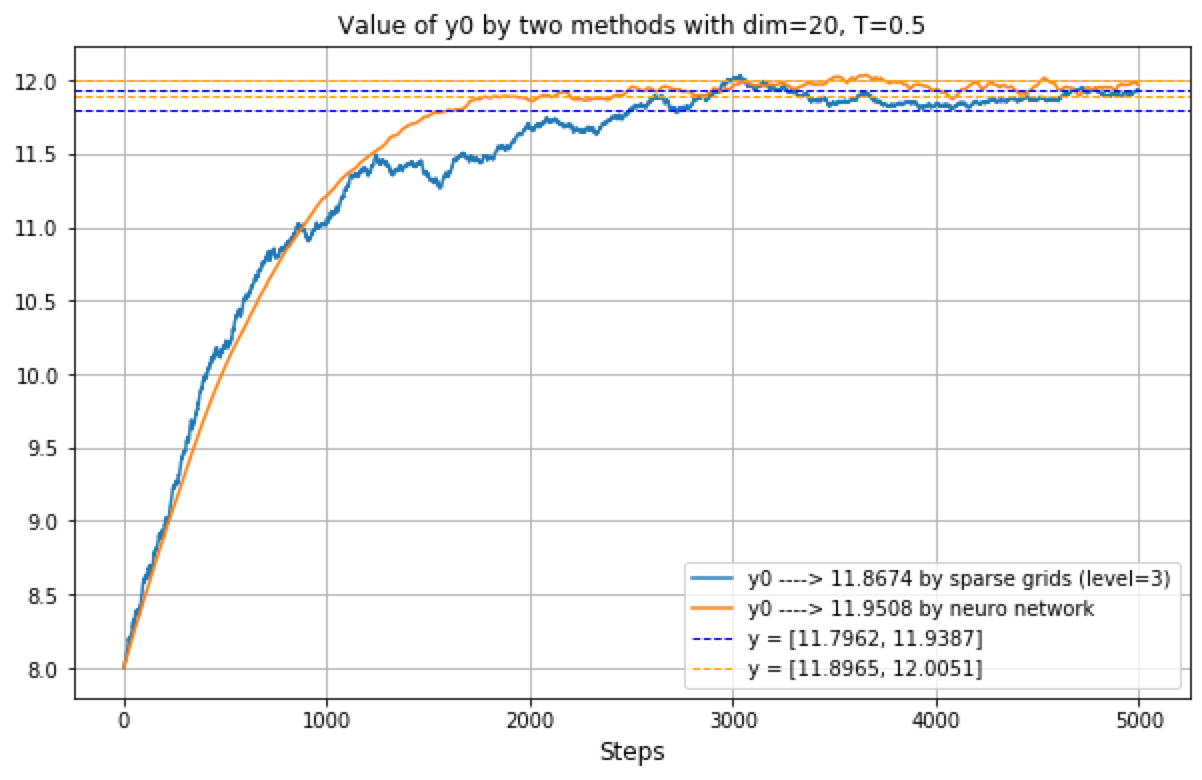}\caption{$\hat{y}_0$  by \emph{direct algorithm} for d=20 and T=0.5.} 
\label{fig financial model 1 L2}
\end{minipage}
\begin{minipage}[t]{0.01\textwidth}
\centering
\quad
\end{minipage}
\begin{minipage}[t]{0.47\textwidth}
\centering
\includegraphics[width=\textwidth]{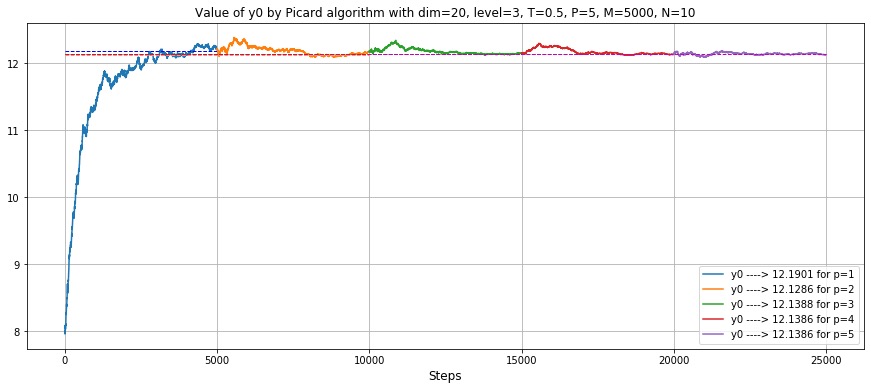}\caption{$\hat{y}_0$   by \emph{Picard algorithm} for d=20, T=0.5, P=4, M=5000, N=10. }
\label{L2 sparse hat financial model}
\end{minipage}
\end{figure}

	For the \emph{direct  algorithm} in this example, we set the parameters $N=10, M=5000$, $\mu = 0.06, \sigma =0.2$, $R^l = 0.04, R^b =0.06$, $K_1=110,$ $ K_2=130$,  $T=0.5$. 
	\Cref{compare 2 algos L2} compares the results for the  \emph{direct  algorithm} and \emph{deep BSDEs solver}: the approximated value for $y_0$ obtained by  the two methods are very close. The running time for the \emph{deep BSDEs solver} shows almost no increase up to $d\le 25$. For the   \emph{direct  algorithm}, it does increase with the dimensions but it stays reasonable. Actually, it is even competitive 
 when $d\le 25$\footnote{ \textcolor{black}{The numerical experiments were realised by C++ 17 on a MacBook Pro 6-core Intel Core i7, using only one core and compiling with optimisation flag `-O3' in gcc.}
 The \emph{deep BSDEs solver} \cite{weinan2017deep}, using \emph{Tensorflow}, spends most of the time to build the graph for the NN and initialize the variables  when the dimension is small then the learning phase is quick. On the contrary, our algorithm builds the approximation grid space quite efficiently (less than 1 second when $d\le 100, level\le 3$) and then the runtime is spent on the the SG algorithm.}.  \Cref{fig financial model 1 L2} shows the performance of the \emph{direct algorithm} when $d= 20$. Finally in \Cref{L2 sparse hat financial model}, the \emph{Picard algorithm} is tested on this model with parameter $d= 20 , P=4$, and $\hat{y}_0 \approx 12.1386$ at the last iteration. 

\begin{table}[htp]
	\centering
	\resizebox{\textwidth}{20mm}{
	\begin{tabular}{|c|c|c|c|c|c|c| }
		\hline
		\multicolumn{1}{|c|}{\multirow{2}*{dimensions}}& \multicolumn{3}{c|}{direct SGD algo with Sparse grids} &\multicolumn{3}{c|}{Deep BSDEs solver\cite{weinan2017deep}  }\\
		\cline{2-7}
		\multicolumn{1}{|c|}{}&$y_0$ &  95\% CI of $y_0$  & time  &  $y_0$  &  95\% CI of $y_0$ & time \\
		\hline
		d=5 & 8.0966  &  [8.0226, 8.1705]  & 3 s  & 8.1010 & [8.0747, 8.1273] & 115 s  \\
		\hline
		d=10 & 10.9865  & [10.9224, 11.0506]&  12 s  & 10.9216 & [10.8944, 10.9489]& 120 s \\
		\hline
		d=15 & 11.848   & [11.7853, 11.9107] &  33 s & 11.8226 & [11.7750, 11.8702] & 122 s\\
		\hline
		d=20 &   11.8674 & [11.7962, 11.9387] &   61 s &  11.9508 & [11.8965, 12.0051] & 127 s\\
		\hline
		d=25 & 11.7801    & [11.6467, 11.9135]  &   130 s & 11.6416  & [11.5316, 11.7517] &  132 s\\
		\hline
	\end{tabular}}
	\caption{Comparison of the \emph{direct algorithm} and the \emph{deep learning algorithm}.}
	\label{compare 2 algos L2}
\end{table}

\paragraph{A challenging example} 
We now consider a model with an unbounded and complex structure solution, which has been analyzed in \cite{hure2020deep}. The value function in this case is given by:
 	 \begin{equation}
 	u(t, x) = \dfrac{T-t}{d}\sum_{i=1}^{d}( \sin(x_i) \mathbbm{1}_{\{x_i < 0 \}} + x_i  \mathbbm{1}_{\{x_i \geq 0 \}} ) +\cos\left(\sum_{i=1}^d i x_i\right), \quad x\in \mathbb{R}^d .
 \end{equation}
 It corresponds to a BSDE, with underlying process given by
 $\cX_t = x + \frac{1}{\sqrt{d}} \mathrm{I}_d \cW_t $, and $x_0=0.5\mathbbm{1}_d $ and
driver and terminal condition given respectively by
\begin{equation*}
	\begin{split}
		f(t, x, y, z) =& \left(1+(T-t) ( \frac{1}{2d} - C) \right)  A(x)  + (1 - (T-t)C)B(x) + C y, \\
		 =&  \left(1+ \frac{T-t}{2d}  \right)  A(x)  +B(x) + C \cos\left(\sum_{i=1}^d i x_i\right) ,  \quad x\in \mathbb{R}^d, \;y\in \mathbb{R}, \; z \in \mathbb{R}^d,\\
		 g(x) =& u(T, x) =\cos\left(\sum_{i=1}^d i x_i\right),  \quad x\in \mathbb{R}^d,
	\end{split}
\end{equation*}
where
	 \[A(x)= \frac{1}{d} \sum_{i=1}^d  \sin(x_i) \mathbbm{1}_{\{x_i < 0 \}} , B(x)=  \frac{1}{d} \sum_{i=1}^d x_i  \mathbbm{1}_{\{x_i \geq 0 \}}, C = \frac{(d+1)(2d+1)}{12}.\]
	 
\begin{figure}[H]
\centering
\begin{minipage}[t]{0.42\textwidth}
\centering
\includegraphics[width=7cm]{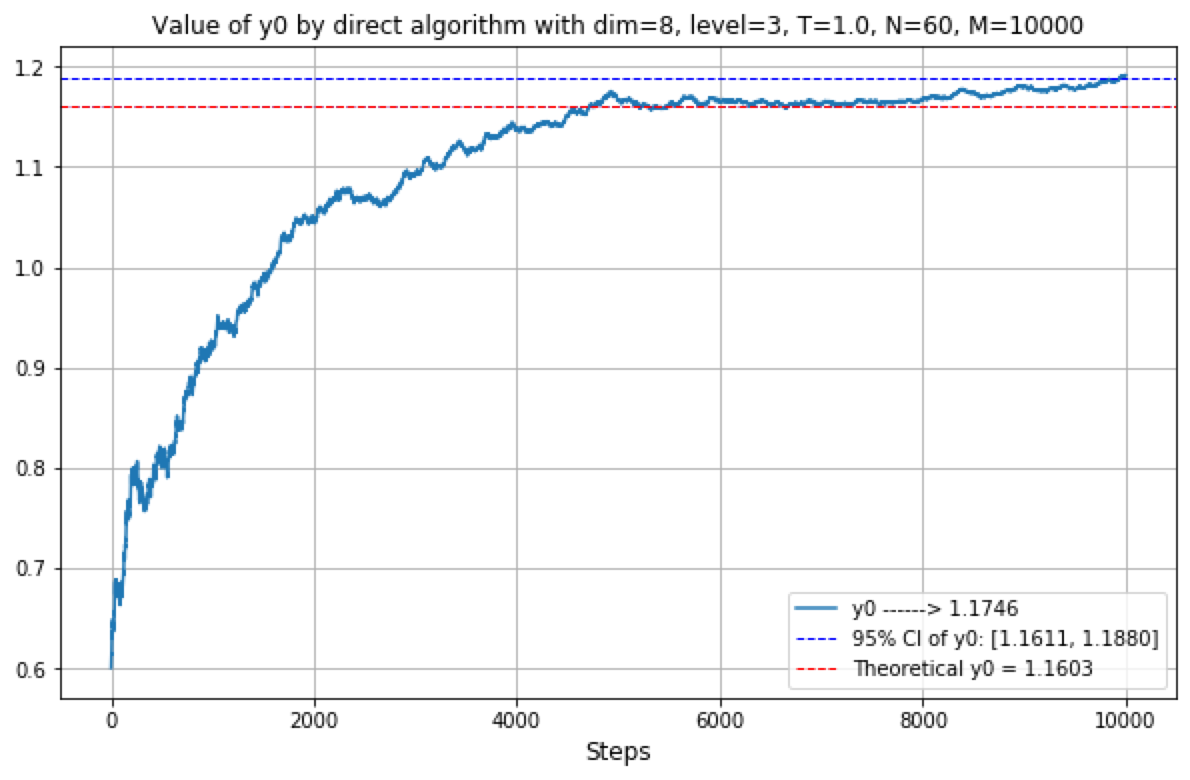} \caption{$\hat{y}_0 \rightarrow 1.1745 $ by direct SGD algorithm  when d=8, N=60, M=10000.}\label{fig come model direct 8}
\end{minipage}
\begin{minipage}[t]{0.1\textwidth}
\centering
\quad
\end{minipage}
\begin{minipage}[t]{0.42\textwidth}
\centering
\includegraphics[width=7cm]{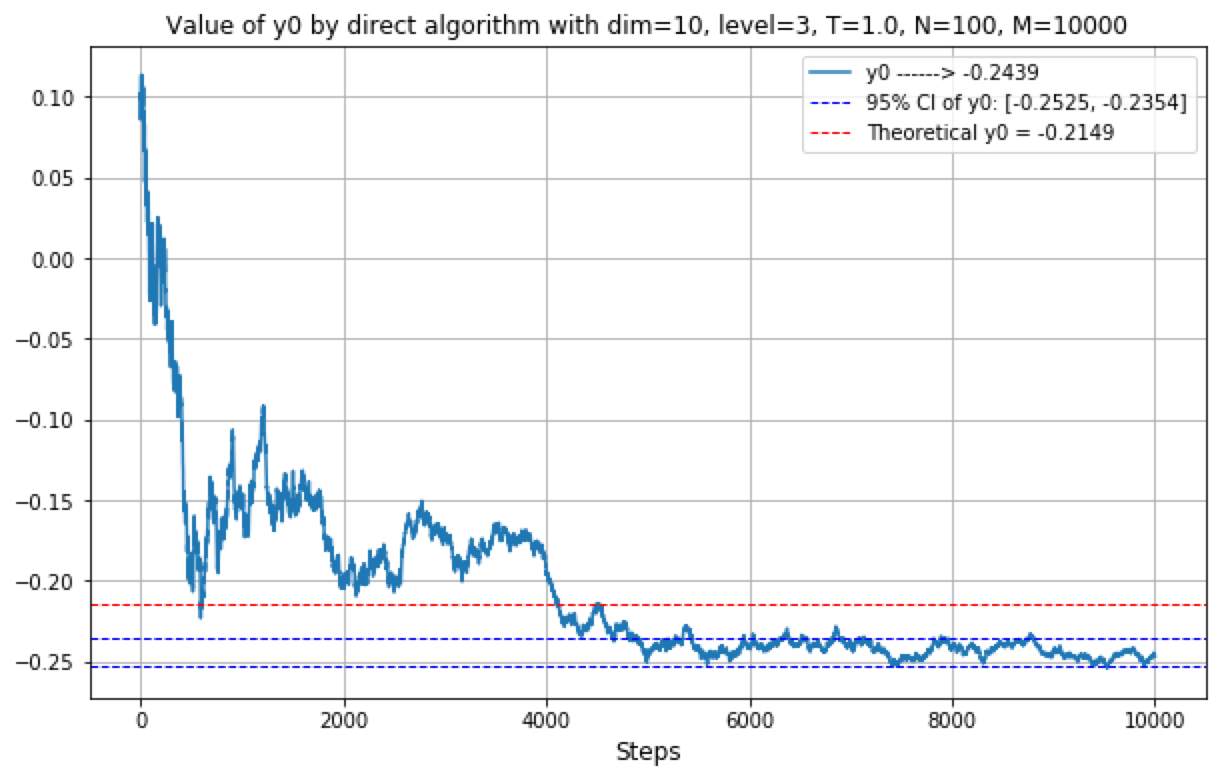} \caption{$\hat{y}_0 \rightarrow -0.2439 $ by direct algorithm  when d=10, N=100, M=10000.}\label{fig come model direct 10}
\end{minipage}
\end{figure}

\begin{figure}[H]
	\centering
	\includegraphics[width = 14cm]{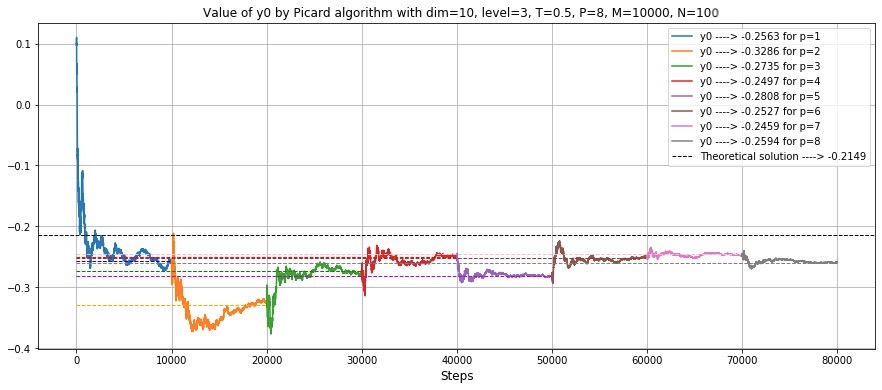} \caption{$\hat{y}_0 \rightarrow -0.2594 $ by picard SGD algorithm  when d=10, P=8, N=100, M=10000.}\label{fig come model picard 10}
\end{figure}

In \Cref{compare 4 algos}, we compare the approximation of $y_0$ by using five different algorithms to the theoretical solution. When the dimension $d \le 2$, all the algorithms perform well. However, \textcolor{black}{as already mentioned in \cite{hure2020deep}} the \emph{deep learning algorithm} \cite{han2017overcoming}  fails when $d\ge 3$ (no matter the chosen initial learning rate and the activation function for the hidden layers, among the tanh, ELU, ReLu and sigmoid ones; besides, taking 3 or 4 hidden layers does not improve the results.) The two deep learning schemes of \cite{hure2020deep}  and our  algorithms with sparse grids still works well when $d \le 8$. \Cref{fig come model direct 8} shows the performance of the \emph{direct algorithm} , $\hat{y}_0$  converges to $1.1745$ when $d=8$, it is close to the theoretical solution $1.1603$, and the 95\% confidence interval of $\hat{y}_0$ is $ [1.1611, 1.1881] $. When the dimension $d=10$,  all the algorithms failed at providing correct estimates of the solution as shown in the table, but the errors of our algorithms appear to be smaller than the errors of deep learning methods. \Cref{fig come model direct 10} and \Cref{fig come model picard 10} show the performance of the \emph{direct algorithm}, \emph{Picard algorithm} respectively. \\

	 
\begin{table}[htp]
	\centering
	\resizebox{\textwidth}{23mm}{
	\begin{tabular}{|c|c|c|c|c|c|c| } 
		\hline
		\multicolumn{1}{|c|}{\multirow{3}*{dimensions}}&  \multirow{3}*{\tabincell{c}{Theoretical\\ solution} }& \multicolumn{2}{c|}{ \multirow{2}*{\tabincell{c}{SGD algo with $L^2$ sparse\\ grids and hat functions} } }& \multicolumn{2}{c|}{ \multirow{2}*{\tabincell{c}{DL scheme \\of HPW \cite{hure2020deep}}} } & \multirow{3}*{\tabincell{c}{DL scheme \\of HJE \cite{han2017overcoming} } } \\
		   && \multicolumn{2}{c| }{}    & \multicolumn{2}{c| }{}    & \\
		 \cline{3-6}
		& &  \emph{direct\; algo} & \emph{Picard algo} &  DBDP1 &  DBDP2 &\\
		\hline
		d=1 & 1.3776  & 1.3790    & 1.3825   & 1.3720 & 1.3736 & 1.3724  \\
		\hline
		d=2 & 0.5707 & 0.5795 & 0.5794 & 0.5715 & 0.5708 & 0.5715  \\
		\hline
		d=5 & 0.8466  &  0.8734  & 0.8606 & 0.8666  & 0.8365 & NC  \\
		\hline
		d=8 & 1.1603  &  1.1745   &  1.1801  & 1.1694  & 1.0758 & NC  \\
		\hline
		d=10 & -0.2149 &  -0.2439 &  -0.2594 & -0.3105 & -0.3961 & NC  \\
		\hline
	\end{tabular}
	} 
	\caption{The approach value of $\hat{y}_0$ of different methods when $T=1$.}
	\label{compare 4 algos}
\end{table}

%
%
%
%
%
%

\section{Study of the discrete optimization problems}\label{study:optimization:problem}

In this section, we study theoretically the \emph{direct algorithm} and the \emph{Picard algorithm} in order to prove the results stated in Section \ref{Description of the numerical methods} and \ref{subse picard periodic}. We first obtain forward and backward estimates on perturbed BSDEs. This allows in particular to derive the wellposedness of the \emph{direct algorithm}. However, most of the work concentrates on the \emph{Picard algorithm} in Section \ref{Study of the Picard algorithm}. A careful study of the iterated SGD algorithms allows to prove the convergence announced in Theorem \ref{th main conv result}. Finally, we prove Theorem \ref{thm full complexity} concerning the complexity of the method in the case of periodic coefficients and using pre-wavelet basis.

\subsection{Preliminary estimates}
In this subsection, we prove general technical estimates for the backward component of a BSDE that will be used in the proof of the convergence of the numerical methods under study.
We will essentially compare two processes with dynamics given by \eqref{eq de bar Y:continuous:time} but taken at two different starting points and controlled processes.
\\
The first process, denoted by $V$, is a scheme built with a random driver $F$ satisfying: 
\begin{Assumption}\label{ass random driver}
\begin{enumerate}
\item For all $(y,z)\in \R\times\R^d$, $F(\cdot, y, z)$ is progressively measurable.
\item There exists some deterministic constant $C \ge 0$ such that for any $t \in [0, T]$ and any $(y,y',z,z') \in \R^2\times\R^{2d}$
$$|F(t, y ,z)-F(t,y',z')| \le C\left( |y-y'| + |z-z'| \right).$$
\end{enumerate}
\end{Assumption}

\vspace{4mm}
\noindent For $Z \in \cS^2_d$ and $\zeta \in \cL^2(\cF_0)$, we thus define
\begin{align}\label{eq de random driver scheme}
V^{\zeta,Z}_t = \zeta -\int_0^tF(\bar{s},V^{\zeta,Z}_{\bar{s}},Z_{\bar{s}}) \ud s + \int_0^t Z_{\bar{s}} \ud W_s,
\end{align}
where we introduced the notation $\bar{s} := t_n$ for $t_n \le s <t_{n+1}$.
The second one, denoted $\tilde{V}$, corresponds to the the true solution to the BSDE
\begin{align}\label{eq de perturbed scheme}
\tilde{V}^{\zeta,Z}_t = \zeta - \int_0^t  \tilde{F}(s,\tilde{V}^{\zeta,Z}_{s}, Z_{s}) \ud s + 
\int_0^t Z_{s}  \ud W_s \;,
\end{align}
\noindent where $\tilde{F}$ satisfies the same assumptions as $F$ above.

\begin{Proposition} 
\label{pr key basic estimate} 
Let Assumption \ref{ass random driver} hold for $F$ and $\tilde{F}$.
For $(\zeta,\zeta') \in \cL^2(\cF_0)\times \cL^2(\cF_0)$ and $(Z,Z')\in \cS^2_d\times\cS^2_d$, we consider $V^{\zeta,Z}$ and $\tilde{V}^{\zeta',Z'}$ as defined in \eqref{eq de random driver scheme}-\eqref{eq de perturbed scheme} and we set 
$\delta F :=  \tilde{F}(\cdot,{V}^{\zeta,Z} ,Z) - F(\cdot,{V}^{\zeta,Z} ,Z)$, $  \eta^f_s  = \tilde{F}(s, \tilde{V}^{\zeta' , Z'}_{s}, Z'_{s} ) -  \tilde{F}(s, \tilde{V}^{\zeta' ,Z' }_{\bar{s}}, Z'_{\bar{s}} )$ and $\eta_s^z = Z'_s - Z'_{\bar{s}} $. Then, under the above assumptions on $F$ and $\tilde{F}$, it holds
\begin{enumerate}
\item Forward estimate:
\begin{align*}
\esp{\sup_{t \in [0,T]} | \tilde{V}^{\zeta',Z'}_t -V^{\zeta,Z}_t |^2} &\le C\left(\esp{|\zeta-\zeta' |^2 
+h \sum_{i=0}^{N-1}  |Z_{t_{n}}-Z'_{t_{n}}|^2}
+ \esp{\int_0^T(|\delta F_{\bar{s}}|^2+ |\eta^f_s|^2 + |\eta^z_s|^2 )\ud s }
\right). 
\end{align*}

\item Backward estimate:
\begin{align*}
\esp{\sup_{t\in [0,T]} | \tilde{V}^{\zeta',Z'}_t -V^{\zeta,Z}_t|^2
+ h \sum_{i=0}^{N-1}  |Z_{t_{n}}-Z'_{t_{n}}|^2
} \le 
C\esp{|\tilde{V}^{\zeta',Z'}_T -V^{\zeta,Z}_T|^2 + \int_0^T(|\delta F_{\bar{s}}|^2+|\eta^f_s|^2 + |\eta^z_s|^2)\ud s}.
\end{align*}
\end{enumerate}
\end{Proposition}
\proof  
\\
1. Denote $\Delta V:= \tilde{V}^{\zeta',Z'}-V^{\zeta,Z}  $, $\Delta Z = Z' - Z$, $\Delta F =  \tilde{F}(\cdot,\tilde{V}^{\zeta',Z'} ,Z') - \tilde{F}(\cdot,{V}^{\zeta, Z} ,Z)$ and $\Delta \Gamma_s = \Delta Z_{\bar{s}} + \eta^z_s$. Applying It\^o's formula, we compute
\begin{align}\label{eq basic ito step}
|\Delta V_t |^2 = |\Delta V_0 |^2 + \int_0^t\set{-2\Delta V_s(\Delta F_{\bar{s}} + \delta F_{\bar{s}} +\eta^f_s)+ |\Delta \Gamma_s|^2}\ud s + 2\int_0^t\Delta V_s \Delta \Gamma_s\ud W_s\,.
\end{align}
Since $\tilde{F}$ is Lipschitz-continuous, we have
\begin{align*}
2|\Delta V_s(\Delta F_{\bar{s}} + \delta F_{\bar{s}} + \eta^f_s)|\le (4+L) \sup_{0\le r \le s}|\Delta V_r|^2 + L |\Delta Z_{\bar{s}}|^2 + |\delta F_{\bar{s}}|^2+ |\eta_s^f|^2
\end{align*}
which combined with  \eqref{eq basic ito step} leads to
\begin{align}\label{eq basic step one}
\esp{\sup_{0\leq s\le t}|\Delta V_s |^2} &\le
\esp{ |\Delta V_0 |^2 + C \int_0^t \set{  \sup_{0\le r \le s}|\Delta V_r|^2 +  |\Delta Z_{\bar{s}}|^2 +|\delta F_{\bar{s}}|^2+ |\eta_s^f|^2+  |\eta_s^z|^2}\ud s 
} \nonumber
\\
&\;\;+ 2\, \esp{
 \sup_{0\leq r \le t} |\int_0^r \Delta V_s (\Delta Z_{\bar{s}} +\eta^z_s)\ud W_s|}.
\end{align}
Applying the Burkholder-Davis-Gundy inequality, we obtain
\begin{align*}
\esp{\sup_{r \in[0,t]} |\int_0^r\Delta V_s (\Delta Z_{\bar{s}} +\eta^z_s)\ud W_s|}
&\le C \esp{|\int_0^t|\Delta V_s (\Delta Z_{\bar{s}} +\eta^z_s)|^2 \ud s|^\frac12}
\\
&\le C\left( \esp{\sup_{0\leq s\le t} |\Delta V_s|^2 + \int_0^t( |\Delta Z_{\bar{s}} |^2 + |\eta^z_s|^2 )\ud s } \right)
\end{align*}
where we used Young's inequality for the last inequality.
Inserting the previous inequality into \eqref{eq basic step one}, we get
\begin{align*}
\esp{\sup_{0\leq r\le t} |\Delta V_t|^2} \le  |\Delta V_0 |^2 + C\int_0^t \esp{\sup_{0\leq r\le s} |\Delta V_r|^2 + |\Delta Z_{\bar{s}}|^2 + |\eta^f_s|^2 +|\delta F_{\bar{s}}|^2+|\eta^z_s|^2 } \ud s\,.
\end{align*}
The proof for this step is concluded by applying Gr\"onwall's Lemma.
\vspace{2mm}

\noindent 2.  From \eqref{eq basic ito step}, we compute
\begin{align*}
\esp{|\Delta V_t|^2 + \int_t^T |\Delta \Gamma_s|^2 \ud s}
&\le
\esp{|\Delta V_T|^2 + 2 \int_t^T \Delta V_s(\Delta F_{\bar{s}}  + \delta F_{\bar{s}} +\eta^f_s) \ud s}.
\end{align*}
We observe that, since $\tilde{F}$ is Lipschitz continuous,
\begin{align*}
\Delta V_s \Delta F_{\bar{s}} \le C(|\Delta V_s|^2 + |\Delta V_{\bar{s}}|^2 + |\Delta V_s \Delta Z_{\bar{s}}| ).
\end{align*}
For $\alpha > 0$, to be fixed later on, we  get using Young's inequality,
\begin{align*}
\Delta V_s \Delta F_{\bar{s}} &\le C\left((1+\frac1{\alpha})|\Delta V_s|^2 + |\Delta V_{\bar{s}}|^2 + \alpha|\Delta Z_{\bar{s}}|^2 \right),
\\
&\le  C\left((1+\frac1{\alpha})|\Delta V_s|^2 + |\Delta V_{\bar{s}}|^2 + \alpha|\Delta \Gamma_{s}|^2 + |\eta^z_s|^2 \right).
\end{align*}
For $\alpha$ small enough, we thus obtain
\begin{align}\label{eq step 1 backward}
\esp{|\Delta V_t|^2 + \frac12\int_t^T |\Delta \Gamma_s|^2 \ud s}
&\le
\esp{|\Delta V_T|^2 + C \int_t^T\left( |\Delta V_s|^2 +|\Delta V_{\bar{s}}|^2 
+ |\delta F_{\bar{s}}|^2
+ |\eta^f_s|^2 + |\eta^z_s|^2 \right) \ud s}.
\end{align}
Applying Gr\"onwall's Lemma leads to, for all $t \le T$,
\begin{align}\label{eq step 2 backward}
\esp{|\Delta V_t|^2 } \le C\esp{ \mathscr{B}_T + \int_t^T |\Delta V_{\bar{s}}|^2  \ud s }
\text{ with } \mathscr{B}_T :=  |\Delta V_T|^2 + \int_0^T (|\eta^f_s|^2 + |\eta^z_s|^2+ |\delta F_{\bar{s}}|^2) \ud s\;.
\end{align}
In particular, for $n \le N$ and $t_n \in \pi$, we have
\begin{align*}
\esp{|\Delta V_{t_n}|^2 } \le C\esp{ \mathscr{B}_T + h\sum_{j=n}^{N-1} |\Delta V_{t_j}|^2}
\end{align*}
which in turn, using the discrete-time Gr\"onwall Lemma, leads to $\max_{n\le N}\esp{|\Delta V_{t_n}|^2 }  \le C \mathbb{E}[\mathscr{B}_T]$.
Combining this inequality with \eqref{eq step 1 backward} and \eqref{eq step 2 backward}, we obtain
\begin{align*}
\esp{|\Delta V_t|^2 + \frac12\int_t^T |\Delta \Gamma_s|^2 \ud s} \le C \mathbb{E}[\mathscr{B}_T]\,,\, t \le T.
\end{align*}
To conclude the proof one applies the Burkholder-Davis-Gundy inequality as in step 1.
\eproof

\subsection{Application to the \emph{direct algorithm}} \label{subse application direct algo}
We here prove the results announced in Section \ref{subse direct algo}. We start by proving the analytic expression of the main quantities appearing in the \emph{direct algorithm}, recall Definition \ref{de implemented direct algo}. 

\vspace{4mm}
\noindent \textbf{Proof of Lemma \ref{pr expression deriv direct algo}}
By standard computations, recall \eqref{eq de Z}, for any $0\leq n\leq N-1$,
\begin{align*}
 \nabla_{\mathfrak{y}^k} Z^{\mathfrak{u}}_{t_n} = 0,  \quad  \text{ and } \quad  \nabla_{\mathfrak{z}^{n,k}} (Z^\mathfrak{u}_{t})^l =  \psi_n^{k}(X_{t_n}) \1_{\set{t=t_n}} \mathbf{e}^l \,.
\end{align*}
This leads to, for $0 \le q \le N-1$,
\begin{align} \label{eq deriv term in Z:1}
\nabla_{\mathfrak{z}^{n,k}} (Z^\mathfrak{u}_{t_q}\cdot\Delta W_q) =  \psi_n^{k}(X_{t_n})    \Delta W_n\1_{\set{q=n}}
\end{align}
\noindent and 
\begin{align}
 \nabla_{\mathfrak{z}^{n,k}} f(Y^{\mathfrak{u}}_{t_q},Z^\mathfrak{z}_{t_q}) = \nabla_y f(Y^{\mathfrak{u}}_{t_q}, Z^{\mathfrak{u}}_{t_q}) \nabla_{\mathfrak{z}^{n,k}}Y^{\mathfrak{u}}_{t_q} + \psi_n^{k}(X_{t_n})   \nabla_z f(Y^{\mathfrak{u}}_{t_q},Z^\mathfrak{u}_{t_q})\1_{\set{q=n}}. \label{eq deriv term in Z:2}
 \end{align}
Now, differentiating  both side of \eqref{eq de bar Y} with respect to the variable $\mathfrak{y}^k, 1\le k\le K^y $ yields
\begin{align*}
\nabla_{\mathfrak{y}^k} Y^{\mathfrak{u}}_{t_n} =  \psi_y^k(X_0) \prod_{j = 0}^{n-1} \left(1- h \nabla_y f(Y^{\mathfrak{u}}_{t_j}, Z^\mathfrak{u}_{t_j}) \right) \text{ for } n \ge 0.
\end{align*}
From \eqref{eq de bar Y}, by differentiation, we obtain $\nabla_{\mathfrak{z}^{n, k}} {Y}^{\mathfrak{u}}_{0} = 0$ and using \eqref{eq deriv term in Z:1} and \eqref{eq deriv term in Z:2}, for $q \ge 1$,
\begin{align*}
\nabla_{\mathfrak{z}^{n,k}} {Y}^{\mathfrak{u}}_{t_{q}} = \nabla_{\mathfrak{z}^{n, k}} {Y}^{\mathfrak{u}}_{t_{q-1}}\left(1- h \nabla_y f({Y}^{\mathfrak{u}}_{t_{q-1}},Z^{\mathfrak{z}}_{t_{q-1}} ) \right) 
+ \psi_n^{k}(X_{t_n})\1_{\set{n=q-1}}   \left( \Delta W_{q-1}^\top  - h \nabla_z f( Y^{\mathfrak{u}}_{t_{q-1}}, Z^\mathfrak{u}_{t_{q-1}})
\right ) \,, 
\end{align*}
which in turn yields
\begin{align*}
&\nabla_{\mathfrak{z}^{n,k}} Y^{\mathfrak{u}}_{t_q} = 0 \text{ for } q \le n\;,\quad 
\nabla_{\mathfrak{z}^{n,k}} Y^{\mathfrak{u}}_{t_{n+1}} =  \psi^{k}_n(X_{t_n})  \left( \Delta W_{n}^{\top}  - h\nabla_z f( Y^{\mathfrak{u}}_{t_{n}}, Z^\mathfrak{u}_{t_{n}})
\right )
\end{align*}
\noindent and for $q \ge n+2$,
\[ \nabla_{\mathfrak{z}^{n,k}} Y^{\mathfrak{u}}_{t_{q}} =  \psi_{n}^{k}(X_{t_n})   \left( \Delta W_{n}  - h\nabla_z f( Y^{\mathfrak{u}}_{t_{n}}, Z^\mathfrak{u}_{t_{n}}) \right ) \prod_{j = n+1}^{q-1} \left(1- h \nabla_y f( Y^{\mathfrak{u}}_{t_j}, Z^\mathfrak{u}_{t_j}) \right).\]
 This concludes the proof.
\eproof


\vspace{4mm}
\noindent Recall that, the time discretization error that will appear in our estimates is classically given by
\begin{align}
\cE_{\pi} &=  \esp{ \sum_{n=0}^{N-1}  \int_{t_{n}}^{t_{n+1}}\left( |\cY_s-\cY_{t_{n}}|^2 + |\cZ_s-\cZ_{t_{n}}|^2 + |\cX_{t_{n}} - \widehatbis{X}_{t_n}|^2 \right)ds } \,. \label{eq error disc time bis}
\end{align}

The approximation error due to the restriction to the functional space, expressed in  \eqref{eq error approx space brut}, is also given by
\begin{align}
\cE_{\psi} &= \esp{|u(0,\cX_0) - \bar{u}_0(\cX_0)|^2 + \sum_{n=1}^{N-1} h | (\sigma^\top\nabla_x u)(t_n,\widehatbis{X}_{t_n})- \bar{v}_n(\widehatbis{X}_{t_n})|^2} \,,
\label{eq error approx space}
\end{align}
where, $\bar{v}_n$ is the $L^{2}(\mathbb{R}^d, \mathbb{P}_{\widehatbis{X}_{t_n}})$-projection of the map $(\sigma^\top\nabla_x u)(t_n,\cdot)$ onto $\mathscr{V}^z_n$, $0 \le n \le N-1$ and $\bar{u}_0$  is the $L^{2}(\mathbb{R}^d, \mathbb{P}_{\cX_0})$-projection of the map $u(0,\cdot)$ onto $\mathscr{V}^y$.
We denote $\bar{\mathfrak{y}}$ the coefficient associated to the decomposition of $\bar{u}$, namely
\begin{align}\label{y:reference:solution}
\R^d \ni x \mapsto \bar{u}_0(x) = \sum_{k = 1}^{K^y} \bar{\mathfrak{y}}^k \psi_0^{k}(x) \in \R, 
\end{align}
 and also $\bar{\mathfrak{z}}^n$ the coefficient associated to the decomposition of $\bar{v}_n$, namely
\begin{align}\label{theta:reference:solution}
\R^d \ni x \mapsto \bar{v}_n(x) = \sum_{k = 1}^{K^z_n} \bar{\mathfrak{z}}^{n,k} \psi_n^{k}(x) \in \R^d, \, 0\leq n \leq N-1.
\end{align}
For later use, we introduce a \emph{reference solution} 
\begin{equation}\label{the:reference:solution}
\bar{\mathfrak{u}}=(\bar{\mathfrak{y}},\bar{\mathfrak{z}}) \in \R^{K^y}\times \R^{d \bar{K}^z}  \,.
\end{equation}

We first discuss the well-posedness of the optimization problem \eqref{eq de theoretical direct algo}.
\begin{Lemma} \label{le well-posedness}
Under Assumption \ref{ass X0}, Assumption \ref{assumption:coefficient} (i) and Assumption \ref{as psi function}, it holds
\begin{align*}
\argmin_{\mathfrak{u} \in \R^{K^y}\times \R^{d \bar{K}^z}} \mathfrak{g}(\mathfrak{u}) \neq \emptyset \;.
\end{align*}
\end{Lemma}

\proof
Let $\mathfrak{u}=(\mathfrak{y},\mathfrak{z}) \in  \R^{K^y}\times \R^{d \bar{K}^z} $. Using the backward estimate of Proposition \ref{pr key basic estimate} with
$V^{\zeta,Z} := Y^{\mathfrak{u}}$ and $\tilde{V}^{\zeta', Z'} := Y^{0}$ ($\zeta'=0$, $Z'\equiv0$) yields
\begin{align*}
\vvvert \mathfrak{u} \vvvert^2= \esp{|Y^{\mathfrak{u}}_0|^2 +\sum_{n=0}^{N-1} h |Z^{\mathfrak{u}}_{t_n}|^2} &\le C \esp{|Y^{\mathfrak{u}}_T - Y^{0}_T|^2} \le C(1 + \esp{|g(\widehatbis{X}_T)-Y^{\mathfrak{u}}_T|^2}).
\end{align*}

Under Assumption \ref{as psi function}, we thus deduce that the continuous function $\R^{K^y}\times \R^{d \bar{K}^z}  \ni \mathfrak{u} \mapsto \mathfrak{g}(\mathfrak{u})$ is coercive. As a consequence, it admits a global minimizer so that the optimization problem \eqref{eq de theoretical direct algo} is well-posed.
\eproof


The following proposition can be seen as a version of the results in  \cite{han2020convergence} (see Theorem 1 \& 2) adapted to our context. Let us note that our setting is simpler as we do not deal with fully coupled Forward Backward SDEs.

\begin{Proposition}\label{pr conv direct algo}
Under Assumption \ref{ass X0}, Assumption \ref{assumption:coefficient} (i), (ii) and Assumption \ref{as psi function}, there exists a positive constant $C$ such that for any
 \[\mathfrak{u}^\star:=(\mathfrak{y}^\star,\mathfrak{z}^\star) \in \argmin_{\mathfrak{u} \in \R^{K^y}\times \R^{d \bar{K}^z} } \mathfrak{g}(\mathfrak{u}),\] 
 
 \noindent it holds
\begin{align}\label{eq main control direct algo}
\quad \esp{|u(0,\cX_0)-Y_0^{\mathfrak{u}^\star}|^2 +  h \sum_{n=0}^{N-1}  |\cZ_{t_{n}}-Z^{\mathfrak{u}^\star}_{t_{n}}|^2}\le C\left(  \cE_{\pi} + \cE_{\psi}
\right) \,
\end{align}
and
\begin{align} \label{eq co bound error ref solution}
\mathfrak{g}(\mathfrak{u}^\star) \le \mathfrak{g}(\bar{\mathfrak{u}}) \le C\left(  \cE_{\pi} + \cE_{\psi} \right)\;.
\end{align}

\end{Proposition}

\proof We use the backward estimate of Proposition \ref{pr key basic estimate} with
$V^{\zeta,Z} := Y^{\mathfrak{u}^\star}$ and $\tilde{V}^{\zeta',Z} := Y^{\bar{\mathfrak{u}}}$ where $\bar{\mathfrak{u}}$ is given in
\eqref{the:reference:solution}, to obtain
\begin{align*}
\esp{\sup_{t \in [0,T]}|Y^{\mathfrak{u}^\star}_t - Y^{\bar{\mathfrak{u}}}_t|^2 +  h \sum_{n=0}^{N-1}  \esp{|Z^{\bar{\mathfrak{u}}}_{t_{n}}-Z^{\mathfrak{u}^\star}_{t_{n}}|^2}}
&\le
C\esp{|Y^{\mathfrak{u}^\star}_T - Y^{\bar{\mathfrak{u}}}_T|^2}
\\
&\le C\esp{|Y^{\mathfrak{u}^\star}_T - g(X^x_T)|^2 + |g(X^x_T) - Y^{\bar{\mathfrak{u}}}_T|^2}.
\end{align*}
By optimality of $\mathfrak{u}^\star$, we get
\begin{align}\label{eq first conv direct algo theo}
\esp{\sup_{t \in [0,T]}|Y^{\mathfrak{u}^\star}_t - Y^{\bar{\mathfrak{u}}}_t|^2 +  h \sum_{n=0}^{N-1}  \esp{|Z^{\bar{\mathfrak{u}}}_{t_{n}}-Z^{\mathfrak{u}^\star}_{t_{n}}|^2}}
&\le C\esp{ |g(X^x_T) - Y^{\bar{\mathfrak{u}}}_T|^2}\;.
\end{align}
We now use the forward estimate of  Proposition \ref{pr key basic estimate} with
$V^{\zeta,Z} := Y^{\bar{\mathfrak{u}}}$ and $\tilde{V}^{\zeta',Z} := \cY$. We obtain
\begin{align}\label{eq second conv dir algo theo}
\mathfrak{g}(\bar{\mathfrak{u}}) = \esp{ |g(X^x_T) - Y^{\bar{\mathfrak{u}}}_T|^2}
\le 
C\esp{  |\cY_0 - Y^{\bar{\mathfrak{u}}}_0|^2+  \int_{0}^{T}\left(|\cY_s-\cY_{\bar{s}}|^2+ |\cZ_s-\cZ_{\bar{s}}|^2\right) \ud s+\sum_{n= 0}^{N-1}  h |\cZ_{t_{n}}-Z^{\bar{\mathfrak{u}}}_{t_{n}}|^2}\,.
\end{align}

\noindent Observe now that in our current smooth coefficients framework $\cZ_{t_n} = (\sigma^\top \nabla_x u)(t_n,\cX_{t_n})$ so that one has
\begin{align}
\esp{|\cZ_{t_n}-Z^{\bar{\mathfrak{u}}}_{t_n}|^2} & =  \mathbb{E}\Big[|(\sigma^\top\nabla_x u)(t_n,\cX_{t_n})-Z^{\bar{\mathfrak{u}}}_{t_n} |^2 \Big] \nonumber \\
& \leq 2 \Big( \mathbb{E}\Big[|(\sigma^\top \nabla_x u)(t_n,\cX_{t_n})- (\sigma^\top \nabla_x u)(t_n,\widehatbis{X}_{t_n}) |^2 \Big] + \mathbb{E}\Big[|(\sigma^\top\nabla_x u)(t_n,\widehatbis{X}_{t_n})-Z^{\bar{\mathfrak{u}}}_{t_n} |^2 \Big]\Big) \nonumber\\
& \leq C \Big( \mathbb{E}\Big[|\cX_{t_n} - \widehatbis{X}_{t_n}|^2 \Big] + \mathbb{E}\Big[|(\sigma^\top\nabla_x u)(t_n,X^{x_0}_{t_n})-Z^{\bar{\mathfrak{u}}}_{t_n} |^2 \Big]  \Big) \label{ineq:control:error:Z}
\end{align}
\noindent where we used the Lipschitz regularity of $x\mapsto (\sigma^\top \nabla_x u)(t_n, x)$ uniformly with respect to the variable $t_n$. 
Combining \eqref{eq second conv dir algo theo} and \eqref{ineq:control:error:Z} leads to
\begin{align*}
\mathfrak{g}(\bar{\mathfrak{u}})
\le
C\left(  \cE_{\pi} + \cE_{\psi}  \right)\;
\end{align*}
which, since $\mathfrak{g}(\mathfrak{u}^\star) \le \mathfrak{g}(\bar{\mathfrak{u}})$, proves \eqref{eq co bound error ref solution}.

The above estimate together with \eqref{eq first conv direct algo theo}  leads to
\begin{align*}
\esp{\sup_{t \in [0,T]}|Y^{\mathfrak{u}^\star}_t - Y^{\bar{\mathfrak{u}}}_t|^2 +  h \sum_{n=0}^{N-1}  \esp{|Z^{\bar{\mathfrak{u}}}_{t_{n}}-Z^{\mathfrak{u}^\star}_{t_{n}}|^2}}
&\le  C\left(  \cE_{\pi} + \cE_{\psi}  \right)\,.
\end{align*}

\noindent Then, using the inequalities $|\cZ_{t_n}-Z^{\mathfrak{u}^\star}_{t_n}|^2 \le 2|Z^{\mathfrak{u}^\star}_{t_n} - Z^{\bar{\mathfrak{u}}}_{t_n}|^2+2|\cZ_{t_n}-Z^{\bar{\mathfrak{u}}}_{t_n}|^2$, $0\leq n \leq N-1$,  
and $|Y^{\mathfrak{u}^\star}_t- \cY_t|^2 \le 2  |Y^{\mathfrak{u}^\star}_t - Y^{\bar{\mathfrak{u}}}_t|^2+ 2|\cY_t- Y^{\bar{\mathfrak{u}}}_t|^2$
yields \eqref{eq main control direct algo} and concludes the proof.
\eproof


\subsection{Study of the \emph{Picard algorithm}} 
 \label{Study of the Picard algorithm}

We introduce the following mean squared error:
\begin{align}\label{eq de error picard}
\cE_p := \esp{\vvvert \mathfrak{u}^p_M - \bar{\mathfrak{u}} \vvvert^2}\,, 0 \le p \le P\,,
\end{align}
where the sequence $(\mathfrak{u}^{p}_{M})_{0\leq p \leq P}$ is given by Definition \ref{de picard scheme full}, $\vvvert\cdot\vvvert$ is given by Definition \ref{de double - triple norm} and $\bar{\mathfrak{u}}$ is the reference solution introduced in \eqref{the:reference:solution}. In this subsection, our aim is to establish an upper bound for the quantity $\cE_P$ that will allow us to prove Theorem \ref{th main conv result}.

\subsubsection{Preliminary estimates}


\begin{Proposition} \label{error of Picard}
Suppose that Assumption \ref{ass X0}, Assumption \ref{assumption:coefficient} (i), (ii) and Assumption \ref{as psi function} hold. If $T(1+ 2L^2(1+ h)) < 1$ and
$\delta_h:=\frac{ 8 L^2 T}{1-T(1+ 2L^2(1+ h))} < 1$, then for any $\varepsilon>0$ such that $\delta_{h,\varepsilon}:= \delta_h (1+\varepsilon)<1$ there exists a positive constant $C_\varepsilon$ such that for any positive integer $P$
\begin{align}
\cE_P \le \delta_{h,\varepsilon}^P \cE_0 +  C_\varepsilon \big( \cE_{\mathrm{RM}}  + \cE_{\psi} +  \cE_{\pi} \big)
\end{align}
\noindent with the notation
\begin{align}
\cE_{\mathrm{RM}} := \max_{1 \le p \le P}\esp{ \vvvert \mathfrak{u}^{p}_{M} - \Phi(\mathfrak{u}^{p-1}_{M}) \vvvert^2}.
\end{align}
\end{Proposition}

\proof 
From the decomposition, 
\begin{align*}
\mathfrak{u}^p_{M}- \bar{\mathfrak{u}} = &\Phi_{M}(\mathfrak{u}^{p-1}_{M}) - \Phi(\mathfrak{u}^{p-1}_{M}) + \Phi(\mathfrak{u}^{p-1}_{M})  -  \bar{\mathfrak{u}} 
\end{align*}
\noindent we obtain, for any $\varepsilon > 0$,
\begin{align*}
\cE_p \le (1+ \frac1\varepsilon)\cE_{\mathrm{RM}} + (1+\varepsilon)\esp{ \vvvert \Phi(\mathfrak{u}^{p-1}_{M})  -  \bar{\mathfrak{u}}  \vvvert^2}.
\end{align*}
Then, using Lemma \ref{le control f(y,z)} below, we get
\begin{align*}
\cE_p \le \delta_{h,\varepsilon} \cE_{p-1} +  C_\varepsilon \left(  \cE_{\psi} + \cE_{\pi}  + \cE_{\mathrm{RM}}  \right)
\end{align*}
\noindent up to a modification of $\varepsilon$. By an induction argument, we  derive
\begin{align*}
\cE_P \le \delta_{h,\varepsilon}^P \cE_0 +  C_\varepsilon \big( \cE_{\mathrm{RM}} + \cE_{\psi} +  \cE_{\pi} \big)
\end{align*}
which concludes the proof.
\eproof

\begin{Lemma}  \label{le control f(y,z)}
Suppose that Assumption \ref{ass X0}, Assumption \ref{assumption:coefficient} (i), (ii) and Assumption \ref{as psi function} hold. If $T(1+ 2L^2(1+ h)) < 1$ and
$\delta_h:=\frac{ 8 L^2 T}{1-T(1+ 2L^2(1+ h))} < 1$, then, for any $\varepsilon>0$ there exists a positive constant $C_\varepsilon$ ($\varepsilon \mapsto C_\varepsilon$ being non-increasing) such that for any $\tilde{\mathfrak{u}} \in  \R^{K^y}\times \R^{d\bar{K}^z}$  it holds 
\begin{align*}
\vvvert \Phi(\tilde{\mathfrak{u}})  -  \bar{\mathfrak{u}} \vvvert^2
\le \delta_{h}(1+\varepsilon) \vvvert \tilde{\mathfrak{u}}  -  \bar{\mathfrak{u}} \vvvert^2
+  C_\varepsilon \left(  \cE_{\psi} + \cE_{\pi}  \right).
\end{align*}
\end{Lemma}

\proof 

\emph{Step 1:} We  denote $\check{\mathfrak{u}} = \Phi(\tilde{\mathfrak{u}})  $ where
$\check{\mathfrak{u}} = (\check{\mathfrak{y}},\check{\mathfrak{z}})$ and $\tilde{\mathfrak{u}} = (\tilde{\mathfrak{y}},\tilde{\mathfrak{z}})$ belongs to $\mathbb{R}^{K^y} \times \mathbb{R}^{d \bar{K}^z}$. We first observe that, recalling \eqref{eq de Phi}, \eqref{eq rewritte norms} and \eqref{eq de U starting point},
\begin{align*}
\vvvert \Phi(\tilde{\mathfrak{u}})  -  \bar{\mathfrak{u}} \vvvert^2 &=
\esp{ |Y^{\bar{\mathfrak{u}}}_0-Y^{\check{\mathfrak{u}}}_0|^2 + \int_0^T |Z^{\check{\mathfrak{u}}}_{t}-Z^{\bar{\mathfrak{u}}}_{t}|^2 \ud t}
\\
&=  \esp{|U^{\tilde{\mathfrak{u}}, \check{\mathfrak{u}}}_T
- U^{\tilde{\mathfrak{u}},\bar{\mathfrak{u}}}_T|^2}.
\end{align*}
Moreover, by optimality of $\check{\mathfrak{u}}$
\begin{align*}
\esp{|U^{\tilde{\mathfrak{u}},\check{\mathfrak{u}} }_T
- U^{\tilde{\mathfrak{u}},\bar{\mathfrak{u}}}_T|^2}
&\le 
2 \left(
\esp{|U^{\tilde{\mathfrak{u}},\check{\mathfrak{u}} }_T
- g({X}_T)|^2}
+
\esp{|g({X}_T)
- U^{\tilde{\mathfrak{u}},\bar{\mathfrak{u}}}_T|^2}
\right)
\\
&\le 4 \esp{|g({X}_T)
- U^{\tilde{\mathfrak{u}},\bar{\mathfrak{u}}}_T|^2}.
\end{align*}

We now compute, for any $\varepsilon > 0$,
\begin{align*}
\esp{|g({X}_T)
- U^{\tilde{\mathfrak{u}},\bar{\mathfrak{u}}}_T|^2}
\le
(1+\frac1\varepsilon)
\esp{|g({X}_T)
- U^{\bar{\mathfrak{u}},\bar{\mathfrak{u}}}_T|^2}
+
(1+\varepsilon)\esp{|U^{\bar{\mathfrak{u}},\bar{\mathfrak{u}}}_T
- U^{\tilde{\mathfrak{u}},\bar{\mathfrak{u}}}_T|^2}\,,
\end{align*}
which, combined with the previous inequality, yields
\begin{align} \label{eq end step 1}
 \vvvert \Phi(\tilde{\mathfrak{u}})  -  \bar{\mathfrak{u}} \vvvert^2  \le 4(1+\frac1\varepsilon)
\esp{|g({X}_T)
- U^{\bar{\mathfrak{u}},\bar{\mathfrak{u}}}_T|^2}
+
4(1+\varepsilon)\esp{|U^{\bar{\mathfrak{u}},\bar{\mathfrak{u}}}_T
- U^{\tilde{\mathfrak{u}},\bar{\mathfrak{u}}}_T|^2}\,.
\end{align}
Since $U^{\bar{\mathfrak{u}},\bar{\mathfrak{u}}} = Y^{\bar{\mathfrak{u}}}$, we can give an upper bound for the first term appearing on the right-hand side of the above inequality by using \eqref{eq co bound error ref solution},
\begin{align}\label{eq bias part}
\esp{|g({X}_T)
- U^{\bar{\mathfrak{u}},\bar{\mathfrak{u}}}_T|^2} = \mathfrak{g}(\bar{\mathfrak{u}})
\le C \left(  \cE_{\psi} + \cE_{\mathrm{\pi}}  \right).
\end{align}
\emph{Step 2: } 
We now turn to the study of the second term appearing on the right hand side of \eqref{eq end step 1}.
Recalling the dynamics \eqref{eq scheme picard}, denoting $\delta U := U^{\tilde{\mathfrak{u}},\bar{\mathfrak{u}}} - Y^{{\bar{\mathfrak{u}}}} $, $\delta Z = Z^{\tilde{\mathfrak{u}}}-Z^{\bar{\mathfrak{u}}}$, $\delta Y  = Y^{\tilde{\mathfrak{u}}} - Y^{\bar{\mathfrak{u}}}$ and $\delta f_{t_n} = f(Y^{\tilde{\mathfrak{u}}}_{t_n}, Z^{\tilde{\mathfrak{u}}}_{t_n}) - f(Y^{\bar{\mathfrak{u}}}_{t_n}, Z^{\bar{\mathfrak{u}}}_{t_n})$, $0\leq n \leq N-1$, using the Cauchy-Schwarz inequality and the Lipschitz-regularity of the map $f$, we get
\begin{align}
\esp{|\delta U_T|^2} &\le T\int_0^T\esp{|\delta f_{\bar{s}}|^2} \ud s  \le 2L^2 T \sum_{n=0}^{N-1}h \esp{|\delta Y_{t_n}|^2 + |\delta Z_{t_n}|^2} \ud s. 
\label{eq first upper bound}
\end{align}
For all $0 \leq n \le N-1$, one has
\begin{align*}
\delta Y_{t_{n+1}}  =& \delta Y_{t_n} - h \delta f_{t_n} + \delta Z_{t_n} \Delta W_n
\end{align*}

\noindent which in turn, setting $\Delta M_n := 2(\delta Y_{t_n} - h \delta f_{t_n}) \delta Z_{t_n} \Delta W_n$, yields
\begin{align*}
|\delta Y_{t_{n+1}}|^2 = |\delta Y_{t_n}|^2 - 2 h \delta Y_{t_n} \delta f_{t_n} + h^2 |\delta f_{t_n} |^2 + | \delta Z_{t_n} \Delta W_n|^2 + \Delta M_n.
\end{align*}
Using the fact that $\mathbb{E}[\Delta M_n]=0$ and the Lipschitz regularity of the map $f$, we deduce
\begin{align*}
\esp{|\delta Y_{t_{n+1}}|^2} & \leq \esp{|\delta Y_{t_n}|^2 + h |\delta Y_{t_n}|^2 +  (h+ h^2) |\delta f_{t_n} |^2 + h | \delta Z_{t_n}|^2}\\
& \leq  \esp{|\delta Y_{t_n}|^2 + (h + 2L^2(h+ h^2) )  (|\delta Y_{t_n}|^2 + | \delta Z_{t_n}|^2)}.
\end{align*}
Summing the previous inequality, we obtain
\begin{align*}
\esp{|\delta Y_{t_{n}}|^2}  \le |\delta Y_{0}|^2 + \esp{\sum_{j=0}^{N-1}(h + 2L^2(h+ h^2) ) (|\delta Y_{t_j}|^2 + | \delta Z_{t_j}|^2)}
\end{align*}
so that, multiplying both side of the previous inequality by $h$ and summing again, we get
\begin{align*}
\sum_{n=0}^{N-1}h\esp{|\delta Y_{t_{n}}|^2} & \le T|\delta Y_{0}|^2 + T(1+ 2L^2(1+ h)) \left(\mathbb{E}\Big[\sum_{j=0}^{N-1}h |\delta Y_{t_j}|^2\Big] + \mathbb{E}\Big[\sum_{j=0}^{N-1}h | \delta Z_{t_j}|^2\Big] \right)
\\
& 
\le \frac{T}{1-T(1+ 2L^2(1+ h))} \left( |\delta Y_{0}|^2 + (1+ 2L^2(1+ h))\mathbb{E}\Big[\sum_{j=0}^{N-1}h | \delta Z_{t_j}|^2\Big] \right)
\end{align*}

\noindent where, for the last inequality, we used the fact that $T(1+ 2L^2(1+ h))<1$. Combining the previous inequality with \eqref{eq first upper bound}, we obtain
\begin{align}\label{eq end step 3}
\esp{|\delta U_T|^2}
&\le
  \frac{2L^2 T}{1-T(1+ 2L^2(1+ h))}  \left( T  |\delta Y_{0}|^2
+
 \sum_{n=0}^{N-1} h\esp{|\delta Z_{t_n}|^2}
\right).
\end{align}
We finally complete the proof by combining \eqref{eq end step 1}, \eqref{eq bias part} and \eqref{eq end step 3}.
%
\eproof

\subsubsection{Study of the approximation error of the stochastic gradient descent algorithm} \label{subse error of RM}


In this subsection, our aim is to study the approximation error of $\Phi$ by $\Phi_M$ where $\Phi_M$ has been introduced in Definition \ref{de sgd for phi}.

\begin{Lemma} \label{le tool robbins-monro} Suppose that Assumption \ref{ass X0}, Assumption \ref{assumption:coefficient} (i) and Assumption \ref{as psi function} hold. Let $\tilde{\mathfrak{u}}$ be a fixed $\R^{K^y} \times \R^{d \bar{K}^z}$-valued random vector and set $\check{\mathfrak{u}} =(\check{\mathfrak{y}},\check{\mathfrak{z}}) := \Phi(\tilde{\mathfrak{u}})$ and $\mathfrak{u}_M =({\mathfrak{y}}_M,\mathfrak{z}_M):= \Phi_M(\mathfrak{u}_0,\mathfrak{X}_0,\mathfrak{W},\tilde{\mathfrak{u}})$, $M$ being a positive integer and where $(\mathfrak{u}_0,\mathfrak{X}_0,\mathfrak{W})$ (recall Definition \ref{de sgd for phi}) is independent of $\tilde{\mathfrak{u}}$. We denote by $\mathbb{E}_{\tilde{\mathfrak{u}}}[\cdot]$, the conditional expectation with respect to the sigma-field $\sigma(\tilde{\mathfrak{u}})$ generated by $\tilde{\mathfrak{u}}$. Then, for any positive integer $M$, the random vector $(\mathfrak{y}_M, \mathfrak{z}_M)$ satisfies:
 \begin{align} 
  \mathbb{E}_{\tilde{\mathfrak{u}}} \Big[|\mathfrak{y}_{M} - \check{\mathfrak{y}} |^2 \Big]
 \le 
L_{K,M} (1+|\check{\mathfrak{y}}|^2) 
\quad\textnormal{ and }\quad
 \mathbb{E}_{\tilde{\mathfrak{u}}} \Big[|\mathfrak{z}^{n,\cdot}_{M, l} - \check{\mathfrak{z}}^{n,\cdot}_{l} |^2 \Big]
 \le 
L_{K,M}\left(\frac{1}{h}+|\check{\mathfrak{z}}^{n,\cdot}_{l}|^2\right) \label{eq final step 2 robbins monro}
\end{align}
\noindent for any $l \in \left\{1, \cdots, d\right\}$, with
\begin{align}
 L_{K,M} &:=  \varrho_0 \varrho_1 \Big(1+ \mathbb{E}\Big[  |\mathfrak{u}_{0}|^2\Big]  \Big) \sum_{m=1}^{M} \exp\left(- {4} \frac{{\alpha}_{K}}{\beta_{K}} (\Gamma_{M}-\Gamma_{m}) \right) \gamma_m^2 \label{de LpsiM}
 \\
\Gamma_{m} &:= \sum_{k=1}^m\gamma_k \;,\; m \ge 1\,,
\end{align}
\noindent where the constants $\varrho_0$, $\varrho_1$ are defined respectively in equation \eqref{de varrho 0} and \eqref{de varrho 1} below.
\noindent Moreover, it holds
\begin{align}\label{eq useful robbins-monro}
  \mathbb{E}_{\tilde{\mathfrak{u}}}\Big[\vvvert \Phi_M(\tilde{\mathfrak{u}}) - \Phi(\tilde{\mathfrak{u}})\vvvert^2\Big]
  \le
 \kappa_{K} L_{K,M}\left(1+ d N\right) + \frac{\kappa_{K}}{\alpha_{K}} L_{K,M} \vvvert\Phi(\tilde{\mathfrak{u}})\vvvert^2\,.
\end{align}
\end{Lemma}
\proof

\noindent \emph{Step 1: } We prove the estimate for the difference $\mathfrak{z}^{n,\cdot}_{M, l} - \check{\mathfrak{z}}^{n,\cdot}_{l}$. The proof for $\mathfrak{y}_{M} - \check{\mathfrak{y}}$ follows from similar arguments and we omit some technical details. From \eqref{def:funct:H:2nd:part}, one gets
\begin{align*}
|H^{n,l}(\cX_0,W,\tilde{\mathfrak{u}},\mathfrak{z}^{n,\cdot}_l)|^2 & = \frac{4}{(\beta_K \sqrt{h})^2}|\mathfrak{G}^{\tilde{\mathfrak{u}}} - \omega^{n,\cdot}_l \cdot \mathfrak{z}^{n,\cdot}_l|^2|\omega^{n,\cdot}_l|^2\le 
\frac{8}{\beta^2_K}(|\mathfrak{G}^{\tilde{\mathfrak{u}}}|^2 + h|\tilde{\omega}^{n,\cdot}_l|^2 | \mathfrak{z}^{n,\cdot}_l|^2)|\tilde{\omega}^{n,\cdot}_l|^2\;.
\end{align*}
Under the boundedness Assumption \ref{assumption:coefficient} (i), recalling \eqref{eq de Gamma}, we obtain
\begin{align}
|\mathfrak{G}^{\tilde{\mathfrak{u}}}|^2 \le2\left( |g|^2_\infty+T^2 |f|^2_{\infty} \right)
\end{align}
\noindent so that setting 
\begin{align}
\label{de varrho 0}
\varrho_0 := 8\max \set{ (|g|_\infty+T |f|_{\infty})^2,2}
\end{align}
\noindent and from the lower bound \eqref{eq de beta psi}, we obtain
\begin{align}\label{eq 2nd moment H n l}
\mathbb{E}_{\tilde{\mathfrak{u}}}\Big[|H^{n,l}(\cX_0,W,\tilde{\mathfrak{u}},\mathfrak{z}^{n,\cdot}_l)|^2\Big]
\le
\varrho_0(1 + h|\mathfrak{z}^{n,\cdot}_l - \check{\mathfrak{z}}^{n,\cdot}_l|^2 +h|\check{\mathfrak{z}}^{n,\cdot}_l|^2)\;.
\end{align}

\noindent \emph{Step 2:} We now introduce the natural filtration of the algorithm namely $\cF= (\cF_m)_{0 \le m \le M}$, defined by $\cF_m = \sigma(\mathfrak{u}_0,\cX_0^{(k)}, W^{(k)}, 1 \le k \le m)$, $m\geq1$, and $\cF_0=\sigma(\mathfrak{u}_0)$. From the dynamics \eqref{RM iteration 2}, we directly get  
 \begin{align*}
 |\mathfrak{z}^{n,\cdot}_{l , m+1} - \check{\mathfrak{z}}^{n,\cdot}_{l} |^2  =    |\mathfrak{z}^{n,\cdot}_{l, m} - \check{\mathfrak{z}}^{n,\cdot}_{l} |^2 &- 2 \gamma^z_{m+1} H^{n,l}(\cX_0^{(m+1)},W^{(m+1)}, \tilde{\mathfrak{u}}, \mathfrak{z}^{n,\cdot}_{l, m} )\cdot( \mathfrak{z}^{n,\cdot}_{l, m} - \check{\mathfrak{z}}^{n,\cdot}_{l}) \\
 &+ (\gamma^z_{m+1})^2| H^{n,l}(\cX_0^{(m+1)},W^{(m+1)},\tilde{\mathfrak{u}},\mathfrak{z}^{n,\cdot}_{l , m})|^2 
 \end{align*}
 
 \noindent so that introducing the sequence of $\cF$-martingale increments,
 \begin{align*}
 \Delta M_{m+1} := \Big(\frac1{\beta_K \sqrt{h}}\nabla_{\!\mathfrak{z}^{n,\cdot}_{l}} \mathfrak{H}^{n,l}(\tilde{\mathfrak{u}},\mathfrak{u}_m) - H^{n,l}(\cX_0^{(m+1)},W^{(m+1)},\tilde{\mathfrak{u}},\mathfrak{z}^{n,\cdot}_{l , m})\Big) \cdot(\mathfrak{z}^{n,\cdot}_{l, m} - \check{\mathfrak{z}}^{n,\cdot}_{l}),\,  m\geq 0,
 \end{align*}
 we have
 \begin{align*}
 |\mathfrak{z}^{n,\cdot}_{l, m+1} - \check{\mathfrak{z}}^{n,\cdot}_{l} |^2  
 =   |\mathfrak{z}^{n,\cdot}_{l , m} - \check{\mathfrak{z}}^{n,\cdot}_{l} |^2 
 &- \frac2{\beta_K \sqrt{h}} \gamma^z_{m+1} \nabla_{\!\mathfrak{z}^{n,\cdot}_{l}} \mathfrak{H}^{n,l}(\tilde{\mathfrak{u}},\mathfrak{u}_m)  \cdot( \mathfrak{z}^{n,\cdot}_{l, m} - \check{\mathfrak{z}}^{n,\cdot}_{l} ) 
 \\
 &+ 2 \gamma^z_{m+1} \Delta M_{m+1}  +(\gamma^z_{m+1})^2 | H^{n, l }(\cX_0^{(m+1)},W^{(m+1)},\tilde{\mathfrak{u}},\mathfrak{z}^{n,\cdot}_{l, m})|^2\,.
 \end{align*}
  Now, from the previous equality, \eqref{eq convex prop:2} and the fact that $\mathbb{E}_{\tilde{\mathfrak{u}}}[\Delta M_{m+1} | \cF_m] = 0$, recall \eqref{eq link h and H}, we obtain
\begin{align}
\mathbb{E}_{\tilde{\mathfrak{u}}}\Big[ |  \mathfrak{z}^{n,\cdot}_{l , m+1} - \check{\mathfrak{z}}^{n,\cdot}_{l} |^2 \Big] & \leq  
\mathbb{E}_{\tilde{\mathfrak{u}}}\Big[|   \mathfrak{z}^{n,\cdot}_{l , m} - \check{\mathfrak{z}}^{n,\cdot}_{l}|^2\Big] \Big(1- {4} \frac{h\alpha_{K}}{\beta_K \sqrt{h}}\gamma^z_{m+1}\Big) + (\gamma^z_{m+1})^2 \mathbb{E}_{\tilde{\mathfrak{u}}}\Big[| H^{n,l}(\cX_0^{(m+1)},W^{(m+1)},\tilde{\mathfrak{u}},\mathfrak{z}^{n,\cdot}_{l , m})|^2\Big] \nonumber \\
 & \leq  \mathbb{E}_{\tilde{\mathfrak{u}}}\Big[ |   \mathfrak{z}^{n,\cdot}_{l , m} - \check{\mathfrak{z}}^{n,\cdot}_{l} |^2\Big]
 \Big(1-{4}  \frac{h\alpha_{K}}{\beta_K \sqrt{h}}\gamma^z_{m+1} + \varrho_0 h(\gamma^z_{m+1})^2 \Big) 
 + \varrho_0 (\gamma^z_{m+1})^2(1+h|\check{\mathfrak{z}}^{n,\cdot}_{l}|^2) \label{temp:ineq:sa}
 \end{align}
 
 \noindent where, for the last inequality, we used \eqref{eq 2nd moment H n l} together with the fact that, since $(\cX_0^{(m+1)},W^{(m+1)})$ is independent of $\mathcal{F}_m$ and $\mathfrak{z}_m$ is $\mathcal{F}_m$-measurable, one has
 $$
 \mathbb{E}_{\tilde{\mathfrak{u}}}\Big[ | H^{n, l}(\cX_0^{(m+1)},W^{(m+1)},\tilde{\mathfrak{u}},\mathfrak{z}^{n,\cdot}_{l , m})|^2 | \mathcal{F}_{m}\Big] = \mathbb{E}_{\tilde{\mathfrak{u}}}\Big[| H^{n,l}(\cX_0,W,\tilde{\mathfrak{u}},\mathfrak{z}^{n,\cdot}_{l})|^2 \Big]_{| \mathfrak{z}^{n, \cdot}_{l}= \mathfrak{z}^{n,\cdot}_{l,m}}.
 $$
 Observe now that by the very definition \eqref{eq de gamama} of the sequence $(\gamma^z_m)_{m\geq1}$ and using the fact that $\alpha_K/\beta_K\leq 1$ and $\varrho_0\geq 16$, one gets $1-{4}  \frac{h\alpha_{K}}{\beta_K \sqrt{h}}\gamma^z_{m+1} + \varrho_0 h(\gamma^z_{m+1})^2  = 1-(4 \frac{\alpha_{K}}{\beta_K}\gamma_{m+1} - \varrho_0 \gamma^2_{m+1})\geq 1- (\sqrt{\varrho_0} \gamma_{m+1} - \varrho_0 \gamma^2_{m+1})\geq 1- 1/4 = 3/4$. As a consequence, $\Pi_m := \prod_{k=1}^{m} (1-{4}  \frac{h\alpha_{K}}{\beta_K \sqrt{h}}\gamma^z_{k} + \varrho_0 h(\gamma^z_{k})^2 ) = \prod_{k=1}^{m} (1- {4} \frac{\alpha_{K}}{\beta_{K}}\gamma_{k} + \varrho_0 \gamma_{k}^2)$ is a product of positive terms. From \eqref{temp:ineq:sa}, we thus deduce
 \begin{align}
 	  \mathbb{E}_{\tilde{\mathfrak{u}}}\Big[ |\mathfrak{z}^{n,\cdot}_{l,m+1}&  - \check{\mathfrak{z}}^{n,\cdot}_{l} |^2 \Big]\nonumber
	  \\
	  & \le \Pi_{m+1} \mathbb{E}_{\tilde{\mathfrak{u}}}\Big[ | \mathfrak{z}^{n,\cdot}_{l,0} - \check{\mathfrak{z}}^{n,\cdot}_{l} |^2 \Big]
	  +  \varrho_0 (\frac1{h}+|\check{\mathfrak{z}}^{n,\cdot}_{l}|^2) \sum_{q=1}^{m+1} \dfrac{\Pi_{m+1}}{\Pi_{q}} \gamma_{q}^2 \nonumber \\
	  & \leq  2\varrho_1 \exp\Big(- {4} \frac{\alpha_{K}}{\beta_{K}} \Gamma_{m+1}\Big) \Big(\mathbb{E}_{\tilde{\mathfrak{u}}}\Big[  |\mathfrak{z}^{n,\cdot}_{l,0}|^2\Big] + |\check{\mathfrak{z}}^{n,\cdot}_{l} |^2  \Big) \nonumber \\
	  & \quad + \varrho_1\varrho_0 (\frac1{h}+|\check{\mathfrak{z}}^{n,\cdot}_{l}|^2) \sum_{q=1}^{m+1} \exp\left(- {4} \frac{\alpha_{K}}{\beta_{K}} (\Gamma_{m+1}-\Gamma_{q}) \right) \gamma_q^2\nonumber \\
	  & \leq    \varrho_0 \varrho_1  \Big(1+\mathbb{E}_{\tilde{\mathfrak{u}}}\Big[  |\mathfrak{z}^{n,\cdot}_{l,0}|^2\Big]\Big)  \left(\frac1h+|\check{\mathfrak{z}}^{n,\cdot}_{l}|^2\right)
 \sum_{q=1}^{m+1} \exp\left(- {4} \frac{\alpha_{K}}{\beta_{K}} (\Gamma_{m+1}-\Gamma_{q}) \right) \gamma_q^2 \label{L2:control:distance:to:target}
 \end{align}
 \noindent where we used the standard inequality $1+ x \leq e^x$, the fact that $\rho_0\geq 2$ and introduced the quantity 
 \begin{align}\label{de varrho 1}
\varrho_1:= \exp(\varrho_0 \sum_{m\geq1} \gamma^2_m).
 \end{align}
 This concludes the proof for \eqref{eq final step 2 robbins monro}.
 
\noindent \emph{ Step 3:} Recalling Definition \ref{de double - triple norm} and using \eqref{eq link triple double} as well as Assumption \ref{as psi function}, we directly deduce 
 \begin{align*}
 \mathbb{E}_{\tilde{\mathfrak{u}}}\Big[\vvvert \Phi_M({\tilde{\mathfrak{u}}}) - \Phi({\tilde{\mathfrak{u}}})\vvvert^2 \Big]
 \le \mathbb{E}_{\tilde{\mathfrak{u}}}\Big[\kappa_{K}|\mathfrak{y}_M-\check{\mathfrak{y}}|^2 + h\kappa_{K} \sum_{n=1}^{N-1}\sum_{l=1}^d |\mathfrak{z}^{n,\cdot}_{l,M} -\check{\mathfrak{z}}^{n,\cdot}_{l} |^2 \Big]
 \end{align*}
 
 \noindent so that, using  \eqref{eq final step 2 robbins monro}
 \begin{align*}
 \mathbb{E}_{\tilde{\mathfrak{u}}}\Big[ \vvvert \Phi_M(\tilde{\mathfrak{u}}) - \Phi(\tilde{\mathfrak{u}})\vvvert^2\Big]
 &\le  
 \kappa_{K} L_{K,M}(1+ |\check{\mathfrak{y}} |^2)
+ h\kappa_{K} L_{K,M}\sum_{n=1}^{N-1}\sum_{l=1}^d (\frac1h+ |\check{\mathfrak{z}}^{n,\cdot}_{l} |^2) \,,
\\
& \le 
 \kappa_{K} L_{K,M} +
 \frac{\kappa_{K}}{\alpha_{K}} L_{K,M}\|\check{\mathfrak{y}}\|_y^2
+ h\kappa_{K} L_{K,M} Nd \frac1h+ \frac{\kappa_{K}}{\alpha_{K}} L_{K,M} \|\check{\mathfrak{z}}\|_z^2 \,,
\\
&\le \kappa_{K} L_{K,M}\left(1+ d N\right) + \frac{\kappa_{K}}{\alpha_{K}} L_{K,M} \vvvert\check{\mathfrak{u}}\vvvert^2\,,
 \end{align*}
 which concludes the proof.
\eproof

\vspace{2mm}
The following result provides an upper-bound for the quantity $L_{K,M}$  for a given specification of the
learning step that is useful to study the complexity of the global algorithm.

\begin{Lemma} \label{le gamma spec LpipsiM}
Let Assumption \ref{as psi function} hold. For $\gamma>0$, $\rho \in (\frac12,1)$, 
set $\gamma_m := \gamma m^{-\rho}$, $m \ge 1$. If the number of steps $M$ in the stochastic gradient descent algorithm satisfies 
\textcolor{black}{ 
\begin{align}\label{eq constraint alpha-beta-M 1}
\gamma\frac{{\alpha}_{ K}}{\beta_{K}}M^{1-\rho} \ge \frac{\sqrt{2}}2\;,  
\end{align}
}
\noindent then, there exists some positive constant $C:=C(\rho, \gamma)$ such that

\begin{align}\label{eq bound LpsiM}
L_{K,M} \le C \left(e^{- 2\sqrt{2}\ln(2)\gamma\frac{\alpha_K }{\beta_{K}}M^{1-\rho}} 
+  \frac{ \beta_{K}}{\alpha_K M^\rho }
\right).
\end{align}

\end{Lemma}

\vspace{2mm}
\begin{Remark}\label{re condition beta}
\begin{enumerate}
\item In practice, $L_{K,M}$ should go to zero with respect to the optimal parameters. Thus, we must have that $\frac{ \beta_{K}}{\alpha_K M^\rho}$ goes to zero  and
$\frac{\alpha_K}{\beta_{K}} M^{1-\rho} $ goes to infinity at the same time as $M$ goes to infinity. This will be carefully discussed in \Cref{subse Com analysis}. With these constraints, we will naturally have that \eqref{eq constraint alpha-beta-M 1} is satisfied.
\item A careful analysis of the proof below shows that $\lim_{\rho \rightarrow  0.5^+} C(\rho, \gamma) = +\infty$, which comes from the dependence of $C(\rho, \gamma)$ with respect to $\sum_{1}^{+\infty} \frac1{m^{2 \rho}}$. However, one must bear in mind that $\rho$ is a fixed (but optimised) parameter.
\end{enumerate}
\end{Remark}

\vspace{2mm}
\noindent \textbf{Proof of Lemma \ref{le gamma spec LpipsiM}.}
We have, since $\Gamma_m = \gamma\sum_{q=1}^m \frac1{q^\rho} $, for $m\ge 1$,
\begin{align}
\frac{\gamma}{1-\rho}\left( m^{1-\rho} - 1\right)  \le \Gamma_m
\le   \frac{\gamma}{1-\rho}\left( m^{1-\rho} - 1\right) +  \gamma 
\end{align}

\noindent by standard computations based on comparison between series and integral, leading to
\begin{align}
\Gamma_m - \Gamma_M \le \frac{\gamma}{1-\rho}\left( m^{1-\rho} - M^{1-\rho}\right) 
+  \gamma .
\end{align}
Recalling \eqref{de LpsiM}, we employ the following decomposition
\begin{align}
 L_{K,M} &=  \varrho_0 \varrho_1  \Big(1+ \mathbb{E}\Big[  |\mathfrak{u}_{0}|^2\Big]  \Big) \sum_{m=1}^{M} \exp\left(- {4} \frac{{\alpha}_{K}}{\beta_{K}} (\Gamma_{M}-\Gamma_{m}) \right) \gamma_m^2 
 \\
 &\le 
  \varrho_0 \varrho_1   \Big(1+ \mathbb{E}\Big[  |\mathfrak{u}_{0}|^2\Big]  \Big)  \exp\left({4}\gamma \frac{{\alpha}_{K}}{\beta_{K}} \right) 
  \big(A_M + B_M +\gamma_M^2 \big)
\label{eq first to conclude}
\end{align}

with
\begin{align*}
A_M &:=  \sum_{m=1}^{\floor{M/2}} \exp\left( {4} \frac{{\alpha}_{K}}{\beta_{K}} \frac{\gamma}{1-\rho}\set{ m^{1-\rho} - M^{1-\rho} } \right)\frac{\gamma^2}{m^{2 \rho}},
\\
B_M &:=  \sum_{m=\floor{M/2}+1}^{M-1} \exp\left( {4} \frac{{\alpha}_{K}}{\beta_{K}} \frac{\gamma}{1-\rho}\set{ m^{1-\rho} - M^{1-\rho} } \right)\frac{\gamma^2}{m^{2 \rho}}.
\end{align*}
For the first term $A_M$, we observe that, for $m \le \floor{M/2} \le M/2$ and $\frac12<\rho<1$,
\begin{align}
 m^{1-\rho} - M^{1-\rho} \le - \frac{\sqrt{2}}{2} \ln(2) (1-\rho)M^{1-\rho}\;.
\end{align}
We then compute
\begin{align}
A_M &\le  \exp\left(- 2\sqrt{2}\ln(2) \gamma \frac{{\alpha}_{K}}{\beta_{K}} M^{1-\rho}  \right) \sum_{m=1}^{\floor{M/2}} \frac{\gamma^2}{m^{2 \rho}} \nonumber
\\
&\le  C_{\rho, \gamma}  \exp\left(- 2\sqrt{2}\ln(2) \gamma \frac{{\alpha}_{K}}{\beta_{K}} M^{1-\rho}  \right). \label{eq upper bound A}
\end{align}
We now study the term $B_M$ which reads
\begin{align*}
B_M = \gamma^2 \exp\left( - {4} \frac{{\alpha}_{K}}{\beta_{K}} \frac{\gamma}{1-\rho} M^{1-\rho}  \right)  \sum_{m=\floor{M/2}+1}^{M-1}\lambda(m)
\end{align*}
where, the map $\lambda$ is defined for $x \ge 1$ by
\begin{align}
\lambda(x) := \exp\left(  {4} \frac{{\alpha}_{K}}{\beta_{K}} \frac{\gamma}{1-\rho} x^{1-\rho}  \right)\frac{1}{x^{2 \rho}}.
\end{align}
We observe that $\lambda$ is increasing on $[\floor{M/2}+1,+\infty)$ when  \eqref{eq constraint alpha-beta-M 1} holds. This leads to
\begin{align*}
\sum_{m=\floor{M/2}+1}^{M-1}\lambda(m) &\le \int_{M/2}^{M} \lambda(x) \ud x\\
& \leq \frac{2^\rho}{M^{\rho}} \int_{M/2}^{M} \exp\left(  {4} \frac{{\alpha}_{K}}{\beta_{K}} \frac{\gamma}{1-\rho} x^{1-\rho}  \right)\frac{1}{x^{\rho}} \ud x
\\
&\le  \frac{ 2^{\rho-  {2} } \beta_{K}}{\gamma{M^\rho \alpha}_{K}}
\exp\left(  {4} \frac{{\alpha}_{K}}{\beta_{K}} \frac{\gamma}{1-\rho} M^{1-\rho}  \right)  
\end{align*}
which in turn yields
$ B_M \le  \frac{ \gamma\beta_{K}}{M^\rho \alpha_{K}}.$

Inserting the previous inequality and estimate \eqref{eq upper bound A} into \eqref{eq first to conclude} concludes the proof,  since $\mathbb{E}\Big[  |\mathfrak{u}_{0}|^2\Big] < \infty$, recall Definition \ref{de picard scheme full}, step 1. and  $\frac{\alpha_K}{\beta_K} \le 1$ recall Assumption \ref{as psi function} and \eqref{eq de beta psi}.
\eproof

\begin{Remark}
Let us importantly point out that if one choses $\gamma_{m}= \gamma /m$, with $\gamma>0$, then from standard comparison between series and integral $\Gamma_{m} - \Gamma_M \leq \gamma (\ln(m/M)+1)$ so that repeating the computations of the proof of Lemma \ref{le gamma spec LpipsiM}, one has to consider the two disjoint cases $\gamma < \frac{\beta_K}{4\alpha_K}$ and $\gamma > \frac{\beta_K}{4\alpha_K}$ in order to provide an upper bound for the quantity of interest $L_{K, M}$. Only the latter allows to obtain the best convergence rate of order $1/M$. However, in practice, the user does not know the exact value of $\frac{\beta_{K}}{4\alpha_{K}}$ so that one will often consider higher values of $\gamma$ than requested which will have the undesirable effect to deteriorate the upper-bound as suggested by the value of $\varrho_1$ in \eqref{de varrho 1}. Moreover, as shown in Section \ref{se con and com analysis}, the value $\frac{\beta_{K}}{4\alpha_{K}}$ actually goes to infinity when the prescribed approximation error $\varepsilon$ goes to zero so that the latter condition becomes more and more stringent. 

\end{Remark}

\noindent We now give an upper bound for the error $\cE_P$ defined by \eqref{eq de error picard}, with respect to all the algorithm's parameters. These parameters will be chosen in the next section taking into account the precise specification of the functional approximation space. 
\begin{Proposition} \label{pr main abstract error}
Suppose that Assumption \ref{ass X0}, Assumption \ref{assumption:coefficient} (i), (ii) and Assumption \ref{as psi function} hold. Assume that there exists a positive constant $\eta$ {(independent of $N$, $M$ and the basis functions $(\psi)$)} such that
\begin{align}\label{eq ass cs control global}
\frac{\kappa_{K}}{h \wedge \alpha_K } L_{K,M}  \le \eta \;.
 \end{align}
\noindent If $T(1+ 2L^2(1+ h)) < 1$ and
$\delta_{h,\eta}:=\frac{ 16L^2 T (1+\eta)}{1-T(1+ 2L^2(1+ h))} < 1$, then for any $\varepsilon>0$ there exists a positive constant $C_\varepsilon$ such that
$$
\cE_{\mathrm{RM}} \leq C_\varepsilon \frac{\kappa_{K}}{h \wedge \alpha_K } L_{K,M} $$
so that, with the notations of Proposition \ref{error of Picard}, it holds
\begin{align}\label{global:control:error}
\cE_P \le \delta_{h, \varepsilon}^P \cE_0 + C_\varepsilon \left(
\frac{\kappa_{K}}{h \wedge \alpha_K } L_{K,M}  + \cE_{\psi} +  \cE_{\pi} \right).
\end{align}
\end{Proposition}

\vspace{2mm}
\begin{Remark} In practice, $\eta$ will be fixed to be a small constant as the term in the left hand side of \eqref{eq ass cs control global} should be asymptotically zero.
\end{Remark}

\vspace{2mm}

\noindent \textbf{Proof of Proposition \ref{pr main abstract error}.}

\noindent \emph{Step 1:} From Proposition \ref{error of Picard}, we see that to obtain \eqref{global:control:error}, it remains to control
\begin{align*}
\cE_{\mathrm{RM}} =\max_{1 \le p \le P} \esp{\vvvert \Phi_M(\mathfrak{u}^{p-1}_M) - \Phi(\mathfrak{u}^{p-1}_M)\vvvert^2}.
\end{align*}
Using Lemma \ref{le tool robbins-monro}, we have
\begin{align}\label{eq  robbins-monro here}
  \esp{\vvvert \Phi_M(\mathfrak{u}^{p-1}_M) - \Phi(\mathfrak{u}^{p-1}_M)\vvvert^2} 
  \le
\kappa_{K} L_{K,M}(1
+  d N )
+ \frac{\kappa_{K}}{\alpha_{K}} L_{K,M} \esp{\vvvert \Phi(\mathfrak{u}^{p-1}_M) \vvvert^2}.
\end{align}
To conclude the proof, we will in the next step, provide an  upper bound for the term $\esp{\vvvert \Phi(\mathfrak{u}^{p-1}_M) \vvvert^2}$, uniformly with respect to $p$.
\\

\noindent \emph{Step 2:} We denote by $C_\varepsilon$ a constant that may change from line to line along with $\varepsilon$.
From Young's inequality and Lemma \ref{le control f(y,z)}, we obtain
\begin{align*}
 \esp{\vvvert \Phi(\mathfrak{u}^p_{M}) \vvvert^2}
 &\le
 (1+\varepsilon) \esp{\vvvert \Phi(\mathfrak{u}^p_{M}) - \bar{\mathfrak{u}}\vvvert^2} + (1+\frac{1}{\varepsilon}) \vvvert \bar{\mathfrak{u}}\vvvert^2
 \\
 &\le \delta_{h}(1+\varepsilon)  \esp{\vvvert \mathfrak{u}^p_{M} -  \bar{\mathfrak{u}} \vvvert^2} + C_\varepsilon(\cE_\pi + \cE_{\psi}+\vvvert \bar{\mathfrak{u}}\vvvert^2)
 \\
  &\le \delta_{h}(1+\varepsilon) \esp{\vvvert \mathfrak{u}^p_{M}\vvvert^2} 
  + C_\varepsilon(\cE_\pi + \cE_{\psi}+\vvvert \bar{\mathfrak{u}}\vvvert^2)
\end{align*}
up to a modification of $\varepsilon$. From the previous inequality, we readily obtain $ \esp{\vvvert \Phi(\mathfrak{u}^0_{M}) \vvvert^2} \leq C_\varepsilon$. Now, if $p$ is a positive integer, noting again that 
$$
\esp{\vvvert \mathfrak{u}^p_{M}\vvvert^2}
\le 2\esp{\vvvert\Phi_M(\mathfrak{u}^{p-1}_{M}) - \Phi(\mathfrak{u}^{p-1}_{M})\vvvert^2}
+ 2 \esp{\vvvert \Phi(\mathfrak{u}^{p-1}_{M})\vvvert^2}
 $$
 we obtain
 \begin{align*}
 \esp{\vvvert \Phi(\mathfrak{u}^p_{M}) \vvvert^2}
 \le
  2\delta_{h}(1+\varepsilon)\esp{\vvvert  \Phi(\mathfrak{u}^{p-1}_{M}) \vvvert^2} &+ 2  \delta_{h}(1+\varepsilon) \esp{\vvvert\Phi_{M}(\mathfrak{u}^{p-1}_{M}) - \Phi(\mathfrak{u}^{p-1}_{M})\vvvert^2} 
  \\
  &+C_\varepsilon \bigg( \cE_\pi + \cE_{\psi}+\vvvert \bar{\mathfrak{u}}\vvvert^2\bigg)
 \end{align*}
\noindent so that using \eqref{eq  robbins-monro here}, we get
 \begin{align*}
 \esp{\vvvert \Phi(\mathfrak{u}^p_{M}) \vvvert^2}
 \le&
  2\delta_{h}(1+\varepsilon)\left(1+\frac{\kappa_{K}}{\alpha_{K}} L_{K,M}\right) \esp{\vvvert  \Phi(\mathfrak{u}^{p-1}_{M}) \vvvert^2}
  \\
&+ 2  \delta_{h}(1+\varepsilon)  \kappa_{K} L_{K,M}\left(1+ d N  \right)
  +C_\varepsilon \bigg( \cE_\pi + \cE_{\psi}+\vvvert \bar{\mathfrak{u}}\vvvert^2\bigg)\;.
 \end{align*}
 From \eqref{eq ass cs control global} and the fact that $\delta_{h, \eta}$, we can set $\varepsilon$ such that $2\delta_{h}(1+\varepsilon)(1+\frac{\kappa_{K}}{\alpha_{K}} L_{K,M}) \le \delta_{h,\eta} (1+\varepsilon)<1$ so that from the above inequality, by induction, for any positive integer $p$, we get
 \begin{align*}
 \esp{\vvvert \Phi(\mathfrak{u}^p_{M}) \vvvert^2} \le   (\delta_{h,\eta}(1+\varepsilon))^{p}\esp{\vvvert \Phi(\mathfrak{u}^0_{M}) \vvvert^2}+ C_\varepsilon(\cE_\pi + \cE_{\psi}+\vvvert \bar{\mathfrak{u}}\vvvert^2 + \delta_h) \le C_\varepsilon\;,
 \end{align*}
 which concludes the proof.
\eproof


\vspace{2mm}
\noindent We now have all the ingredients to give the proof of the main result announced in Section \ref{subse picard algo} on the upper bound for the global convergence error at the initial time.

\paragraph{Proof of Theorem \ref{th main conv result}} From the very definition \eqref{eq de error algorithm} of the global error, we deduce
\begin{align}\label{eq start decompose}
\cE_{\mathrm{MSE}} \le 2\esp{|u(0,\cX_0) - \bar{u}(\cX_0)|^2 + |Y^{\bar{\mathfrak{u}}}_0 - Y^{\mathfrak{u}^P_M}_0|^2}
\end{align}
\noindent where we used the notations introduced in \eqref{eq de Y starting point}, \eqref{y:reference:solution} and \eqref{the:reference:solution}. A fortiori, we have
\begin{align} \label{eq interm step decompose}
\esp{|u(0,\cX_0) - \bar{u}_0(\cX_0)|^2} \le \cE_\psi \quad \text{ and } \quad
\esp{|Y^{\bar{\mathfrak{u}}}_0 - Y^{\mathfrak{u}^P_M}_0|^2} \le \cE_P
\end{align}
recalling \eqref{eq error approx space}, \eqref{eq rewritte norms} and \eqref{eq de error picard}. Combining \eqref{eq interm step decompose} and \eqref{eq start decompose} yields
\begin{align}
\cE_{\mathrm{MSE}}  \le 2 \left( \cE_\psi  + \cE_P \right)\,.
\end{align}
We eventually conclude the proof by invoking Proposition \ref{pr main abstract error} together with Lemma \ref{le gamma spec LpipsiM}. \eproof

\subsection{Convergence  and complexity analysis for sparse grid approximations}
\label{se con and com analysis}

For this part, we work in the setting of Section \ref{subse picard periodic}. Our goal is to prove the theoretical upper-bound on the algorithm's complexity stated in Theorem \ref{thm full complexity}. 


\vspace{2mm}
\noindent We first state the following useful estimate.
\begin{Lemma}\label{le main estimate for frame}
Suppose that Assumption \ref{ass periodic} is satisfied. Let $\phi:\R^d \rightarrow \R$ be a non-negative measurable function whose support is included in $\cO$ and $\widecheck{\phi}$ be its $1$-\emph{periodisation} defined by \eqref{eq de one periodisation}. Then, it holds
\begin{align*}
\mathfrak{C}^{-1}\esp{\phi(U)} \le \esp{\widecheck{\phi}(X_{t_n})} = \esp{\phi(\widehat{X}_{t_n})} \le \mathfrak{C} \esp{\phi(U)}
\end{align*}
where $U$ has law $\cU((0,1)^d)$ and $\mathfrak{C}$ is given in Lemma \ref{Aronson:two:sided:bounds}.
\end{Lemma}

\proof
We denote by $x' \mapsto p_X(t_n, x')$ the density function of $X_{t_n}$ given by the Euler-Maruyama scheme taken at time $t_n$ and starting from $\cX_0$ with  law $\cU((0,1)^d)$ at time $0$. Note that we have
$p_X(t_n,x') = \int p^{\pi}(0, t_n, x, x') \1_{(0,1)^d}(x) \ud x $ and using
\eqref{Aronson:bounds:transition:density},
\begin{align*}
\mathfrak{C}^{-1} \int p(\mathfrak{c} t_n, x - x')\1_{(0,1)^d}(x) \ud x  \leq p_X(t_n,x') \leq \mathfrak{C} \int p(\mathfrak{c} ^{-1}t_n, x'-x) \1_{(0,1)^d}(x) \ud x \,.
\end{align*}
Then,
\begin{align*}
 \esp{\phi(\widehat{X}_{t_n})} &= \esp{\widecheck{\phi}(X_{t_n})}  \ge \mathfrak{C}^{-1} \int \widecheck{\phi}(x')\int p(\mathfrak{c}  t_n, x - x')\1_{(0,1)^d}(x)\ud x \ud x'
 \end{align*}
 so that introducing the notation $\Xi = \xi + W_{ \mathfrak{c}  t_n}$,
\begin{align}
\int \widecheck{\phi}(x')\int p(\mathfrak{c}  t_n, x - x')\1_{(0,1)^d}(x)\ud x \ud x' &= \esp{\widecheck{\phi}(\Xi)}
= \esp{{\phi}(\widehat{ \Xi})}
 \end{align}
 The proof is then concluded by observing that $\cL(\widehat{\Xi}) = \cL(U)$. The proof of the upper-bound follows from similar arguments.
\eproof


%


\subsubsection{Sparse grid approximation error}

\noindent We now provide some upper-bound estimates for the sparse grid approximation error.


\begin{Theorem} \label{th main sg error}
Under Assumption \ref{ass periodic}, there exists a positive constant $C:=C(T, b, \sigma, d, \lambda_0)$ such that 
\begin{align}
\cE_{\psi} \le C 2^{-4\ell} \ell^{d-1}.
\end{align}
\textcolor{black}{
To obtain an error of order $\varepsilon^2$ for the quantity $\cE_{\psi}$ one may thus set $\ell_\varepsilon = \log_2(\varepsilon^{-\frac12}|\log_2(\varepsilon)|^{\frac{d-1}{ 4 }})$ so that, for each $n=1, \cdots, N-1$, the number of basis functions required satisfies
\begin{align*}
K_\varepsilon = \varepsilon^{-\frac12}|\log_2(\varepsilon)|^{  \frac{5(d-1)}4 }.
\end{align*}
}
\end{Theorem}
\proof 
 For $1 \le i \le d$, $n=0, \cdots, N-1$, setting $v_i(t_n, x) = (\sigma^{\top}\nabla_x u)_i(t_n, x)\1_{\set{x \in \cO}}$, $x \in \mathbb{R}^d$, from \eqref{bound:multi:indices:derivatives:pde:solution} we have
 \begin{align}\label{eq bound hmix}
 \|u\|_{H^{k}_{mix}(\mathcal{O})} + \max_{1 \le i \le d}\|v_i\|_{H^{k}_{mix}(\mathcal{O})} \le C.
 \end{align}
 Moreover, from \eqref{identity:periodic:function:euler:scheme}, we obtain
 \begin{align*}
 (\sigma^\top\nabla_x u)_i(t_n,X_{t_n}) = v_i(t_n, \widehat{X}_{t_n})\,,
 \end{align*}
 thus
\begin{align*}
\esp{|(\sigma^\top\nabla_x u)(t_n,X_{t_n})-Z^{\mathfrak{u}}_{t_n}|^2} & =  \sum_{i=1}^d  \esp{\Big|v_i(t_n, \widehat{X}_{t_n}) - \sum_{k=1}^{K}\mathfrak{z}^{n, k}_i \psi^{k}_{n}(\widehat{X}_{t_n})\Big|^2} 
\\
&\le \mathfrak{C} \sum_{i=1}^d  \esp{\Big|v_i(t_n, U) - \sum_{k=1}^{K}\mathfrak{z}^{n, k}_i \psi^{k}_{n}(U)\Big|^2} 
\end{align*}
\noindent with $U\sim \mathcal{U}((0,1)^d)$ and where we use the upper-estimate given in Lemma \ref{le main estimate for frame} to obtain the last inequality.
We also recall that
\begin{align}
\esp{|u(0,\cX_0)-Y^{\mathfrak{u}}_0|^2} = \esp{|u(0,U) - \sum_{k=1}^{K}\mathfrak{y}^k \psi^k_y(U)|^2} \,.
\end{align}

\noindent From \eqref{eq error approx space brut} and the previous estimates, we thus deduce
\begin{align}
\cE_{\psi}  &\le C\left(
\inf_{\xi \in \mathscr{V}^y} \|\xi-u(0, .)\|^2_{L^2(\cO)} 
+
 \max_{0 \le n \le N-1} 
\sum_{i=1}^d\inf_{\xi \in \mathscr{V}^z_n} \|\xi-v_{i}(t_n, .)\|^2_{L^2(\cO)} 
\right) \,,
\nonumber \\
 & \leq C 2^{-4\ell} \ell^{d-1}, \label{eq control A1}
\end{align}
where for the last inequality we used \eqref{eq control sparse basic} and \eqref{eq bound hmix}.
\\
Finally, setting $\ell_\varepsilon = \log_2(\varepsilon^{-\frac12}|\log_2(\varepsilon)|^{\frac{d-1}{ \textcolor{black}{4} }})$ yields $\cE_\psi = O(\varepsilon^2)$ as $\varepsilon \downarrow 0$
and from \eqref{eq dim approx space} (see also Remark \ref{re indexation}) we deduce that in this case
$K_\varepsilon = \varepsilon^{-\frac12}|\log_2(\varepsilon)|^{  \frac{5(d-1)}4 }$.
\eproof

\subsubsection{Norm equivalence constants}

We now provide some estimates for the value of $\alpha_{K}$ and $\kappa_{K}$ appearing in Assumption \ref{as psi function}.
\begin{Proposition} \label{pr minimal bound alpha} Suppose that Assumption \ref{ass periodic} holds. 
In the setting of Section \ref{subse picard periodic}, there exists a constant $\mathfrak{k} \leq 1$ such that
\begin{align*}
\mathfrak{k} =  \alpha_{K} = \frac1{\kappa_{K}} \;.
\end{align*}
\end{Proposition}
\proof 
For any $\mathfrak{u} = (\mathfrak{y}, \mathfrak{z}) \in \R^{K^y} \times \R^{d \bar{K}^z}$, any $0 \le n \le N-1$ and any $l \in \left\{1, \cdots, d\right\}$, from \eqref{eq de Z} we have
\begin{align*}
\esp{|(Z^{\mathfrak{u}}_{t_n})^l|^2}
= \esp{\Big| \sum_{k = 1}^{K} \mathfrak{z}^{n, k}_l  \psi^{k}_n(\widehat{X}_{t_n})\Big|^2}\,,
\end{align*}
Using Lemma \ref{le main estimate for frame}, we obtain
\begin{align}\label{eq encadrement 1}
 \mathfrak{C}^{-1}\esp{\Big| \sum_{k = 1}^{K} \mathfrak{z}^{n, k}_l  \psi^{k}_n(U)\Big|^2}
 \le
 \esp{|(Z^{\mathfrak{u}}_{t_n})^l |^2}
 \le
  \mathfrak{C}\esp{\Big| \sum_{k = 1}^{K} \mathfrak{z}^{n, k}_l  \psi^{k}_n(U)\Big|^2}\;.
\end{align}
%
%
Note that in our setting, $\psi_{n}^{k} = \chi^k$, so it holds
\begin{align}
\int_{\cO}\Big|\sum_{k = 1}^{K} \mathfrak{z}^{n, k}_l  \psi_{n}^{k}(z)\Big|^2
 \ud z &= \int_{\cO}\Big|\sum_{k = 1}^{K} \mathfrak{z}^{n, k}_l  \chi^{k}(x)\Big|^2
\ud x\;. \label{eq min z 2}
\end{align}
Since the basis functions $(\chi^k)_{1 \le k \le K}$ forms a Riesz basis \cite{griebel1995tensor}, there exists a constant $\underline{c}\geq 1$ such that it holds
\begin{align}
\underline{c}^{-1} \sum_{k = 1}^{K} |\mathfrak{z}^{n, k}_l |^2  \le 
\int_{[0, 1]^d}\Big|\sum_{k = 1}^{K} \mathfrak{z}^{n, k}_l  \chi^{k}(x)\Big|^2 \ud x
\le \underline{c}
\sum_{k = 1}^{K} |\mathfrak{z}^{n, k}_l |^2 
\label{eq min z 3}
\end{align}
Combining \eqref{eq encadrement 1}-\eqref{eq min z 2}-\eqref{eq min z 3} and taking into account all the component of $Z^{\mathfrak{u}}_{t_n}$, we compute
\begin{align}\label{eq conclusion for z}
\underline{c}^{-1} \mathfrak{C}^{-1} d \sum_{k=1}^{K}  |\mathfrak{z}^{n, k}_l |^2  \le 
\esp{|(Z^{\mathfrak{u}}_{t_n}) |^2}
\le \underline{c}  \mathfrak{C} d \sum_{k=1}^{K} |\mathfrak{z}^{n, k}_l |^2 \;.
\end{align}
We also observe that
\begin{align}
\esp{|Y^{\mathfrak{u}}_0|^2}
= \esp{\Big| \sum_{k = 1}^{K} \mathfrak{y}^{k}  \chi^{k}(\cX_0)\Big|^2}\,.
\end{align}
Since $\cX_0 \sim \cU((0,1)^d)$, we similarly deduce that
\begin{align}
\underline{c}^{-1} \sum_{k = 1}^{K} |\mathfrak{y}^{n, k}_l |^2  \le 
\esp{|Y^{\mathfrak{u}}_0|^2}
\le \underline{c}
\sum_{k = 1}^{K} |\mathfrak{y}^{n, k}_l |^2. 
\label{eq conclusion y}
\end{align}
The proof is concluded by combining \eqref{eq conclusion for z} and \eqref{eq conclusion y} with \eqref{eq rewritte norms} and setting $\mathfrak{k} = \underline{c}^{-1} \mathfrak{C}^{-1}d^{-1}$.
\eproof
%
%


\subsubsection{Complexity analysis}\label{subse Com analysis}

\begin{Lemma}\label{le bound beta}
Under Assumption \ref{ass periodic}, one can set
\begin{align}
\beta_{K} = C_d(1+  2^\ell \ell^{d - 1})
\end{align}
\noindent for some positive constant $C_d$ which depends on the PDE dimension $d$.
\end{Lemma}
\proof

\noindent \emph{Step 1: } For any $n \in \set{0,\dots,N-1}$, any $\ell \in \set{1,\dots,d}$, one has
\begin{align}
\esp{|\tilde{\omega}^{n,\cdot}_\ell |^4}  &= \esp{\left(\sum_{k=1}^{K} |\tilde{\omega}^{n,k}_\ell|^2\right)^2} \,.
\end{align}
Using Jensen's inequality, we obtain
\begin{align}
\esp{|\tilde{\omega}^{n,\cdot}_\ell |^4}  &\le K \esp{\sum_{k=1}^{K} |\tilde{\omega}^{n,k}_\ell|^4}  \le 3 {d^2}K  \sum_{k=1}^{K} \esp{|\psi_{n}^{k}(\widehat{X}_{t_n})|^4} \;. \label{eq first majo}
\end{align}
From now on, we use the indexation related to the sparse grid description introduced in Remark \ref{re indexation}, namely, we write
\begin{align}\label{eq start majo}
\sum_{k=1}^{K}  \esp{|\psi_{n}^{k}(\widehat{X}_{t_n})|^4} &= \sum_{(\mathbf{l},\mathbf{i}) \in \cC}\esp{|\psi_{n}^{(\mathbf{l},\mathbf{i})}(\widehat{X}_{t_n})|^4 }
\\
&= \sum_{(\mathbf{l},\mathbf{i}) \in \cC}\esp{|\chi^{(\mathbf{l},\mathbf{i})}(\widehat{X}_{t_n})|^4 }
\end{align}
in our setting.

\noindent \emph{Step 2: }
 Using Lemma \ref{le main estimate for frame}, we observe
\begin{align}\label{eq control beta interm 1}
\esp{|\chi^{(\mathbf{l},\mathbf{i})}(\widehat{X}_{t_n})|^4 }  &\le \mathfrak{C} 
\esp{|\chi^{(\mathbf{l},\mathbf{i})}(U)|^4}\;.
\end{align}


\noindent Moreover, we compute
 \begin{align*}
 \int |\phi^{(l,i)}(x)|^4 \ud x = \int |\phi(2^l x - i)|^4 \ud x \le C 2^{-l}
 \end{align*}
 
 \noindent and from the definition of $\chi^{(l_j,i_j)}$, we deduce
 \begin{align*}
 \int |\chi^{({l}_j,{i}_j)}(x)|^4 \ud x \le C 2^{l} \;.
 \end{align*}

\noindent Combining the previous inequality with  \eqref{eq control beta interm 1} leads to
 \begin{align*}
\esp{|\chi^{(\mathbf{l},\mathbf{i})}(\widehat{X}_{t_n})|^4 }  \leq {C} 2^{|\mathbf{l}|_1}\;
 \end{align*}
 
 \noindent and inserting the previous estimate in \eqref{eq start majo} yields
 \begin{align} \label{eq first step 3}
 \sum_{k=1}^{K}  \esp{|\psi_{n}^{k}(\widehat{X}_{t_n})|^4} &\le {C} \sum_{(\mathbf{l},\mathbf{i}) \in \cC} 2^{|\mathbf{l}|_1}.
 \end{align}
\noindent \emph{Step 3: } We now quantify the term appearing in the right-hand side of \eqref{eq first step 3}, namely
 \begin{align}
Q :=  \sum_{(\mathbf{l},\mathbf{i}) \in \cC} 2^{|\mathbf{l}|_1}
 = 1 + \sum_{k=1}^\ell \sum_{\mathbf{l} \in \mathbb{N}^d} 2^{|\mathbf{l}|_1}\1_{\set{\zeta_d(\mathbf{l})=k}}\;.
 \end{align}
 We denote by $\| \mathbf{l} \|_0 = |\set{j|l_j = 0}|$. For $\mathbf{l} \neq \mathbf{0}$,
 recall that $\zeta_d(\mathbf{l}) = |\mathbf{l}|_1 + \|\mathbf{l}\|_0 -(d-1)$ (from the definition of $\zeta_d$). Thus, 
 \begin{align*}
 Q &= 1 + \sum_{k=1}^\ell 
 \sum_{q=0}^{d-1}
 \sum_{\mathbf{l} \in (\mathbb{N}_{>0})^{d-q}} 
 2^{|\mathbf{l}|_1}\1_{\set{|\mathbf{l}|_1=k+ d-1 -q  \text{ and } \|\mathbf{l}\|_0 = q } }
 \\
 &= 1 + \sum_{k=1}^\ell \sum_{q=0}^{d-1} 2^{k+ d-1-q} \mathrm{C}_{k + d-q-2}^{d-q-1}
 \mathrm{C}_d^q
 \end{align*}
 recall that $|\set{ \mathbf{l} \in (\mathbb{N}_{>0})^{d-q} \,| \,|\mathbf{l}|_1=k+ d-1 -q}| =\mathrm{C}_{k + d-q-2}^{d-q-1} $.
 Introducing $\theta = d-1-q$, we get
 \begin{align}
 Q = 1 + \sum_{k=1}^\ell d2^k + \sum_{\theta = 1}^{d-1} 2^{1+\theta} \mathrm{C}^{d-1 - \theta}_{d} \sum_{k=1}^\ell 2^{k-1} \mathrm{C}_{k-1+\theta}^{\theta}
 \end{align}
 From Lemma 3.6 in \cite{bungartz2004sparse}, we know that $\sum_{k=1}^\ell 2^{k-1} \mathrm{C}_{k-1+\theta}^{\theta} = 2^\ell(\frac{\ell^\theta}{\theta!} + O_d(\ell^{\theta-1}))$. We thus obtain
 \begin{align*}
 Q = 2^\ell\left(\frac{2^d}{(d-1)!}\ell^{d-1} + O_d(\ell^{d-2})\right)
\end{align*}
\noindent which combined with \eqref{eq first step 3} yields
\begin{align} \label{eq last step 3}
 \sum_{k=1}^{K}  \esp{|\psi_{n}^{k}(\widehat{X}_{t_n})|^4} \le C_d 2^\ell \ell^{d-1}.
 \end{align}
Combining the previous inequality with \eqref{eq first majo}, we obtain
\begin{align}
\esp{|\tilde{\omega}^{n,\cdot}_\ell|^4} \le C_d 2^{2\ell} \ell^{2d-2}\;.
\end{align}
Using similar arguments, as the basis function are chosen to be the same in our setting, we also have
\begin{align*}
\esp{|\theta|^4} \le C_d 2^{2\ell} \ell^{2d-2}\;.
\end{align*}
The proof is then concluded recalling the definition of $\beta_{K}$ in \eqref{eq de beta psi}.
\eproof


\vspace{2mm}
\noindent We now turn to the analysis of the convergence and complexity of the full Picard algorithm. The following corollary is a preparatory result and expresses the main convergence results in terms of the parameters $P$, $M$, $\ell$ and $h$.

\begin{Corollary} \label{co conv results new}
Suppose that Assumption \ref{ass periodic} as well as \eqref{eq constraint alpha-beta-M 1} and \eqref{eq ass cs control global}  hold. Set $\gamma_m = \gamma/ m^{\rho}$, for some $\rho \in (1/2, 1)$ and $\gamma > 0$. If $T(1+ 2L^2(1+ h)) < 1$ and $\delta_{h,\eta}=\frac{16L^2 T (1+2\eta)}{1-T(1+ 2L^2(1+ h))} < 1$, then, with the notations of Proposition \ref{error of Picard}, for any $\varepsilon>0$ such that $\delta_{h, \varepsilon}:= \delta_{h, \eta } (1+\varepsilon)<1$ there exist constants $C_\varepsilon:=C(\varepsilon, T, b, \sigma, d, \gamma,\rho)\geq1$, $c := c(T, b, \sigma, d, \gamma) > 0$ such that it holds
\begin{align} \label{global:mean:squared:error:before:complexity}
\cE_P \le \delta_{h, \varepsilon}^P \cE_0 + C_{\varepsilon} 
\left(
Ne^{- c \frac{M^{1-\rho}}{1 + 2^\ell\ell^{d-1}}} 
+  \frac{ 1 + 2^\ell\ell^{d-1}}{h M^\rho } + 2^{-4\ell}\ell^{d-1} + h
\right) \;.
\end{align}

\end{Corollary}
\proof Combining Proposition \ref{pr main abstract error} with Theorem \ref{th main sg error} and \eqref{eq control disc error}, we obtain
\begin{align}
\cE_P \le \delta_{h, \varepsilon}^P \cE_0 + C_\varepsilon \left(
\frac{\kappa_{K}}{h \wedge \alpha_{K} } L_{K,M}  + 2^{-4\ell}\ell^{d-1} + h \right).
\end{align}
From Proposition \ref{pr minimal bound alpha}, we have that $\frac{\kappa_{K}}{h \wedge \alpha_{K} } \le \frac{C}{h}$, which combined with Lemma \ref{le gamma spec LpipsiM} gives
\begin{align}
\frac{\kappa_{K}}{h \wedge \alpha_{K} } L_{K,M} &\le \frac{ C_{\rho,\gamma} }{h} \left(e^{- \gamma 2 \sqrt{2} \ln(2)\frac{\mathfrak{k} }{\beta_{K}}M^{1-\rho}} 
+  \frac{ \beta_{K}}{\mathfrak{k} M^\rho }
\right)
\\
&\le \frac{ C_{\rho,\gamma} }{h} \left(e^{- c \frac{M^{1-\rho}}{1 + 2^\ell\ell^{d-1}}} 
+  \frac{ 1 + 2^\ell\ell^{d-1}}{M^\rho }
\right)
\end{align}
for some positive constant $c$, where we used Lemma \ref{le bound beta} for the last inequality.
\eproof

We are now ready to establish the complexity of the full Picard algorithm.

\paragraph{Proof of Theorem \ref{thm full complexity}}
\emph{Step 1: Setting the parameters $P$, $N$, $M$, $\ell$ and $\rho$.}\\ We will chose the parameters $P$, $N$, $M$, $\ell$ and $\rho$ in order to achieve a global error $\cE_P$ of order $\varepsilon^2$, as this error controls $\cE_{\mathrm{MSE}}$. We first set $P = 2 |\log_\delta(\varepsilon)|$ and \textcolor{black}{$N_\varepsilon= \ceil*{T\varepsilon^{-2}} $ so that $h_\varepsilon=T/N_\varepsilon\leq \varepsilon^2$}. From Theorem \ref{th main sg error}, we also know that setting $\ell_\varepsilon = \log_2(\varepsilon^{-\frac12}|\log_2(\varepsilon)|^{\frac{d-1}{4 } })$, we obtain $\cE_\psi= O_d(\varepsilon^2)$ and $K_\varepsilon = O_d(\varepsilon^{-\frac12}|\log_2(\varepsilon)|^{\textcolor{black}{\frac{5(d-1)}4 }})$.

We now set $M$ such that the term $\frac{K_\varepsilon}{M^\rho h_\varepsilon}$ is of order $\varepsilon^2$, which leads to
\begin{align*}
M_\varepsilon = O_d(\varepsilon^{-\frac9{2 \rho}}|\log_2(\varepsilon)|^{ \textcolor{black}{ \frac{5(d-1)}{4\rho }} })\,.
\end{align*}
For $\iota > 1$, we set $\bar{\rho} = \frac9{10 \iota}$ with the constraint $\bar{\rho} > \frac12$ and we verify that
\begin{align*}
\frac{M_\epsilon^{1-\bar{\rho}}}{K_\varepsilon} \ge c \varepsilon^{5(1-\iota)}|\log_2(\varepsilon)|^{\textcolor{black}{\frac54}(d-1)(\frac1{\bar{\rho}} - 2)}
\end{align*}
for some constant $c>0$. This leads to
\begin{align*}
e^{- c \frac{M_\varepsilon^{1-\bar{\rho}}}{1 + 2^{\ell_\varepsilon} \ell_\varepsilon^{d-1}}}  = o(\varepsilon^4) 
\end{align*}
and we also have that \eqref{eq constraint alpha-beta-M 1} is satisfied.\\
\emph{Step 2: Computing the complexity $\cC_\varepsilon$.} Recalling Remark \ref{re num complexity direct algo}, we see that the overall complexity $\cC_\varepsilon$ to reach the prescribed approximation accuracy $\varepsilon^2$ satisfies 
\begin{align*}
\cC_\varepsilon &= P_\varepsilon N_\varepsilon K_\varepsilon M_\varepsilon
\\
&= O_d \left(
|\log_{2}(\varepsilon)| \varepsilon^{-2} \varepsilon^{-\frac12}|\log_2(\varepsilon)|^{\textcolor{black}{\frac{5(d-1)}4}}
\varepsilon^{-5\iota }|\log_2(\varepsilon)|^{ \frac{25(d-1)\iota}{18}}
\right)
\\
&=O_d(\varepsilon^{-\frac52(1+2\iota )}|\log_{2}(\varepsilon)|^{1+ \frac{45+50\iota}{36}(d - 1)} )\;,
\end{align*}
which concludes the proof.

\eproof
\section{Appendix}


\label{se algo param}

We gather below all the  parameters values used in the various algorithms, examples and basis functions settings.
In particular, we recall that the domain specification is given in \eqref{eq de BM domain} and \eqref{eq de GBM domain} The learning rates are given by \eqref{eq de learning rate}.
Denoting
\[ \Gamma(\lambda) := ( \alpha(\lambda ), \beta_1(\lambda), \beta_0(\lambda), m_0(\lambda)) \in \mathbb{R}^4,  \]
we set the  parameters in all the approximating space  $\mathscr{V}^z_n , 1\le n\le N-1$ to be the same: namely $\Gamma(\mathfrak{z}^{n, \cdot}) = \Gamma(\mathfrak{z}^{1, \cdot })$ for all $n\ge 2$. Thus, in the table below, the parameters of the learning rates are simply denoted by: 
\[ \Gamma := \{ \Gamma(\mathfrak{y}), \Gamma(\mathfrak{z}^{0, \cdot} ), \Gamma(\mathfrak{z}^{1, \cdot }) \} \in \mathbb{R}^{3\times 4}.  \]

\begin{table}[H]
	\centering
	\resizebox{\textwidth}{30mm}{
	\begin{tabular}{|c|c|c|c|c|c|c|c|c|c| } 
		\hline
		Algorithms & Basis functions &  dim & N & M & T &  Initial $\mathfrak{z}^{n, k}$ & p & $\Gamma$ & $\cE_{\mathrm{MSE}} $ \\
		\hline
		\multirow{8}*{\tabincell{c}{ \emph{Picard Algorithm} }} & \multirow{8}*{Pre-wavelets} & \multirow{8}*{3} & \multirow{8}*{10}  & \multirow{8}*{100000}   &  \multirow{8}*{0.3}  &  \multirow{8}*{0} &  1 &  \tabincell{c}{ (0.6, 0, 3, 15000),  \\ (0.6, 0, 20, 10000) }  & 0.0286 \\
		\cline{8 - 10}
		&  &  &   &   &    &   &  2 &  \tabincell{c}{ (0.6, 0, 2, 8000),  \\ (0.6, 0, 10, 8000 } & 0.0247  \\
		\cline{8 - 10}
		&  &  &   &   &    &   &  3 &  \tabincell{c}{ (0.6, 0, 1, 8000),  \\ (0.6, 0, 5, 8000 } & 0.0219 \\
		\cline{8 - 10}
		&  &  &   &   &    &   &  4 &  \tabincell{c}{ (0.6, 0, 0.5, 8000),  \\ (0.6, 0, 3, 8000) } & 0.0207 \\
		\cline{8 - 10}
		&  &  &   &   &    &   &  5 &  \tabincell{c}{ (0.6, 0, 0.5, 8000),  \\ (0.6, 0, 3, 8000) } & 0.0201 \\
		\hline
	\end{tabular}
	} 
	\caption{Parameters for the periodic example} \label{Parameters periodic}
\end{table}

\begin{table}[H]
	\centering
	\resizebox{\textwidth}{32mm}{
	\begin{tabular}{|c|c|c|c|c|c|c|c| } 
		\hline
		\multicolumn{1}{|c|}{\multirow{1}*{Examples}}&  a & dim & M &  N &   Initial $y_0$  & \multirow{1}*{ Learning rate } \\
		\hline
		\multirow{1}*{\tabincell{c}{Quadratic}} & 1  & 5  & 2000 &  10   &  0.5 & 0.01    \\
		\hline
		\multirow{1}*{\tabincell{c}{Limits to \\ Picard algorithm}} & -0.4  & 2 & 5000 & 20   &  0.3 &  0.002   \\
		\cline{2 - 7}
		 & -1.5 & 2  & 5000 &  20  &  0  &  0.001  \\
		\hline
		\multirow{7}*{\tabincell{c}{Financial \\ example}} & \multirow{7}*{N.A.}  & 2  &6000 &  20  &  2  &  0.005    \\
		\cline{3 - 7}
		& & 4  & 6000&  20  & 5   &   0.005  \\
		\cline{3 - 7}
		 &   & 5  &  5000&   20  & 5   &   0.005      \\
		\cline{3 - 7}
		& & 10  & 5000  & 20  & 8   &   0.005     \\
		\cline{3 - 7}
		&  & 15  &5000 &  20  & 8   &   0.005     \\
		\cline{3 - 7}
		&  & 20  &5000 &   20  & 8   &   0.005    \\
		\cline{3 - 7}
		&  & 25  & 5000 &  20  & 8   &   0.005   \\
		\hline
	\end{tabular}
	} 
	\caption{Parameters by model for the \emph{deep learning method} with $ layers = 4, batchsize = 64$}
	 \label{Parameters deep learning}
\end{table}

\begin{table}[htp]
	\centering
	\resizebox{\textwidth}{110mm}{
	\begin{tabular}{|c|c|c|c|c|c|c|c|c| } 
		\hline
		\multicolumn{1}{|c|}{\multirow{2}*{Examples}}& \multicolumn{1}{c|}{\multirow{2}*{Basis functions}}&  \multirow{2}*{ dim }& \multirow{2}*{ N }& \multirow{2}*{r} & \multicolumn{3}{c|}{ \multirow{1}*{Initial value}}  & \multirow{2}*{ $\Gamma$ } \\
		 \cline{6-8}
		& &   & & & $\mathfrak{y} $ & $ \mathfrak{z}^{0, \cdot} $ &  $\mathfrak{z}^{n, k}$ & \\
		\hline
		\multirow{4}*{\tabincell{c}{Quadratic \\ model}} & Pre-wavelets & 5 & 10  & 2   & 0.5  &  -0.2 & 0 &  \tabincell{c}{ (1, 0, 1, 100), \\ (0.8, 0, 1, 100),  \\ (0.86, 0.02, 0.05, 100) }  \\
		\cline{2 - 9}
		  & \multirow{1}*{Hat} & 
		  100 & 10  & 3.2   &   5.5  &  0 & 0  &  \tabincell{c}{ (0.9, 0, 1, 100), \\ (0.7, 0, 1, 100),  \\ (1, 0.003, 0.01, 1000) }  \\
		\hline
		\multirow{20}*{\tabincell{c}{Financial \\ example}} & \multirow{3}*{Pre-wavelets} & 2 & 10  & 2   &  2  &  0 &  0&  \tabincell{c}{ (1, 0, 1.5, 1000), \\(1, 0, 20, 1000),     \\ (1, 0.003, 0.01, 1000)}  \\
		\cline{3 - 9}
		& & 4& 10  & 2   &  5  &  0 &  0&  \tabincell{c}{(1, 0, 1, 1000), \\(1, 0, 20, 1000),     \\ (1, 0.001, 0.01, 1000) }  \\
		\cline{2 - 9}
		 & \multirow{13}*{Hat}  & 5 & 10 & 2   &  5  & 0  & 0 &  \tabincell{c}{ (0.95, 0, 0.3, 100), \\ (1, 0, 5, 100), \\ (1, 0.001, 0.01, 100) }  \\
		\cline{3 - 9}
		& & 10 & 10 & 2.5   &   8 &  0 & 0 &  \tabincell{c}{ (0.9, 0, 0.35, 300), \\ (1, 0, 5, 300), \\ (1, 0.001, 0.01, 300) }  \\
		\cline{3 - 9}
		&  & 15& 10  & 2.5   &  8  & 0  & 0 &  \tabincell{c}{ (0.85, 0, 0.3, 500), \\ (1, 0, 5, 500), \\ (1, 0.001, 0.01, 500) }  \\
		\cline{3 - 9}
		&  & 20 & 10 & 2.5   &  8  & 0  & 0 &  \tabincell{c}{(0.8, 0, 0.2, 1000), \\ (1, 0, 5, 1000), \\ (1, 0.001, 0.01, 1000) }  \\
		\cline{3 - 9}
		&  & 25 & 10 & 2.8   &  8  &  0 & 0 &  \tabincell{c}{ (0.7, 0, 0.2, 1500), \\ (1, 0, 5, 1500), \\ (1, 0.001, 0.01, 1500) }  \\
		\hline
		\multirow{13}*{\tabincell{c}{The \\challenging\\ example}} & \multirow{13}*{Hat}  & 1 & 10 & 2   &  0.5  &  0  &  0 &  \tabincell{c}{ (1, 0, 0.5, 100), \\ (1, 0, 3, 100),  \\ (1, 0.1, 0.1, 100) }  \\
		\cline{3 - 9}
		&   & 2  & 20 & 2   &  0.1  & 0  &  0&  \tabincell{c}{ (1, 0, 0.5, 100), \\ (1, 0, 5, 100),  \\ (1, 0.2, 0.5, 100) }  \\
		\cline{3 - 9}
		&  & 5 & 40 & 2   &   0.4 &-1   & 0 &  \tabincell{c}{ (1, 0, 0.5, 300), \\ (0.95, 0, 5, 500),  \\ (1, 0.2, 0.1, 500) }  \\
		\cline{3 - 9}
		&  & 8  & 60& 2.2  &   0.6 & 1  & 0.08 &  \tabincell{c}{ (1, 0, 0.35, 500), \\ (1, 0, 5, 500),  \\ (1, 0.2, 1, 500) }  \\
		\cline{3 - 9}
		& & 10 & 100 & 2.5  &  0.1  & -1  & -0.1 &  \tabincell{c}{  (1, 0, 0.35, 500), \\ (0.95, 0, 5, 500),  \\ (1, 0.2, 1, 500) }  \\
		\hline
	\end{tabular}
	} 
	\caption{Parameters by model for the \emph{direct algorithm}}
	 \label{Parameters direct}
\end{table}

\begin{table}[htp]
	\centering
	\resizebox{\textwidth}{110mm}{
	\begin{tabular}{|c|c|c|c|c|c|c|c|c|c| } 
		\hline
		\multicolumn{1}{|c|}{\multirow{2}*{Examples}}& \multicolumn{1}{c|}{\multirow{2}*{Basis functions}}&  \multirow{2}*{ dim }& \multirow{2}*{ r }&  \multicolumn{3}{c|}{ \multirow{1}*{Initial value}}  & \multirow{2}*{ a } & \multirow{2}*{ p } & \multirow{2}*{ $\Gamma$ } \\
		 \cline{5-7}
		& &   & & $\mathfrak{y} $ & $ \mathfrak{z}^{0, \cdot} $ &  $\mathfrak{z}^{n, k}$  & & & \\
		\hline
		\multirow{10}*{\tabincell{c}{Quadratic \\ model}} & \multirow{8}*{Pre-wavelets} &  \multirow{8}*{5}  &  \multirow{8}*{2}   &  \multirow{8}*{0.5}  &  \multirow{8}*{-0.2} & \multirow{8}*{0} &  \multirow{8}*{1} &  1 & \tabincell{c}{ (1, 0, 0.8, 100), \\ (0.9, 0, 1, 100),  \\ (0.84, 0.02, 0.05, 100) }  \\
		\cline{ 9 -10}
		&  & &&&& & &2 &   \tabincell{c}{ (1, 0, 0.3, 100), \\ (0.9, 0, 0.4, 100),  \\ (0.84, 0.01, 0.02, 100) }  \\
		\cline{9 - 10}
		&  & &&&& &  & $ p \ge 3 $& \tabincell{c}{ (1, 0, 0.2, 100), \\ (0.9, 0, 0.2, 100),  \\ (0.84, 0.005, 0.01, 100) }  \\
		\cline{2 - 10}
		  & \multirow{1}*{Hat} & \multirow{1}*{25}  & \multirow{1}*{2.8 }  & \multirow{1}*{ 2}   & \multirow{1}*{0.1}  & \multirow{1}*{0} &\multirow{1}*{1} &$1\le p\le 3$ & \tabincell{c}{ (0.9, 0, 0.5, 100), \\ (0.8, 0, 0.8, 100),  \\ (1, 0.003, 0.01, 1000) }  \\
		\hline
		\multirow{25}*{\tabincell{c}{Financial \\ example}} & \multirow{12}*{Pre-wavelets} & \multirow{12}*{4}  & \multirow{12}*{2}   & \multirow{12}*{5}   & \multirow{12}*{0.1}  & \multirow{12}*{-0.01} & \multirow{12}*{N.A.} & 1&  \tabincell{c}{ (1, 0, 1, 1000), \\ (1, 0, 20, 1000),  \\ (1, 0.001, 0.01, 1000) }  \\
		\cline{9- 10}
		& &   &    &    &   &  & &2&  \tabincell{c}{(1, 0, 0.3, 1000), \\ (1, 0, 5, 1000),  \\ (1, 0.0005, 0.005, 1000)  }  \\
		\cline{9- 10}
		& &   &    &    &   &  & &3&  \tabincell{c}{ (1, 0, 0.2, 1000), \\ (1, 0, 5, 1000),  \\ (1, 0.0003, 0.003, 1000) }  \\
		\cline{9- 10}
		& &   &    &    &   &  & &4, 5&  \tabincell{c}{ (1, 0, 0.15, 1000), \\ (1, 0, 3, 1000),  \\ (1, 0.0002, 0.002, 1000) }  \\
		\cline{9- 10}
		& &   &    &    &   &  & & $p\ge 6$ &  \tabincell{c}{ (1, 0, 0.1, 1000), \\ (1, 0, 2, 1000),  \\ (1, 0.0001, 0.001, 1000)  }  \\
		\cline{2- 10}
		 & \multirow{10}*{Hat}  &\multirow{10}*{ 20}  & \multirow{10}*{2.5}   &   \multirow{10}*{8} & \multirow{10}*{0} & \multirow{10}*{0} & \multirow{10}*{N.A.} & 1& \tabincell{c}{ (0.9, 0, 0.6, 1000), \\ (1, 0, 5, 1000),  \\ (1, 0.001, 0.01, 1000) }  \\
		 \cline{9 - 10}
		 &&&  &    &  & &  &2 & \tabincell{c}{ (0.9, 0, 0.2, 1000), \\ (1, 0, 4, 1000),  \\ (1, 0.001, 0.01, 1000) }  \\
		 \cline{9 - 10}
		 & &&  &    &  & &  & 3& \tabincell{c}{ (0.9, 0, 0.15, 1000), \\ (1, 0, 3, 1000),  \\ (1, 0.0005, 0.005, 1000) }  \\
		 \cline{9 - 10}
		 &&&   &    &  & &  &  $p\ge 4$& \tabincell{c}{ (0.9, 0, 0.1, 1000), \\ (1, 0, 2, 1000),  \\ (1, 0.0005, 0.005, 1000) }  \\
		\hline
		\multirow{4}*{\tabincell{c}{Limits to \\ Picard \\ algorithm}} & \multirow{4}*{Pre-wavelets}  &  \multirow{4}*{2}  &  \multirow{4}*{2}   & 0.3   & 0 & 0   & -0.4&  $1\le p \le 9$  &   \tabincell{c}{ (1, 0, 1, 300), \\ (1, 0, 1, 300),  \\ (1, 0.6, 0.1, 300) }  \\
		\cline{5 - 10}
		&   & &   &  0  &  0 & 0&  -1.5 & $1\le p \le 9$ &  \tabincell{c}{ (1, 0, 2, 300), \\ (1, 0, 1, 300),  \\ (1, 0.6, 0.1, 300)  }  \\
		\hline
%
	\end{tabular}
	} 
	\caption{Parameters by model for the \emph{Picard algorithm}}
	 \label{Parameters picard 1}
\end{table}

\begin{table}[htp]
	\centering
	\resizebox{\textwidth}{80mm}{
	\begin{tabular}{|c|c|c|c|c|c|c|c|c|c| c| } 
		\hline
		\multicolumn{1}{|c|}{\multirow{2}*{Examples}}& \multicolumn{1}{c|}{\multirow{2}*{Basis functions}}&  \multirow{2}*{ dim }& \multirow{2}*{ r }&  \multicolumn{3}{c|}{ \multirow{1}*{Initial value}}  & \multirow{2}*{ N } & \multirow{2}*{ p } & \multirow{2}*{ $\Gamma$ } \\
		 \cline{5-7}
		& &   & & $\mathfrak{y} $ & $ \mathfrak{z}^{0, \cdot} $ &  $\mathfrak{z}^{n, k}$  & & & \\
		\hline
		\multirow{25}*{\tabincell{c}{The \\challenging \\ example}} & \multirow{25}*{Hat}  & 1  & 2   &  0.5  & 0& 0   & 10 & $1\le p\le 5$ &  \tabincell{c}{(1, 0, $0.5*(0.8)^{p-1}$, 100), \\ (1, 0, 3$*(0.8)^{p-1}$, 100),  \\ (1, 0.1$*(0.8)^{p-1}$, 0.1$*(0.8)^{p-1}$, 100)  }  \\
		\cline{3 - 10}
		&   & 2  & 2   &  0.1 & 0 & 0  & 20 & $1\le p\le 5$ &  \tabincell{c}{  (1, 0, 0.5$*(0.8)^{p-1}$, 100), \\ (1, 0, 5$*(0.8)^{p-1}$, 100),  \\ (1, 0.2$*(0.8)^{p-1}$, 0.5$*(0.8)^{p-1}$, 100)  }  \\
		\cline{3 - 10}
		&  & 5  & 2   &0.4 &  -2 & 0 & 40  & $1\le p\le 5$ &  \tabincell{c}{ (1, 0, 0.5$*(0.8)^{p-1}$, 300), \\ (0.95, 0, 5$*(0.8)^{p-1}$, 500),  \\ (1, 0.2$*(0.8)^{p-1}$, 0.1$*(0.8)^{p-1}$, 500)}  \\
		\cline{3 - 10}
		&  & 8  & 2.2    &0.6 & 1  &  0.08 &  60 & $1\le p\le 5$&  \tabincell{c}{ (1, 0, 0.35$*(0.8)^{p-1}$, 500), \\ (1, 0, 5$*(0.8)^{p-1}$, 500),  \\ (1, 0.2$*(0.8)^{p-1}$, $(0.8)^{p-1}$, 500) }  \\
		\cline{3 - 10}
		& & \multirow{13}*{10}  & \multirow{13}*{2.5}    & \multirow{13}*{0.1} &  \multirow{13}*{-1} &  \multirow{13}*{-0.1} &  \multirow{13}*{100} & 1 & \tabincell{c}{(1, 0, 0.25, 500),     \\ (0.95, 0, 4, 500),     \\ (1, 0.15, 0.5, 500) }  \\
		\cline{9 - 10}
		& & &&&&&   &2  & \tabincell{c}{ (1, 0, 0.2, 500),    \\ (0.95, 0, 3, 500),     \\ (1, 0.12, 0.4, 500) }  \\
		\cline{9 - 10}
		& & &&&&&  &3  & \tabincell{c}{ (1, 0, 0.15, 500),    \\ (0.95, 0, 2, 500),     \\ (1, 0.1, 0.3, 500)  }  \\
		\cline{9 - 10}
		& & &&&&&   &  $p\ge 4$ & \tabincell{c}{  (1, 0, 0.1, 500),     \\ (0.95, 0, 1, 500),    \\  (1, 0.08, 0.25, 500)}  \\
		\cline{9 - 10}
		& & &&&&&  & $p\ge 6$  & \tabincell{c}{ (1, 0, 0.05, 500),    \\ (0.95, 0, 0.5, 500),   \\ (1, 0.05, 0.2, 500)}  \\
		\hline
	\end{tabular}
	} 
	\caption{Parameters by model for the \emph{Picard algorithm}}
	 \label{Parameters picard 2}
\end{table}

\bibliography{Bibliography}
\bibliographystyle{acm}

\end{document}